\documentclass[11pt,reqno]{amsart}

\usepackage{a4wide}

\usepackage{amscd,amsmath,amssymb}
\usepackage[mathscr]{eucal}
\usepackage[all]{xy}
\usepackage{graphics}
\usepackage{mathtools}
\usepackage[english]{babel}

\usepackage[charter]{mathdesign}

\usepackage{hyperref}
\usepackage{xcolor}
\hypersetup{colorlinks=true,urlcolor=black,linkcolor=black,citecolor=black}

\makeatletter \@addtoreset{equation}{section}\makeatother

\DeclareMathAlphabet{\mathpzc}{OT1}{pzc}{m}{it}

\newtheorem{theorem}{Theorem}[section]
\newtheorem{lemma}[theorem]{Lemma}
\newtheorem{proposition}[theorem]{Proposition}
\newtheorem{corollary}[theorem]{Corollary}

\newtheorem{remark}[theorem]{Remark}

\newcommand{\Z}{{\mathbb{Z}}}

\newcommand{\ZZ}{{\Z/2\Z}}
\newcommand{\K}{{\mathbb{K}}}
\newcommand{\C}{{\mathbb{C}}}
\newcommand{\rC}{{\mathcal{C}}}
\newcommand{\ddbar}{{\delta_{\rm bar}}}
\newcommand{\ddcurve}{{\delta_{\rm curv}}}
\newcommand{\ddh}{{\partial_{\rm Hoch}}}
\newcommand{\ddw}{{\partial_{\rm curv}}}
\newcommand{\bbh}{{b_{\rm Hoch}}}
\newcommand{\bbw}{{b_{\rm curv}}}
\newcommand{\ddko}{{\delta_{\rm Kos}}}
\newcommand{\ddk}{{\partial_{\rm Kos}}}
\newcommand{\bbk}{{b_{\rm Kos}}}
\newcommand{\ddkg}{{\partial_{\rm Kos}}{(g)}}
\newcommand{\ddwg}{{\partial_{\rm curv}}{(g)}}
\newcommand{\bbkg}{{b_{\rm Kos}}{(g)}}
\newcommand{\bbwg}{{b_{\rm curv}}{(g)}}
\newcommand{\LL}{{\Upsilon}}
\newcommand{\LLL}{{\Upsilon^\dagger}}

\newcommand{\Jaco}{{\,^t{\rm M}}}
\newcommand{\Jacob}{{\Jac'}}
\newcommand{\Omeg}{{{\,^t\Omega}}}
\newcommand{\Jacoii}{{\rm M}^*}
\newcommand{\Omegii}{{{\Omega^*}}}

\newcommand{\gen}{{{\mathtt{g}}}}
\newcommand{\rB}{{{\mathcal{B}}}}
\newcommand{\rD}{{{\mathcal{D}}}}
\newcommand{\rHH}{{\mathcal{HH}}}
\newcommand{\cHH}{\,^t{\mathsf{HH}}}
\newcommand{\ccHH}{{\mathsf{HH}}}
\newcommand{\cHD}{\,^t{\mathsf{HD}}}

\newcommand{\rK}{{\mathcal{K}}}

\newcommand{\rKo}{{\mathcal{K}\!os}}
\newcommand{\cKo}{{\,^t\mathsf{Kos}}}

\newcommand{\uu}{{\it t}}
\newcommand{\bbb}{(\!({\uu})\!)}
\newcommand{\bb}{[{\uu}^{\pm1}]}
\newcommand{\Hom}{{\mathrm{Hom}}}

\newcommand{\Jac}{{\rm Jac}}

\newcommand{\age}{{\mathrm{age}}}
\newcommand{\zv}{{\widehat{\zeta}}}
\newcommand{\orig}{{\underline{0}}}
\allowdisplaybreaks

\title{On Hochschild invariants of Landau-Ginzburg orbifolds}

\dedicatory{{\sc with an Appendix by  A. Basalaev and D. Shklyarov}}

\author{Dmytro Shklyarov}
\address{Technische Universit\"at Chemnitz,  {\tt dmytro.shklyarov@mathematik.tu-chemnitz.de}}

\address{Ruprecht-Karls-Universit\"at Heidelberg,  {\tt abasalaev@mathi.uni-heidelberg.de }
}

\begin{document}

\begin{abstract} 
We develop an approach to calculating  the cup and cap products on  Hochschild cohomology and homology of curved algebras  associated with polynomials and their finite abelian symmetry groups. For polynomials with isolated critical points, the approach yields a complete description of the products. We also reformulate the result for the corresponding categories of equivariant matrix factorizations. In an Appendix written jointly with Alexey Basalaev, we apply the formulas to calculate the Hochschild cohomology of a simple but non-trivial class of so-called invertible LG orbifold models. The resulting algebras turn out to be isomorphic to what has already appeared in the literature on LG mirror symmetry under the name of twisted or orbifolded Milnor/Jacobian  algebras. We conjecture that this holds true for all invertible LG models. In the second part of the  Appendix, the formulas are applied to a different class of LG orbifolds which have appeared in the context of homological mirror symmetry for varieties of general type as  mirror partners of surfaces of genus 2 and higher. In combination with a homological mirror symmetry theorem for the surfaces, our calculation yields a new proof of the fact that the Hochschild cohomology of the Fukaya category of a surface is isomorphic, as an algebra, to the cohomology of the surface.

\end{abstract}

\maketitle

\vspace{-0.5in}
\baselineskip 1.44pc
\tableofcontents

\section{Introduction} The problem of extending the classical Hochschild-Kostant-Rosenberg theorem to various classes of spaces has long been a popular topic of research. In any geometric setting, the ultimate goal of this activity is to describe -- in geometric terms -- the Hochschild homology and cohomology of spaces in question  as completely as possible, including all parts of the rich structure that the Hochschild invariants carry: the cup and cap products, the Gerstenhaber bracket, pairings,  etc. This work is devoted to HKR like isomorphisms in the setting of {\it affine Landau-Ginzburg orbifolds}, i.e. triples $(X,W,G)$ where $X$ is a (say complex) smooth  affine variety, $W\in\C[X]$, and $G$ is a finite group of automorphisms of $(X,W)$. Our results  can be summarized as follows: For a  class of such orbifolds we present a complete description of the Hochschild cohomology and homology together with their cup and cap products.  

We should clarify what kind of Hochschild invariants of LG orbifolds we have in mind.  We adopt here the definition from \cite[Sect.6]{CT}: The Hochschild (co)homology of $(X,W,G)$ is the Hochschild (co)homology of the {\it curved algebra} $(\C[X]\rtimes G,W)$ where $\C[X]\rtimes G$ is the ordinary crossed product algebra associated with the $G$-action. A curved algebra \cite{Pos} is pair $(A, W)$ where $A$ is an associative not necessarily commutative algebra and $W\in A$ is a central element (``curvature''). A variant of the classical Hochschild theory for such objects was developed in \cite{CT,PP}. 

Let us explain  why the use of the physics terminology is legitimate here,  or in other words, what the invariants we study have to do with the actual topological LG model associated with $(X,W,G)$. 
According to a general philosophy \cite{Co,KR,KKP,Ko,KS1}, the truly useful Hochschild invariants associated with a topological string model are those of the corresponding D-brane category which in the LG case is the category ${\rm MF}_G(X,W)$ of $G$-equivariant matrix factorizations of $W$ \cite{ADD1,ADD2,QV}. The point is that the Hochschild (co)homology of  $(\C[X]\rtimes G,W)$ is {\it isomorphic} to the Hochschild (co)homology of ${\rm MF}_G(X,W)$. An explicit isomorphism was found  in \cite{PP,Se} and it turns out to preserve  the cup and cap products (see Section \ref{MFC}). In particular, the Hochschild cohomology of $(\C[X]\rtimes G,W)$ is (or should be) isomorphic to the closed string algebra of the LG model. (We have to admit though that we do not know how to describe the  topological metric on this algebra solely in terms of the Hochschild theory of $(\C[X]\rtimes G,W)$, whereas in terms of  ${\rm MF}_G(X,W)$ it can be done \cite{PV,Shk2}.) 

As mentioned above, this work is about explicit formulas for the products on the Hochschild (co)homology of $(\C[X]\rtimes G,W)$ for special $X$, $W$, and $G$.  Modulo minor details, those special triples look as follows:

\vspace{-0.1in}
\begin{enumerate}
\item $X$ is $\C^N$, that is, $\C[X]=\C[x_1,\ldots,x_N]$; 
\item $W$ has isolated critical points;
\item $G$ acts on $\C[X]$ by rescaling the variables $x_i$; in particular, $G$ is abelian.
\end{enumerate}
We would like to  comment on the situation with more general cases. Firstly, we actually treat a slight generalization of $\C^N$ which looks exotic but turns out useful for applications (as demonstrated in Section \ref{A.2}). However, abstract affine varieties are completely out of reach with our methods.  Secondly, the problem with non-abelian groups is easy to explain: our HKR isomorphisms  depend on the  coordinate system $(x_1,\ldots,x_N)$ and are equivariant only with respect to rescalings. We do not yet know how to bring different coordinate systems under one roof in a way compatible with the cup and cap products. Lastly, the requirement that $W$ have isolated singularities is more of a stylistic issue. In Section \ref{sect5}, while proving our theorem, we work out the general case. The problem is that our findings in the general case are not easy to formulate as a concise statement.

Let us outline the results of this paper in some more detail, focusing on the case of Hochschild  cohomology $\ccHH^*$ and its cup product. 

We assume from now on that $X$, $W$, and $G$ are as in (1)--(3) above. Although we think of $X$ as a variety, everything is linear in this setting and we will use the language of linear algebra and talk of subspaces of $X$ instead of subvarieties,  etc.  For $g\in G$ we will denote by $X^g$ the subspace of $g$-invariants in $X$. It has a unique $g$-invariant complement which we will denote by $X_g$.  Since $G$ is abelian, both $X^g$ and $X_g$ are $G$-stable. Let also $d_g:={\rm dim}\ X_g$. Finally, let  $W^g$ stand for the restriction of $W$ to $X^g$ and ${M}(W^g)$ for its Milnor algebra.

Let us first formulate an {\it additive} HKR theorem for the Hochschild cohomology. 
Consider the following $\ZZ$-graded vector space:
\begin{equation}\label{jacobian}
\Jacoii(X,W,G):=\bigoplus_{g\in G} M(W^g)\otimes {\rm det}(X_g) 
\end{equation}
where $M(W^g)$ and ${\rm det}(X_g)$ are placed in degrees 0 and   $d_g$ mod 2, respectively. The space  carries a degree preserving $G$-action coming from the $G$-action on $M(W^g)$ (induced by that on $\C[X^g]$) and the above-mentioned $G$-action on the subspace $X_g$ (recall that it is $G$-stable).
One has: 

\medskip

\noindent    
\begin{equation*}
\text{\emph{ 
There is an isomorphism of $\ZZ$-graded spaces}}\,\,\,
\ccHH^*(\C[X]\rtimes G,W)\simeq \Jacoii(X,W,G)^G.
\end{equation*}

\medskip

\noindent The claim in fact holds for arbitrary polynomials but with  $M(W^g)={\rm H}^0(\wedge^* T_{X^g}, [W^g,\cdot])$ replaced by ${\rm H}^*(\wedge^* T_{X^g}, [W^g,\cdot])$ where $[\cdot,\cdot]$ is the Schouten-Nijenhuis bracket. There are analogous results for the Hochschild homology. We prove all these claims in Section \ref{sect5} (see Propositions \ref{invsec1}, \ref{invsec2} and Section \ref{sect5.4}) but these facts seem to be known to the experts, in one version or another   \cite{BFK,CT,PV,Se}. 

The question we are really interested in is the following: {\it What is the product on $\Jacoii(X,W,G)^G$ that corresponds under the isomorphism to the cup product on the Hochschild cohomology?} There is another, even more interesting question, namely: {\it Is there a natural $G$-equivariant product on $\Jacoii(X,W,G)$ that gives rise to the correct product on the invariants?} In fact, such a product has to exist for very general reasons which we explain in Section \ref{sect3}. Moreover, we know some of its properties: 

\noindent (1) It preserves the $G$-grading: 
\[
\left(M(W^g)\otimes {\rm det}(X_g)\right)\cup \left(M(W^h)\otimes {\rm det}(X_h)\right)\subset
\left(M(W^{gh})\otimes {\rm det}(X_{gh})\right).
\]

\noindent (2) It is braided super-commutative: 
\begin{equation}\label{bc}
v_g\cup v_h=(-1)^{|v_g||v_h|} v_h\cup h^{-1}(v_g)
\end{equation}
for all $v_g\in M(W^g)\otimes {\rm det}(X_g)$ and   $v_h\in M(W^h)\otimes {\rm det}(X_h)$ ($|\cdot|$ denotes the $\ZZ$-degree). The induced product on the $G$-invariants is then super-commutative.

Let us make a guess as to what this product might look like. In the non-equivariant case the usual product on $M(W)$ can be interpreted as the product on the 0th cohomology of the complex $(\wedge^* T_{X}, [W,\cdot])$ induced by the wedge product on polyvector fields.
If we knew an equivariant analog of $\wedge^* T_{X}$ with its wedge product, we could try to predict the shape of the product on $\Jacoii(X,W,G)$. But, in fact, a good equivariant analog of $\wedge^* T_{X}$ {\it is} known:  it is the Hochschild cohomology of the crossed product $\C[X]\rtimes G$ itself. As shown in \cite{An,PPTT,SW2}, this Hochschild cohomology  is isomorphic, as an algebra, to the $G$-invariant part of 
$
\bigoplus_{g} \wedge^* T_{X^g}\otimes {\rm det}(X_g)
$
equipped with the  product
\[
(\mathfrak{X}_g\otimes \xi_g)\cup(\mathfrak{X}_h\otimes \xi_h)=
\begin{cases}(-1)^{d_g\cdot|\mathfrak{X}_h|}(\mathfrak{X}_g|_{X^{gh}}\wedge \mathfrak{X}_h|_{X^{gh}})\otimes (\xi_g\wedge\xi_h) &
\quad X^g\cap X^h=X^{gh}\,\,\,(\star)\\
0 & \quad\text{otherwise}
\end{cases}
\]
where $\mathfrak{X}_g\in \wedge^* T_{X^g}$, $\xi_g\in {\rm det}(X_g)$,  etc.  The transversality condition $(\star)$ is equivalent to $X_g\oplus X_h=X_{gh}$ and the second wedge product on the right-hand side is simply the canonical isomorphism ${\rm det}(X_g)\otimes {\rm det}(X_h)\ \tilde{\to}\ {\rm det}(X_g\oplus X_h)=\det(X_{gh})$. 

So, here is  our guess: The sought-after  product on $\Jacoii(X,W,G)$ is given by the formula
\begin{equation}\label{napr}
(\varphi_g\otimes \xi_g)\cup(\varphi_h\otimes \xi_h)=
\begin{cases}(\varphi_g|_{X^{gh}}\cdot \varphi_h|_{X^{gh}})\otimes (\xi_g\wedge\xi_h) &
\quad X^g\cap X^h=X^{gh}\\
0 & \quad\text{otherwise}
\end{cases}
\end{equation}
where $\varphi_g\in M(W^g)$ and $\varphi_h\in M(W^h)$. This product is easily seen to satisfy all the properties we want: it is $G$-equivariant, $G$-graded, and braided super-commutative. (But, in fact, it is also super-commutative in the ordinary sense!) Let us consider a simple example. 

Let $X=\C$, $W(x)=x^3$, and $G=\Z/3\Z=\{0,1,2\}$
where $G$ acts by multiplication by the cubics roots of unity: ${n}\mapsto \zeta^n$ for $\zeta$ a primitive root. In this case $X_{{1}}=X_{{2}}=\C$ and therefore, as a $\ZZ$-graded space,
\[
{\rm M}^{\rm even}(X,W,G)=\C\xi_{{0}}\oplus\C x\xi_{{0}},\quad {\rm M}^{\rm odd}(X,W,G)=\C\xi_{{1}}\oplus \C\xi_{{2}}
\]
where  $\xi_i$ stands for a generator of $\det(X_i)$. The product (\ref{napr}) is quite boring in this case: $\xi_{{0}}$ is the unit and all other elements multiply to 0. 

Actually, our guess is {\it wrong}:  (\ref{napr}) is not the sought-after product on $\Jacoii(X,W,G)$ and it does not induce the right cup product on $\Jacoii(X,W,G)^G$. The  relation between the product (\ref{napr}) and the actual one turns out to be very similar to the relation between the classical and quantum cohomology of a symplectic manifold. Namely, $\Jacoii(X,W,G)$ with the above product is the limit of $\Jacoii(X,tW,G)$ with the true product as $t\to0$. 

The true product on $\Jacoii(X,tW,G)$ has the following form:  For  any fixed choice of generators $\{\xi_g\}_{g\in G}$ of the one-dimensional spaces $\det(X_g)$ there are elements $\sigma_{g,h}\in M(W^{gh})$ such that 

\begin{equation}\label{trpr}
(\varphi_g\otimes \xi_g)\cup(\varphi_h\otimes \xi_h)=
\begin{cases}t^{\frac{d_g+d_h-d_{gh}}{2}}(\sigma_{g,h}\cdot \varphi_g|_{X^{gh}}\cdot \varphi_h|_{X^{gh}})\otimes \xi_{gh} &
\quad \frac{d_g+d_h-d_{gh}}{2}\in \Z_{\geq0}\\
0 & \quad\text{otherwise}
\end{cases}
\end{equation}
As we mentioned above, the condition $(\star)$ in (\ref{napr}) is equivalent to $d_g+d_h-d_{gh}=0$. When the condition is satisfied, the power of $t$ in (\ref{trpr}) disappears and the corresponding products survive the limit $t\to0$  (and, in fact, become what we had before). All other products tend to 0. 

Our main result -- Theorem \ref{mainthm} -- gives  explicit, though quite complicated, formula for $\sigma_{g,h}$. Together with (\ref{trpr}) (for  $t=1$) this provides a complete description of the product on $\Jacoii(X,W,G)$ and, consequently, a complete description of the product on $\Jacoii(X,W,G)^G\simeq\ccHH^*(\C[X]\rtimes G,W)$. 

In the one-dimensional example we discussed above, the correct product on  $\Jacoii(X,W,G)$ differs from the naive one in that the product of $\xi_{{1}}$ and $\xi_{{2}}$ is not 0 anymore: up to a renormalization of the generators $\xi_{{i}}$ one has 
\begin{equation}\label{x3}
\xi_{{1}}\cup\xi_{{2}}=\frac{1}{\zeta-1} x \xi_0,\quad \xi_{{2}}\cup\xi_{{1}}=\frac{1}{\zeta^{-1}-1} x \xi_0.
\end{equation}
Even in this simple example one can observe that the product on $\Jacoii(X,W,G)$ is indeed not  super-commutative in general but  braided super-commutative. 

In fact, the formulas (\ref{x3}) can be compared with something that has already appeared in the literature. The point is that $x^3$ is (almost) the simplest example of an {\it invertible polynomial} \cite{KS}, i.e. a quasi-homogeneous polynomial with an isolated critical point at the origin having the same number of variables and monomials. These polynomials  have been studied quite extensively due to an important role they play in LG mirror symmetry and Fan-Jarvis-Ruan-Witten theory \cite{BH,FJ,Kre}. In particular, there already exists an analog of $\Jacoii(X,W,G)$ in this setting which was constructed ``by hand'' in \cite{Kra}  (building on pioneering ideas of \cite{Kau0,Kau1,Kau2}) and later used in \cite{FJJS} to prove LG mirror symmetry at the level of  Frobenius algebras. In the recent work \cite{BTW1}, a more systematic study of these algebras was undertaken and, in particular, an improved version of the original definition was proposed (under the name of {\it the $G$-twisted Jacobian algebra} of $W$) which satisfies various expected properties, e.g. the braided super-commutativity. For $W(x)=x^3$ this $G$-twisted Jacobian algebra can easily be seen to be isomorphic  to $\Jacoii(X,W,G)$. In Appendix \ref{appA.1} the two algebras are compared in less trivial examples and turn out to be isomorphic there as well. We believe this holds true for all invertible polynomials. 

There is yet another class of LG orbifolds against which our formulas can be tested, namely, certain orbifolded cusp singularities (in dimension 3) which are shown in \cite{Ef,Sei} to be homological mirror partners of two dimensional symplectic surfaces. In Appendix \ref{A.2} we calculate  $\Jacoii(X,W,G)^G$ for these LG models and show that the resulting algebras are isomorphic to the cohomology of the corresponding surfaces, as it should be for certain general reasons \cite{Gan,GPS1,GPS2}.

We would like to conclude the Introduction by confessing that we do not yet understand the geometric meaning of our formulas for the above-mentioned ``structure constants'' $\sigma_{g,h}$. Perhaps, they could somehow be related to the Chern characters of matrix factorizations of $W(x)-W(y)$ but we have no evidence to support this idea. In any case, understanding what  $\sigma_{g,h}$ mean geometrically should help to extend the results beyond the limited setting of the present work.

\bigskip

\noindent{\bf Acknowledgements.} I am grateful to Alexey Basalaev for following the progress of the work with constant interest, for testing preliminary formulas (and spotting inconsistencies in some of them) and, of course, for collaborating on Appendix \ref{A}. Many thanks also to Sheel Ganatra for answering my questions about \cite{Gan} and further comments and explanations. Last but not least, I would like to thank Christian Sevenheck, who was the first to read my earlier writings on this topic and provided many helpful suggestions.

\section{Hochschild invariants of curved crossed product algebras}\label{Sect.2}

In this section, and throughout the paper, $\K$ is  an arbitrary field of characteristic 0. 

\subsection{Curved Hochschild calculus}\label{sect2} 
\subsubsection{Outline}\label{sect2.1} 
For any algebra $A$ the Hochschild cohomology functor 
$
\ccHH^*(A,-)
$ from the category of $A$-bimodules to that of graded vector spaces
carries a natural monoidal structure;  the corresponding  maps 
$
\ccHH^*(A,M_1)\otimes \ccHH^*(A,M_2)\to \ccHH^*(A,M_1\otimes_AM_2)
$ are usually referred to as ``{cup products}''. Also, the Hochschild homology functor $\ccHH_*(A,-)$ has a natural structure of an $\ccHH^*(A,-)$-module which is encoded in ``{cap} products''.
The combination of these two structures is what is called  {\it Hochschild calculus} in this paper. 
The aim of Section \ref{sect2} is to discuss a counterpart of the Hochschild calculus for curved algebras. In fact, such a counterpart,  in the much broader context of curved dg categories,  already exists \cite{PP} (albeit without mentioning the products explicitly). In particular, one has a notion of Hochschild (co)homology ``of the second kind''  of a curved algebra with coefficients in a curved  bimodule.  However, for the purposes of this work the full power of the theory developed in \cite{PP} is not needed. We will only be interested in curved  bimodules sitting in degree 0 (We call them simply bimodules, without the adjective ``curved''.) For such bimodules, the theory of \cite{PP} can be streamlined by implementing the language of mixed complexes.

\subsubsection{Mixed complexes}\label{mc}

Recall \cite{Kas}\footnote{Unlike in \cite{Kas}, our mixed complexes are not necessarily left or right bounded.} that a mixed complex is a triple $(\rC, b, B)$ where $\rC=\oplus_{n\in\Z}\rC_n$ is a $\Z$-graded vector space and $b$ and $B$ are operators on $\rC$ of degrees 1 and $-1$, respectively, satisfying
\[
b^2=0,\quad B^2=0,\quad bB+Bb=0.
\] 
The degree of $c\in\rC$ will be denoted by $|c|$. A morphism $(\rC, b, B)\to (\rC', b', B')$ is by definition a degree preserving map $\rC\to\rC'$  commuting with both differentials. Such a morphism is a called a  quasi-isomorphism if the induced morphism of complexes $(\rC, b)\to (\rC', b')$ is a quasi-isomorphism.

Mixed complexes form a tensor category under
\[
(\rC, b, B)\otimes (\rC', b', B');=(\rC\otimes \rC', b\otimes1+1\otimes b', B\otimes1+1\otimes B')
\]
where $\otimes$ on the right-hand side is the usual tensor product in the category of graded vector spaces (namely, $(\rC\otimes \rC')_n=\oplus_{p+q=n}\rC_p\otimes \rC'_q$ and the new differentials pick up signs when applied to the elements of  $\rC\otimes \rC'$, in agreement with the Koszul rule of signs). 

A mixed complex $(\rC, b, B)$ together with a morphism $(\rC, b, B)^{\otimes 2}\to (\rC, b, B)$ satisfying the associativity and unitality conditions will be called a mixed dg algebra. The notion of a  mixed dg module over a mixed dg algebra is defined similarly.

Let us fix now a formal variable ${\uu}$ of degree 2. Given a graded vector space $\rC=\oplus_{n\in\Z}\rC_n$, we will denote by $\rC\bbb$ the graded $\K\bb$-module spanned (over $\K$) by homogeneous formal Laurent series in ${\uu}$ with coefficients in $\rC$:
\begin{equation}\label{Ct}
\rC\bbb^n=\left\{\sum_{i=i_0}^\infty c_i\uu^i\,|\, c_i\in \rC_{n-2i}\right\}.
\end{equation}
The correspondence $\rC\mapsto\rC\bbb$ can be promoted to a functor from the category of mixed complexes to that of graded $\K\bb$-modules, namely:
\begin{equation}\label{per}
(\rC, b, B)\mapsto {\rm H}^*(\rC\bbb, b+ {\uu} B).
\end{equation}
We call ${\rm H}^*(\rC\bbb, b+ {\uu} B)$ the periodic cohomology of $(\rC, b, B)$.  

Let us point out two properties of the functor (\ref{per}) which we will use in the future. 
Firstly,  the functor is lax monoidal. In particular, it  transforms mixed dg algebras into $\K\bb$-linear algebras and mixed dg modules over the former into $\K\bb$-linear  modules over the latter. 
Secondly, the functor transforms quasi-isomorphisms into isomorphisms \cite[Prop.2.4]{GJ0}.

\subsubsection{Hochschild mixed complexes}\label{sect2.3}
 
Recall that the bar resolution of an associative unital algebra $A$ is the complex $(\rB_*(A), \ddbar)$ of $A$-bimodules with 
$
\rB_{-n}(A):= A\otimes A^{\otimes n}\otimes A
$ ($n=0,1,\ldots$) and
\begin{multline*}
\ddbar(a_0[a_1|\ldots |a_n]a_{n+1})=a_0a_1[a_{2}|\ldots |a_n]a_{n+1}+\sum\limits_{i=1}^{n-1}(-1)^{i}a_0[a_1|\ldots
|a_ia_{i+1}|\ldots|a_n]a_{n+1}\\+(-1)^{n}a_0[a_1|\ldots|a_{n-1} ]a_na_{n+1}.
\end{multline*}
where $a_0[a_{1}|\ldots|a_n]a_{n+1}$ is shorthand for $a_0\otimes a_{1}\otimes \ldots \otimes a_n\otimes a_{n+1}$.

Let now $W\in A$ be a central element. Associated with $W$ there is a degree $-1$ differential on $\rB_*(A)$, namely 
\begin{eqnarray*}
\ddcurve(a_0[a_1|\ldots |a_n]a_{n+1})=\sum\limits_{i=0}^{n}(-1)^{i}a_0[a_1|\ldots
|a_i|W|a_{i+1}|\ldots|a_n]a_{n+1}.
\end{eqnarray*}
One has
\begin{equation}\label{bar}
\delta_{\rm curv}^2=0,\quad \ddbar\ddcurve+\ddcurve\ddbar=l_W-r_W
\end{equation}
where $l_W$ (resp. $r_W$) is the operator in $\rB_*(A)$ of left (resp. right) multiplication with $W$. 

\begin{remark}{\rm {\it All} the results and conclusions in this paper remain valid if one start with the {\it normalized} version $(\overline{\rB}_*(A),\ddbar, \ddcurve)$ where $\overline{\rB}_{-n}(A):= A\otimes (A/\K)^{\otimes n}\otimes A$.
}
\end{remark}

Let $M$ be an $(A,W)$-bimodule, i.~e. an $A$-bimodule in which the operators $l_W$ and $r_W$ of left and right multiplications with $W$ coincide: 
\begin{equation}\label{bimod}
l_W-r_W=0.
\end{equation} 
Let $\rB^*(A, M):=\Hom_{A\otimes A^{\rm op}}(\rB_{-*}(A), M)$ where $A^{\rm op}$ denotes the opposite algebra. (Note that $\rB^*(A, M)$ is non-negatively graded.) The pairing of an element $D\in \rB^*(A, M)$ with an element $\underline{a}\in \rB_*(A)$ will be written as $\langle \underline{a}, D\rangle\in M$; so $\langle a_0\,\underline{a}\,a_1, D\rangle=a_0\langle \underline{a}, D\rangle a_1$ for $a_0,a_1\in A$. 

\begin{remark}\label{lr}{\rm We will also use $\langle\cdot,\cdot\rangle$ in a more general sense, namely, to denote the natural pairing $X\otimes \Hom_{A\otimes A}(X,Y)\to Y$ for {\it any} $A$-bimodules $X,Y$. 
}
\end{remark}

The Hochschild cochain mixed complex of $(A,W)$ with coefficients in $M$ is defined by 
\begin{equation*}
\rHH^{*}(A,W;M):=\left(\rB^*(A, M), \,\,\ddh:=\delta_{\rm bar}^\vee, \,\,\ddw:=\delta_{\rm curv}^\vee\right)
\end{equation*}
where $\delta^\vee$ denotes the standard dual of  $\delta$ defined by $\langle -, \delta^\vee(D)\rangle=(-1)^{|D|} \langle \delta(-), D\rangle$. (That it is indeed a mixed complex follows from (\ref{bar}) and (\ref{bimod}).) We will denote its periodic cohomology  by 
$
\cHH^{*}(A,W;M)
$, or simply $\cHH^{*}(A,W)$ when $M=A$.

Let $\rB_*(A, M):= M\otimes_{A\otimes A^{\rm op}}\rB_*(A)$ where the right-hand side is an abbreviation for 
\[
M\otimes\rB_*(A)/\{a_0\,m\,a_1\otimes \underline{a}-m\otimes a_1\,\underline{a}\,a_0\,|\,a_0,a_1\in A, m\in M, \underline{a}\in \rB_*(A)\}.
\] 
(Note that $\rB_*(A, M)$ is non-positively graded.)
The  Hochschild chain mixed complex of $(A,W)$ with coefficients in $M$ is defined as the mixed complex
\begin{equation*}
\rHH_{*}(A,W;M):=\left(\rB_*(A, M), \,\,\bbh:=1\otimes \delta_{\rm bar},\,\, \bbw:=1\otimes\delta_{\rm curv}\right)
\end{equation*}
Its periodic cohomology   will be denoted  by
$
\cHH_{*}(A,W;M)
$ or $\cHH_{*}(A,W)$ when $M=A$.

\begin{remark}\label{HHtHHII}{\rm As we have already mentioned in Section \ref{sect2.1}, ``our''  Hochschild (co)homology  is nothing but a special case  of the Hochschild (co)homology of the second kind of  $(A,W)$ with coefficients in a {\it curved} $(A,W)$-bimodule introduced in \cite{PP}. However, the reader familiar with \cite{PP} will notice that even in the special case our definition does not match the one in \cite{PP} because the variable $t$ does not appear there. The (co)homology $\ccHH^{\rm II,*}(A,W;M)$ and $\ccHH^{\rm II}_*(A,W;M)$ of \cite{PP} are the {\it $\ZZ$-graded} spaces defined as the cohomology of the $\ZZ$-graded complexes 
\begin{eqnarray*}
\left(\bigoplus_{i\,\,\text{even}}\rB^i(A, M)\oplus \bigoplus_{i\,\,\text{odd}}\rB^i(A, M), \,\,\ddh+\ddw\right)
\end{eqnarray*}
and
\begin{eqnarray*}
\left(\prod_{i\,\,\text{even}}\rB_i(A, M)\oplus \prod_{i\,\,\text{odd}}\rB_i(A, M), \,\,\bbh+\bbw\right),
\end{eqnarray*}
respectively.  It is easy to see that there is a straightforward relation between the two definitions, namely, $\cHH$ is just a 2-periodic $\Z$-graded version of  $\ccHH^{\rm II}$: 
\begin{equation*}
\ccHH^{\rm II,even}\simeq  \cHH^{2n},\quad \ccHH^{\rm II, odd}\simeq  \cHH^{2n+1}\qquad \forall n
\end{equation*}
and the same for homology. As a consequence, all the results we obtain in the present work have $\ZZ$-graded counterparts for $\ccHH^{\rm II}$.

}
\end{remark}

\subsubsection{The cup and cap products}\label{cip}

Let 
$
\Delta=\Delta_{\rm bar}: \rB_*(A)\to \rB_*(A)\otimes_A\rB_*(A)
$ be the morphism of $A$-bimodules given by
\begin{equation}\label{coalg0}
\Delta(a_0[a_1|\ldots |a_n]a_{n+1})=\sum_{i=0}^n (a_0[a_1|\ldots |a_i]1)\otimes (1[a_{i+1}|\ldots |a_n]a_{n+1}).
\end{equation}
This morphism is coassociative, i.~e. 
$
(\Delta\otimes1)\Delta=(1\otimes\Delta)\Delta,
$
and is easily seen to be compatible with the differentials $\ddbar$ and $\ddcurve$:
\begin{equation*}
\Delta\cdot \ddbar=(\ddbar\otimes1+1\otimes \ddbar)\cdot\Delta, \quad \Delta\cdot \ddcurve=(\ddcurve\otimes1+1\otimes \ddcurve).
\end{equation*}

Given two $(A,W)$-bimodules $M_1$ and $M_2$, $\Delta$ induces the {\it cup product}
\[
\cup\,=\,\cup_{\rm Hoch}: \rB^*(A, M_1)\otimes \rB^*(A, M_2)
\to \rB^*(A, M_1\otimes_AM_2),
\]
\begin{equation}\label{cup}
\langle \underline{a}, D_1\cup D_2\rangle:=\langle \Delta(\underline{a}), D_1\boxtimes D_2\rangle=(-1)^{|D_1||\underline{a}_{(2)}|}\langle\underline{a}_{(1)},D_1\rangle\otimes \langle\underline{a}_{(2)},D_2\rangle \quad (\forall\, \underline{a}\in \rB_*(A)),
\end{equation}
as well as the {\it cap product}
\[
\cap\,=\,\cap_{\rm Hoch}: \rB_*(A, M_1)\otimes \rB^*(A,M_2)\to \rB_*(A,M_1\otimes_AM_2),
\]
\begin{equation}\label{int}
(m\otimes \underline{a})\cap D:=m\otimes\langle\Delta(\underline{a}), D\boxtimes {\rm id}_{\rB_*(A)}\rangle
=(-1)^{|D||\underline{a}_{(2)}|}(m\otimes\langle\underline{a}_{(1)},D\rangle)\otimes\underline{a}_{(2)}
\end{equation}
where $\underline{a}_{(1)}\otimes \underline{a}_{(2)}:=\Delta(\underline{a})$ (Sweedler's notation) and $\boxtimes$ denotes the natural map 
\begin{equation}\label{sqten}
\Hom_{A\otimes A}(X_1,Y_1)\otimes \Hom_{A\otimes A}(X_2,Y_2)\to\Hom_{A\otimes A}(X_1\otimes_A X_2,Y_1\otimes_A Y_2)
\end{equation}
for $A$-bimodules $X_i$, $Y_i$. 
The following facts follow easily from the definitions:
\begin{proposition}\label{cupcoh} 
(1) The cup product is a morphism of  mixed complexes and hence yields a product 
\[
\cup: \cHH^{*}(A,W;M_1)\otimes_{\K\bb} \cHH^{*}(A,W;M_2)\to \cHH^{*}(A,W;M_1\otimes_AM_2)
\]
\noindent(2) For any  $(A,W)$-bimodules $M_i$, $i=1,2,3$, the following diagram is commutative:
\[\begin{CD}
\rB^{*}(A,M_1)\otimes \rB^{*}(A,M_2)\otimes \rB^{*}(A,M_3) @>{{\rm id}\otimes\cup}>>\rB^{*}(A,M_1)\otimes \rB^{*}(A,M_2\otimes_A M_3)\\
@VV{\cup\otimes{\rm id}}V @VV\cup V\\
\rB^{*}(A,M_1\otimes_A M_2)\otimes \rB^{*}(A,M_3) @>\cup>> \rB^{*}(A,M_1\otimes_AM_2\otimes_A M_3)
\end{CD}\]
In particular, $\rHH^{*}(A,W)$ is an associative mixed dg algebra\footnote{The unit is the element of $\Hom\left(\K, A\right)$ (0-cochain) sending the unit of $\K$ to the unit of $A$.} and for any $(A,W)$-bimodule $M$ $\rHH^{*}(A,W;M)$ is a mixed dg bimodule over $\rHH^{*}(A,W)$.
As a consequence, $\cHH^{*}(A,W)$ is an associative $\K\bb$-linear algebra and  $\cHH^{*}(A,W;M)$ is an $\cHH^{*}(A,W)$-bimodule.

\noindent(3)  The cap product is a morphism of  mixed complexes and therefore yields a product
\[
\cap: \cHH_{*}(A,W;M_1)\otimes_{\K\bb} \cHH^{*}(A,W;M_2)\to \cHH_{*}(A,W;M_1\otimes_AM_2)
\] 
\noindent(4) For $M_i$ as above, the diagram 
\[\begin{CD}
\rB_{*}(A,M_1)\otimes \rB^{*}(A,M_2)\otimes \rB^{*}(A,M_3) @>{{\rm id}\otimes\cup}>>\rB_{*}(A,M_1)\otimes \rB^{*}(A,M_2\otimes_A M_3)\\
@VV{\cap\otimes{\rm id}}V @VV\cap V\\
\rB_{*}(A,M_1\otimes_A M_2)\otimes \rB^{*}(A,M_3) @>\cap>> \rB_{*}(A,M_1\otimes_AM_2\otimes_A M_3)
\end{CD}\]
is commutative. In particular,   $\rHH_{*}(A,W;M)$, for any $(A,W)$-bimodule $M$, is a right mixed dg module over $\rHH^{*}(A,W)$ and $\cHH_{*}(A,W;M)$ is a right  $\cHH^{*}(A,W)$-module.
\end{proposition}

\begin{remark}\label{rmk1}{\rm Just as in the  non-curved case \cite{Ge}  $\cHH^{*}(A,W)$ turns out to be super-commutative. This is a special case of a stronger result which we will discuss in Section \ref{brcom}.}
\end{remark}

\subsubsection{K\"unneth isomorphisms for Hochschild calculi} Let $(A,W)$ and $(A',W')$ be two curved algebras. Their tensor product is defined as the curved algebra $(A\otimes A', W\otimes 1+1\otimes W')$. Observe that for an $(A,W)$-bimodule $M$ and an $(A',W')$-bimodule $M'$ the $A\otimes A'$-bimodule $M\otimes M'$ is actually  an $(A\otimes A', W\otimes 1+1\otimes W')$-bimodule. Thus, it is natural to ask if, and how, the Hochschild calculus of   $(A\otimes A', W\otimes 1+1\otimes W')$ can be ``calculated'' in terms of the Hochschild calculi of $(A,W)$ and $(A',W')$. It is indeed possible, at least for a class of algebras, namely,

\vspace{0.1in}

\noindent \emph{Recall \cite{KS} that an algebra $A$ is called (homologically) smooth if, as a bimodule over itself,  it admits a bounded below resolution by finitely generated projective $A$-bimodules.
}

\begin{proposition}\label{Kunn} If $A$ and $A'$ are smooth then there are natural isomorphisms
\begin{eqnarray}\label{kunn}
\cHH^*(A,W;M)\otimes_{\K\bb} \cHH^*(A',W';M') \simeq \cHH^*(A\otimes A', W\otimes 1+1\otimes W'; M\otimes M'),\nonumber\\
\cHH_*(A,W;M)\otimes_{\K\bb} \cHH_*(A',W';M') \simeq \cHH_*(A\otimes A', W\otimes 1+1\otimes W'; M\otimes M'),
\end{eqnarray}
 compatible with the cup and cap products. 
\end{proposition}
\noindent{\bf Proof} is given in Appendix \ref{B} (page \pageref{Kunnproof}).

\subsection{Hochschild (co)homology of equivariant curved algebras}\label{sect3}

\subsubsection{Outline} Suppose $A$ is acted upon by a finite symmetry group $G$ (i.e. one has a group homomorphism $G\to {\rm Aut}(A)$).   Recall that  the crossed product  $A\rtimes  G$ is defined as the algebra
\[
A\otimes \K[G]=\bigoplus_{g\in G} A\otimes g,\quad (a\otimes g)\cdot(b\otimes h):= ag(b)\otimes gh,\quad a,b\in A,\,\, g,h\in G.
\]
 A $G$-invariant element $W\in A$ gives rise to the curved algebra 
$
(A\rtimes  G,W)
$. 
Our aim in this section is to develop an analog for $(A\rtimes  G,W)$ of some standard technique \cite[...]{An,Ba,CGW,Ka,PPTT,SW2} that allows one to calculate the (co)homology of $A\rtimes  G$ in terms of the Hochschild calculus of $A$ with coefficients in $A\rtimes  G$.

\subsubsection{$G$-twisted Hochschild (co)homology}\label{sect3.1}

The multiplication on $A\rtimes  G$ endows each subspace $A\otimes g\subset A\rtimes  G$ with an $(A,W)$-bimodule structure and induces isomorphisms of $(A,W)$-bimodules 
\begin{equation}\label{twiso}
(A\otimes g)\otimes_A(A\otimes h)\to A\otimes gh.
\end{equation}
Then Proposition \ref{cupcoh} implies that
\[
\rHH^*(A,W;A\rtimes  G)=\bigoplus_{g\in G} \rHH^*(A,W;A\otimes g),\quad \rHH_*(A,W;A\rtimes  G)=\bigoplus_{g\in G} \rHH_*(A,W;A\otimes g)
\]
have natural structures of a mixed dg algebra and a right mixed dg module over the dg algebra, respectively. Notice that the additional $G$-grading on both of them is compatible with the cup and cap products. (Perhaps, it is worthwhile mentioning explicitly that the  $G$-grading is {\it non-homological}, in the sense that it does not affect the Leibniz rule.) Consequently,
\[
\cHH^*(A,W;A\rtimes  G)=\bigoplus_{g\in G} \cHH^*(A,W;A\otimes g),\quad \cHH_*(A,W;A\rtimes  G)=\bigoplus_{g\in G} \cHH_*(A,W;A\otimes g)
\]
have structures of an associative $(\Z\times G)$-graded $\K\bb$-linear algebra and a right $(\Z\times G)$-graded  module over this algebra, respectively.

Furthermore, observe that $G$ acts on $A\rtimes  G$ by conjugation: $h(a\otimes g)=h(a)\otimes hgh^{-1}$ ($h\in G$), and this action induces isomorphisms (of vector spaces) $h: A\otimes g\to A\otimes hgh^{-1}$ which are compatible with the isomorphisms (\ref{twiso}). Also, $G$ acts, in the obvious manner, on the bar resolution $\rB_*(A)$, and this action respects the differentials $\ddbar$ and $\ddcurve$, as well as the coproduct $\Delta_{\rm bar}$. Combining the above two actions one gets well-defined $G$-actions on  $\rB^*(A,A\rtimes  G)$ and $\rB_*(A,A\rtimes  G)$, namely  
\begin{eqnarray}\label{Gactcoh}
\rB^*(A, A\otimes g)\ni D&\mapsto& h(D):=h\circ D\circ h^{-1}\in \rB^*(A, A\otimes hgh^{-1}),\nonumber\\
\rB_*(A, A\otimes g)\ni m\otimes \underline{a}&\mapsto& h(m\otimes \underline{a}):=h(m)\otimes h(\underline{a})\in \rB_*(A, A\otimes hgh^{-1}).
\end{eqnarray}
These $G$-actions are easily seen to be compatible with all the structures we are interested in:

\noindent(1) They preserve the $\Z$-gradings and commute with the differentials, thereby inducing  $G$-actions on the mixed complexes $\rHH^*(A,W;A\rtimes  G)$ and $\rHH_*(A,W;A\rtimes  G)$.

\noindent(2) 
 They commute with the cup and cap products:
\begin{equation}\label{compat}
k(D_1\cup D_2)=k(D_1)\cup k(D_2),\quad  k(\omega\cap D)=k(\omega)\cap k(D)
\end{equation}
for all $k\in G$, $D_1\in\rB^*(A,A\otimes g)$,  $D,D_2\in \rB^*(A,A\otimes h)$, and $\omega\in\rB_*(A,A\otimes g)$.

\subsubsection{Hochschild (co)homology of the curved crossed product} It follows from (\ref{compat}) that the cup and cap products on $\rHH^*(A,W;A\rtimes  G)$ and $\rHH_*(A,W;A\rtimes  G)$ descend to well-defined cup and cap products on the mixed complexes $\rHH^*(A,W;A\rtimes  G)^G$ and $\rHH_*(A,W;A\rtimes  G)_G$ where $(\,)^G$ resp. $(\,)_G$ denote the space of $G$-invariants resp. $G$-coinvariants. One has:
\begin{proposition}\label{Ginv} There are isomorphisms of $\Z$-graded $\K\bb$-modules
\[
\cHH^*(A\rtimes  G,W)\simeq \cHH^*(A,W;A\rtimes  G)^G,\quad \cHH_*(A\rtimes  G,W)\simeq \cHH_*(A,W;A\rtimes  G)_G
\]
compatible with the cup and cap products.
\end{proposition}
\noindent{\bf Proof} is given in Appendix \ref{B} (page \pageref{123}). 

\subsubsection{Braided (a.k.a. $G$-twisted) commutativity}\label{brcom}
We conclude this section by formulating an equivariant curved analog of the classical result of \cite{Ge} on the super-commutativity of the Hochschild cohomology of ordinary algebras. 

\begin{proposition}\label{brcomm} The algebra $\cHH^*(A,W;A\rtimes  G)$ is {\rm braided super-commutative}: For all homogeneous classes $[D_1]\in\cHH^*(A,W;A\otimes g)$,  $[D_2]\in \cHH^*(A,W;A\otimes h)$ 
\begin{equation}\label{bracom}
[D_1]\cup [D_2]=(-1)^{|[D_1]||[D_2]|} [D_2]\cup h^{-1}([D_1]).
\end{equation}
\end{proposition}
\noindent{\bf Proof}  is given in Appendix \ref{B} (page \pageref{Gtwcomproof}). 

\vspace{0.1in}

When $G$ is trivial, the braided super-commutativity is just the ordinary super-commutativity, so 

\begin{corollary}\label{cge} For any $(A,W)$ the algebra $\cHH^{*}(A,W)$ is super-commutative. 
\end{corollary}

%%%%%%%%%

\section{Main results}
\subsection{Setting and notation}\label{SN} 
From now on, we focus on the following case:
\begin{itemize}
\item $A=\K[{X}]:=\K[x_1,x_2,\ldots,x_N][S^{-1}]$ where $S$ is the multiplicative set generated by a finite (or empty) set of polynomials of the form $(x_i-\lambda)$ with $\lambda\in\K^*$. 
\item The ``curvature'' is just a regular function $W=W(x)\in\K[{X}]$ ($x:=(x_1,\ldots, x_N)$). In the statement of the main result we will require $W$ to have isolated critical points by which we understand the condition that the cohomology of the complex $(\Omega^*_X, dW\wedge\cdot)$ vanishes in degrees less than $N$.
\item The group $G$ acts on $\K[X]$ by rescaling the variables: 
\begin{equation}\label{ge}
(\K^*)^N\ni g=(g_1, \ldots, g_{N}): (x_1,\ldots, x_N)\mapsto (g_1x_1,\ldots, g_Nx_N)=:g(x).
\end{equation}
   
\end{itemize}
Let us introduce some notation which will be used in the statement of the main theorem, as well as throughout the proof.

\subsubsection{Difference derivatives}\label{DD} We will write the elements of $\K[X]^{\otimes 2}$ as functions of two sets of $N$ variables where the second set is ${y}=(y_1,\ldots, y_N)$. Similarly, the elements of $\K[X]^{\otimes 3}$ will be written as $f(x,y,z)$.

Observe that because of the special form of the above multiplicative set $S$, we have the following  well-defined maps:
\begin{eqnarray}\label{deltai}
\nabla_i=\nabla^{x\to (x,y)}_i: \K[{X}]\to\K[{X}]^{\otimes 2}\quad (i=1,\ldots, N),\qquad \nabla_i(f):=\frac{l_i(f)-l_{i+1}(f)}{x_i-y_i}.
\end{eqnarray}
where $l_i(f):=f(y_1,\ldots,y_{i-1},x_i,\ldots,x_N)$ for $i=1,\ldots, N+1$.
(In particular, $l_1f=f\otimes 1=f(x)$ and $l_{N+1}f=1\otimes f=f(y)$.)  
Note that  
\begin{equation}\label{diffder}
\sum_{i=1}^N(x_i-y_i)\nabla_i(f)=f(x)-f(y).
\end{equation}

The symbol $\nabla^{y\to(y,z)}\nabla^{x\to(x,y)}$ will have the following meaning: $\nabla^{x\to(x,y)}$ is applied to a function of $x$ and produces a function of $(x,y)$; then $\nabla^{y\to(y,z)}$ is applied to this new function viewed as a function of $y$, with $x$ ``frozen''; the result is a function of $(x,y,z)$.

The symbol $\nabla^{x\to(x,\psi(x))}$, where $\psi(x)$ is some function of $x$, will have the following meaning:  $\nabla^{x\to(x,y)}$ is applied to a function of $x$ and then $\psi(x)$ is substituted for $y$; the result is a new function of $x$. The symbol $\nabla^{x\to(x,\psi'(x))}\nabla^{x\to(x,\psi''(x))}$ is the composition of two operations of this type for two different functions of $x$; the result is again a function of $x$.

\subsubsection{Clifford algebra}\label{CA} We will denote by ${\rm Cl}_N$ the $N$th $\Z$-graded Clifford algebra: 
\[
{\rm Cl}_N=\K[\theta_1,\ldots,\theta_N,{\partial_{\theta_1}},\ldots,{\partial_{\theta_N}}], \quad |\theta_i|=-1,\quad |{\partial_{\theta_i}}|=1\quad \forall i
\] modulo the relations
\begin{eqnarray*}
\theta_i\theta_j=-\theta_j\theta_i,\quad{\partial_{\theta_i}}{\partial_{\theta_j}}=-{\partial_{\theta_j}}{\partial_{\theta_i}},\quad{\partial_{\theta_i}}\theta_j=-\theta_j{\partial_{\theta_i}}+\delta_{ij}.
\end{eqnarray*}
 For $I\subset\{1,\ldots,N\}$ we write
$
\partial_{\theta_{I}}:=\prod_{i\in I}\partial_{\theta_i},   
$
$
{\theta_{I}}:=\prod_{i\in I}{\theta_i}
$
where in both cases  the multipliers are taken in increasing order of the indices.

We endow the subspaces $\K[{\theta}]=\K[\theta_1,\ldots,\theta_N]$ and $\K[\partial_{\theta}]=\K[\partial_{\theta_1},\ldots,\partial_{\theta_N}]$ with left $\Z$-graded ${\rm Cl}_N$-module structures via the isomorphisms  
\[
\K[{\theta}]\simeq {\rm Cl}_N/{\rm Cl}_N\cdot\langle{\partial_{\theta_1}},\ldots,{\partial_{\theta_N}}\rangle, \quad \K[\partial_{\theta}]\simeq {\rm Cl}_N/{\rm Cl}_N\cdot\langle{{\theta_1}},\ldots,{{\theta_N}}\rangle.
\]  
Given a $\Z$-graded space $\rC$, the spaces $\rC[{\theta}]:=\rC\otimes \K[{\theta}]$ and $\rC[\partial_{\theta}]:=\rC\otimes \K[\partial_{\theta}]$ are endowed with ${\rm Cl}_N$-module structures by (graded) $\rC$-linearity.

\subsubsection{Notation related to the group action}\label{GA} For an element $g\in G$ as in (\ref{ge}) we define
\[
I^{g}:=\{\,i\,|\, g_i=1\,\},\quad I_{g}:=\{1,\ldots,N\}\setminus I^{g},\quad d_g:=|I_g|.
\]
We denote by $\K[X^g]$ the quotient of $\K[X]$ by the ideal generated by $x_i$ with $i\in I_{g}$ and by ${\sf res}^g$  the projection $\K[X]\to\K[X^g]$. We write $x^g$ instead of ${\sf res}^g(x)$, i.e. $x^g_i=x_i$ if  $i\in I^{g}$ and $0$ if  $i\in I_{g}$.

\subsection{The main theorem}\label{section3.2}
Let $W^g:={\sf res}^g(W)$ and $M(W^g):=\K[X^g]/(\partial_{x_i}W^g)_{i\in I^g}$.
Consider the following two $(\Z\times G)$-graded $\K\bb$-modules:
\begin{equation}\label{jacom}
\Jaco^*(X,W,G):=\bigoplus_{g\in G} M(W^g)\bb\cdot \xi_g,\quad \Omeg^*(X,W,G):=\bigoplus_{g\in G} M(W^g)\bb\cdot \omega_g
\end{equation}
where
the elements of $M(W^g)$ have degree 0, $\xi_g$ is a formal generator of degree $d_g$, and $\omega_g$ is a formal generator of degree $d_g-N$. We equip (\ref{jacom}) with degree preserving $\K\bb$-linear $G$-actions by requiring that $G$ act on $M(W^g)$ in the natural way and on $\xi_g$ and $\omega_g$ as follows: 
\begin{equation}\label{Gvw} 
G\ni h=(h_1, \ldots, h_{N}): \quad \xi_g\mapsto \prod_{i\in I_g} h^{-1}_i\,\cdot  \xi_g,\quad \omega_g\mapsto \prod_{i\in I^g} h_i\,\cdot\omega_g.
\end{equation}

\noindent Furthermore, for $g,h\in G$ we define $\sigma_{g,h}\in M(W^{gh})$ as the coefficient at $\partial_{\theta_{I_{gh}}}$ in the expression
\begin{equation}\label{cgh}
\frac1{d_{g,h}!}\,\LL\left(\left(\left\lfloor{\rm H}_W(x,g(x),x)\right\rfloor_{gh}+\lfloor{\rm H}_{{W,g}}(x)\rfloor_{gh}\otimes1+1\otimes \lfloor {\rm H}_{{W,h}}(g(x))\rfloor_{gh}\right)^{d_{g,h}}\otimes\partial_{\theta_{I_g}}\otimes\partial_{\theta_{I_h}}\right)
\end{equation} 
where 

\vspace{0.1in}
\noindent (1) 
 ${\rm H}_W(x,g(x),x)$ is the degree $-2$ element of $\K[X]\otimes \K[\theta]^{\otimes2}$ defined as the restriction to the set $\{y=g(x),\, z=x\}$ of the following degree $-2$ element of  $\K[X]^{\otimes 3}\otimes \K[\theta]^{\otimes 2}$ 
\begin{equation}\label{del-2}
{\rm H}_{W}={\rm H}_{W}(x,y,z):=\sum_{1\leq j\leq i\leq N} \nabla^{y\to(y,z)}_j\nabla^{x\to(x,y)}_i(W)\,\theta_i\otimes \theta_j;
\end{equation}
  
\vspace{0.1in}
\noindent (2) 
 ${\rm H}_{{W,g}}\in\K[X]\otimes \K[\theta]$ is the degree $-2$ element of $\K[X][\theta]$ given by the formula 
\begin{eqnarray}\label{hg}
{\rm H}_{{W,g}}= {\rm H}_{{W,g}}(x):=\sum_{{i,j\in I_{g},\,\,  j<i}}\frac{1}{1-g_j}\nabla^{x\to(x,x^g)}_j\nabla^{x\to(x,g(x))}_i(W)\,\theta_j\,\theta_i;
\end{eqnarray}

\vspace{0.1in}
\noindent (3) 
$\lfloor f \rfloor_{g}$ for $f\in\K[X]$ denotes the class of ${\sf res}^{g}(f)\in\K[X^g]$  in $M(W^{g})$; we extend this operator to an operator $\K[X]\otimes V\to M(W^g)\otimes V$ by $V$-linearity; 

\vspace{0.15in}
\noindent (4) 
$d_{g,h}:=\frac{d_g+d_h-d_{gh}}{2}$ and the $d_{g,h}$-th power  in (\ref{cgh}) is computed with respect to the natural product on  $\K[{X}]\otimes \K[\theta]\otimes \K[\theta]$;   
we set
\begin{equation}\label{non-int}
\text{\it $\sigma_{g,h}=0$ if $d_{g,h}$ is not a non-negative integer};
\end{equation}
 
\vspace{0.1in}
\noindent (5) 
 $\LL$ is the $\K[X]$-linear extension of the degree 0 map $\K[{\theta}]^{\otimes2}\otimes \K[\partial_{\theta}]^{\otimes2}\to\K[\partial_{\theta}]$ defined by
\begin{equation}\label{mu}
 p_1(\theta)\otimes p_2(\theta)\otimes q_1(\partial_\theta)\otimes q_2(\partial_\theta)\mapsto(-1)^{|q_1||p_2|}p_1(q_1)\cdot p_2(q_2)
\end{equation}
where $p_i(q_i)$ denotes the action of $p_i(\theta)$ on $q_i(\partial_\theta)$ via the ${\rm Cl}_N$-module structure on $\K[\partial_\theta]$ introduced in Section \ref{CA} and $\cdot$   is the natural product in $\K[\partial_{\theta}]$.

\begin{theorem}\label{mainthm} Assume $W$ has only {\bf isolated critical points}. Then

\noindent(1) For all $g,h\in G$ the $\K\bb$-linear maps
\begin{eqnarray}\label{cupjac}
\cup: M(W^g)\bb\cdot \xi_g\otimes_{\K\bb}M(W^h)\bb\cdot \xi_h&\to& M(W^{gh})\bb\cdot \xi_{gh},\nonumber\\ \lfloor f_1\rfloor_g\cdot \xi_g\otimes \lfloor f_2\rfloor_h\cdot \xi_h&\mapsto& 
t^{d_{g,h}}\cdot \lfloor f_1f_2\rfloor_{gh}\cdot \sigma_{g,h}\cdot \xi_{gh},
\end{eqnarray}
\begin{eqnarray*}
\cap: M(W^g)\bb\cdot \omega_g\otimes_{\K\bb}M(W^h)\bb\cdot \xi_h&\to& M(W^{gh})\bb\cdot \omega_{gh},\nonumber\\
\lfloor f_1\rfloor_g\cdot \omega_g\otimes \lfloor f_2\rfloor_h\cdot \xi_h&\mapsto& t^{d_{g,h}}\cdot \lfloor f_1f_2\rfloor_{gh}\cdot \sigma_{g,h}\cdot \omega_{gh}\\
(\forall f_1,f_2\in\K[X])\nonumber
\end{eqnarray*}
are well-defined and make $\Jaco^*(X,W,G)$ into a $(\Z\times G)$-graded $\K\bb$-linear $G$-equivariant braided super-commutative  associative algebra {\rm(}with unit $\xi_e${\rm)} and $\Omeg^*(X,W,G)$ into a $(\Z\times G)$-graded $\K\bb$-linear $G$-equivariant free rank 1 right  $\Jaco^*(X,W,G)$-module {\rm(}with generator $\omega_e${\rm)}.

\noindent(2)
There is a $G$-equivariant isomorphism 
\begin{equation*}
\text{\scalebox{0.95}{$
\left(\cHH^*(\K[X],W;\K[X]\rtimes  G), \cup, \cHH_*(\K[X],W;\K[X]\rtimes  G), \cap\,\right)\simeq \left(\Jaco^*(X,W,G), \cup, \Omeg^*(X,W,G), \cap\,\right).
$}}
\end{equation*}
\end{theorem}

\noindent Applying Proposition \ref{Ginv}, we obtain
\begin{corollary}\label{maincor}
\[
\text{\scalebox{0.99}{$
\left(\cHH^*(\K[X]\rtimes  G,W), \cup, \cHH_*(\K[X]\rtimes  G,W), \cap\,\right)\simeq \left(\Jaco^*(X,W,G)^G, \cup, \Omeg^*(X,W,G)_G, \cap\,\right).
$}}
\]
\end{corollary}

Let us formulate a variant of the above results for the Hochschild homology of the second kind (see Remark \ref{HHtHHII}).
Let $\Jacoii(X,W,G)$ and $\Omegii(X,W,G)$ stand for the algebra and the module over the algebra  defined just as $\Jaco^*(X,W,G)$ and $\Omeg^*(X,W,G)$ but with $\K\bb$ replaced  by $\K$ and the $\Z$-grading reduced to a $\ZZ$-grading. 
\begin{corollary}\label{maincorii} 
There is a $G$-equivariant  isomorphism 
\begin{equation*}\text{\scalebox{0.98}{$
\left(\ccHH^{\rm II,*}(\K[X],W;\K[X]\rtimes  G), \cup, \ccHH^{\rm II}_*(\K[X],W;\K[X]\rtimes  G), \cap\,\right)\simeq \left(\Jacoii(X,W,G), \cup, \Omegii(X,W,G), \cap\,\right).
$}}
\end{equation*}
Consequently,
\[\text{\scalebox{0.98}{$
\left(\ccHH^{\rm II,*}(\K[X]\rtimes  G,W), \cup, \ccHH^{\rm II}_*(\K[X]\rtimes  G,W), \cap\,\right)\simeq \left(\Jacoii(X,W,G)^G, \cup, \Omegii(X,W,G)_G, \cap\,\right).
$}}
\]
\end{corollary}
%%%%%%%%%%

\subsection{Implications for equivariant matrix factorization categories}\label{3.3} Associated with any curved algebra $(A,W)$ there is a differential $\ZZ$-graded (dg) category  of {\it curved dg $(A,W)$-modules} \cite{PP} whose objects are the pairs $(E, \delta_E)$ where $E=E^{\rm even}\oplus E^{\rm odd}$ is a $\ZZ$-graded finitely generated projective $A$-module and $\delta_E: E\to E$ is an odd morphism  satisfying $\delta_E^2=W\cdot{\rm id}_E$.   When $(A,W)=(\K[X]\rtimes G,W)$, the category is precisely what we denoted by ${\rm MF}_G(X,W)$ in the Introduction. 

The aim of this short section is to present a HKR like theorem for the Hochschild invariants  ${\rm MF}_G(X,W)$.
Although we do not discuss the Hochschild theory of dg categories in this paper, this material is much more standard than in the curved case and is much better represented in the literature; in particular, \cite{PP} can also serve as a reference.

\begin{theorem}\label{mf} 
Assume, in addition, that the only critical value of $W$ is 0.  Then
\[
\left(\ccHH^{*}({\rm MF}_G(X,W)), \cup, \ccHH_*({\rm MF}_G(X,W)), \cap\,\right)\simeq \left(\Jacoii(X,W,G)^G, \cup, \Omegii(X,W,G)_G, \cap\,\right).
\]
\end{theorem}
\noindent As will be  explained in  Section \ref{MFC}, this theorem is essentially a reformulation  
  of Corollary \ref{maincorii}.  
 
\section{Proofs}
\subsection{A multiplicative HKR isomorphism for the Hochschild calculus of LG models}\label{sect4}
\subsubsection{Outline} The aim of this section is to establish an isomorphism of two calculi associated with $(\K[X],W)$: 
\[
\left(\cHH^*(\K[X],W;-), \cup, \cHH_*(\K[X],W;-), \cap\,\right)\simeq
\left(\cKo^*(\K[X],W;-), \cup, \cKo_*(\K[X],W;-), \cap\,\right)
\] 
where $\cKo$ is the periodic cohomology of what we call {\it Koszul mixed complexes} and  $\cup$ and $\cap$ are certain explicit products on $\cKo$.

In this section, $W$ is arbitrary (e.g. $W=0$ is allowed). We use freely the notation and conventions  established Section \ref{DD}, \ref{CA}, \ref{GA}.

\subsubsection{Koszul  mixed complexes}\label{kmc}

The Koszul resolution of $\K[{X}]$ is the complex $(\rK_*(\K[{X}]), \ddko)$ of $\K[{X}]$-bimodules with $\rK_*(\K[{X}])=\K[{X}]^{\otimes 2}[\theta]$ and  
\begin{eqnarray*}
\ddko:=\sum_{i=1}^N(x_i-y_i)\cdot\partial_{\theta_i}.
\end{eqnarray*}
Associated with $W$ there is a degree $-1$ differential on $\rK_*(\K[{X}])$, namely 
\begin{eqnarray*}
\ddcurve:=\sum_{i=1}^N\nabla_i(W)\cdot\theta_i,
\end{eqnarray*}
and it follows from the commutation relations in ${\rm Cl}_N$ and (\ref{diffder}) that
\[
\delta_{\rm curv}^2=0,\quad \delta_{\rm Kos}\delta_{\rm curv}+\delta_{\rm curv}\delta_{\rm Kos}=W(x)-W(y).
\]

Let $M$ be a $(\K[{X}],W)$-bimodule. We set
\[
\rK^*(\K[{X}], M):=\Hom_{\K[{X}]^{\otimes 2}}(\rK_{-*}(\K[{X}]), M), \quad \rK_*(\K[{X}], M):=M\otimes_{\K[{X}]^{\otimes 2}} \rK_*(\K[{X}]).
\]
As in the abstract context,  $M$ gives rise to two mixed complexes 
\begin{eqnarray*}
&&\rKo^*(\K[{X}],W; M)=\left(\rK^*(\K[{X}], M), \,\,\ddk:=\delta^\vee_{\rm Kos}, \,\,\ddw:=\delta^\vee_{\rm curv}\right),\\
&&\rKo_*(\K[{X}],W; M)=\left(\rK_*(\K[{X}], M), \,\,\bbk:=1\otimes\delta_{\rm Kos},\,\,\bbw:=1\otimes \delta_{\rm curv}\right).
\end{eqnarray*}
Their periodic cohomology will be denoted by $\cKo^*(\K[{X}],W; M)$ and $\cKo_*(\K[{X}],W; M)$.

\subsubsection{Koszul vs Hochschild complexes} Consider the morphism of $\K[X]$-bimodules
\[
\Psi: \rB_*(\K[{X}])\to \rK_*(\K[{X}])
\] defined by $\Psi(f_0[]f_{1})=l_1(f_0)l_{N+1}(f_{1})=f_0(x)f_1(y)$ and
\begin{multline}\label{psi}
\Psi(f_0[f_1|\ldots |f_n]f_{n+1})=\\=l_1(f_0)\left[\sum_{1\leq j_1<\ldots< j_n\leq N}\nabla_{j_1}(f_1)\theta_{j_1}\nabla_{j_2}(f_2)\theta_{j_2}\ldots \nabla_{j_n}(f_n)\theta_{j_n}\right]l_{N+1}(f_{n+1})
\end{multline}
for $n\geq1$ where the product on the right-hand side is taken in $\K[X]^{\otimes 2}[\theta]$.
\begin{proposition}\label{main}
One has
\begin{eqnarray}\label{morph}
\Psi\cdot\ddbar=\ddko\cdot\Psi,\quad \Psi\cdot\ddcurve=\ddcurve\cdot\Psi.
\end{eqnarray}
Moreover, $\Psi: (\rB_*(\K[{X}]),\ddbar)\to (\rK_*(\K[{X}]),\ddko)$ is a quasi-isomorphism.
\end{proposition}
\noindent{\bf Proof} is given in Appendix \ref{B} (page \pageref{mainproof}).

\begin{remark}{\rm After having discovered the formula (\ref{psi}), we looked for similar results in the literature and found them. We suspect that for $\K[X]=\K[x_1,\ldots,x_N]$ the map (\ref{psi}) coincides with the map $\Psi_B$ from \cite[Sect.4]{SW1}. The description given in {\it loc.cit.} is quite combinatorial and we will not claim this as fact. 
}
\end{remark}

The previous proposition and the fact that both $(\rK_*(\K[{X}]),\ddko)$ and $(\rB_*(\K[{X}]),\ddbar)$ are K-projective complexes of $\K[X]$-bimodules, and hence K-flat \cite[Sect.10.12]{BL}, imply

\begin{corollary}\label{mainl} For any $(\K[X],W)$-bimodule $M$ the morphisms of mixed complexes
\begin{eqnarray*}
\Psi^*:=\Psi^\vee: \rKo^*(\K[{X}],W; M)\to \rHH^*(\K[{X}],W; M)
\end{eqnarray*}
{\rm(}$\Psi^\vee$ stands for the dual of $\Psi${\rm)} and
\begin{eqnarray*}
\Psi_*:=1\otimes\Psi: \rHH_*(\K[{X}],W; M)\to \rKo_*(\K[{X}],W; M)
\end{eqnarray*}
are quasi-isomorphisms. 
\end{corollary}

\subsubsection{The cup and cap products on Koszul complexes}\label{cuap} Our next goal is to ``transfer'' the cup and cap products to the Koszul mixed complexes. 
The naive idea that there exists a morphism of bimodules $\Delta: \rK_*(\K[{X}])\to \rK_*(\K[{X}])\otimes_{\K[X]} \rK_*(\K[{X}])$ --  an analog of (\ref{coalg0}) -- that is compatible with $\ddko$ and $\ddcurve$ and  matches (\ref{coalg0}) under the quasi-isomorphism $\Psi$ from the preceding section does not work. Calculations show that such a map does not exist even in the one-dimensional case ($N=1$). As we will see in this and the next sections, the situation is more complicated.

\vspace{0.1in}
Consider the sequence of morphisms of $\K[X]$-bimodules 
\[
\Delta_{-2l}: \K[{X}]^{\otimes 2}[\theta]\to \K[X]^{\otimes 3}\otimes\K[\theta]^{\otimes 2},\quad l=0,\ldots, N
\] 
defined in terms of the $(x,y,z)$-coordinates (see Section \ref{DD}) as follows:
\begin{equation}\label{del0}
\Delta_{-2l}: f(x,y)\cdot p(\theta_1,\ldots,\theta_N)\mapsto \frac1{l!}f(x,z)\cdot {\rm H}_{W}(x,y,z)^{l}\cdot p(\theta_1\otimes1+ 1\otimes\theta_1,\ldots,\theta_N\otimes1+ 1\otimes\theta_N)
\end{equation}
where ${\rm H}_{W}$ is the  element (\ref{del-2}). 

Note that  we can also view these maps as maps to $\K[{X}]^{\otimes 2}[\theta]\otimes_{\K[X]}\K[X]^{\otimes 2}[\theta]$ where we identify the latter with $\K[X]^{\otimes 3}\otimes\K[\theta]^{\otimes 2}$ via 
\[
f_1(x,y)p_1(\theta)\otimes f_2(y,z)p_2(\theta)\mapsto f_1(x,y)f_2(y,z)\otimes p_1(\theta)\otimes p_2(\theta).
\]
Keeping this in mind, one has

\begin{proposition}\label{coprod} 
The morphism of $\K[X]\bb$-bimodules
\[
\Delta=\Delta_{\rm Kos}: \K[{X}]^{\otimes 2}\bb[\theta]\to \K[{X}]^{\otimes 2}\bb[\theta]\otimes_{\K[X]\bb}\K[X]^{\otimes 2}\bb[\theta]
\]
 given by the formula
\[
\Delta=\sum_{l=0}^N \Delta_{-2l} \uu^l=e^{t{\rm H}_{W}}\cdot \Delta_0
\]
is compatible with the operator $(\ddko+\uu\ddcurve)$: 
\[
\Delta\cdot(\ddko+\uu\ddcurve)=((\ddko+\uu\ddcurve)\otimes1+1\otimes(\ddko+
\uu\ddcurve))\cdot\Delta.
\]
\end{proposition}
\noindent{\bf Proof} is given in Appendix \ref{B} (page \pageref{coprodproof}). 

\vspace{0.2in}

It is the above operator $\Delta_{\rm Kos}$ that will be playing the role of $\Delta_{\rm bar}$ (\ref{coalg0}) in the ``Koszul theory''. Just as in the abstract setting, $\Delta_{\rm Kos}$ gives rise to cup and cap products which one defines exactly as in   (\ref{cup}) and (\ref{int}): For two $(\K[X],W)$-bimodules $M_1$ and $M_2$ 
\begin{equation*}
\cup\,=\,\cup_{\rm Kos}: \rK^*(\K[{X}], M_1)\bb\otimes_{\K\bb} \rK^*(\K[{X}], M_2)\bb
\to \rK^*(\K[{X}], M_1\otimes_{\K[X]}M_2)\bb,
\end{equation*}
is the map determined by 
\begin{equation}\label{caup1}
\forall \omega\in\rK_*(\K[X]): \langle\omega, \xi_1\cup_{\rm Kos} \xi_2\rangle:=\langle \Delta_{\rm Kos}(\omega), \xi_1\boxtimes \xi_2\rangle=(-1)^{|\xi_1||\omega_{(2)}|}\langle \omega_{(1)}, \xi_1\rangle\otimes \langle \omega_{(2)}, \xi_2\rangle
\end{equation}
and
\begin{equation*}
\cap\,=\,\cap_{\rm Kos}: \rK_*(\K[{X}], M_1)\bb\otimes_{\K\bb} \rK^*(\K[{X}], M_2)\bb\to \rK_*(\K[{X}],M_1\otimes_{\K[X]}M_2)\bb,
\end{equation*}
is the map defined by
\begin{equation*}
(m\otimes \omega)\cap_{\rm Kos} \xi:=m\otimes\langle \Delta_{\rm Kos}(\omega), \xi\boxtimes{\rm id}_{\rK_*(\K[{X}])}\rangle=(-1)^{|\xi||\omega_{(2)}|} (m\otimes \langle \omega_{(1)},\xi\rangle)\otimes \omega_{(2)}.
\end{equation*}
In the above formulas, $\langle \cdot,\cdot\rangle$ and $\boxtimes$ have the same meaning as before, that is, as in Remark \ref{lr} and in (\ref{sqten}), respectively. (Formally speaking, $\langle \cdot,\cdot\rangle$ above is a $\K\bb$-linear extension of the previous definition.)

Proposition \ref{coprod} yields 
\begin{corollary}\label{colcoprod} The maps $\cup_{\rm Kos}$ and $\cap_{\rm Kos}$  
descend to products on the periodic cohomology
\[
\cKo^*(\K[{X}],W; M_1)\otimes_{\K\bb} \cKo^*(\K[{X}],W; M_2)\to \cKo^*(\K[{X}],W; M_1\otimes_{\K[X]}M_2),
\]
\[
\cKo_*(\K[{X}],W; M_1)\otimes_{\K\bb} \cKo^*(\K[{X}],W; M_2)\to \cKo_*(\K[{X}],W; M_1\otimes_{\K[X]}M_2).
\]

\end{corollary}

\begin{remark}{\rm Note that, unlike  $\cup_{\rm Hoch}$ and $\cap_{\rm Hoch}$, the products $\cup_{\rm Kos}$ and $\cap_{\rm Kos}$ are {\it not associative} on the cochain level. However, as we will see shortly, the induced products on the periodic cohomology are associative.
}
\end{remark}
\subsubsection{Comparing the products on the Hochschild and Koszul complexes}

The diagram
\[\begin{CD}
\rB_{*}(\K[X])\bbb @>{\Delta_{\rm bar}}>>\left(\rB_{*}(\K[X])\otimes_{\K[X]} \rB_{*}(\K[X])\right)\bbb\\
@VV{\Psi}V @VV\Psi\otimes\Psi V\\
\rK_*(\K[X])\bb @>\Delta_{\rm Kos}>> \left(\rK_*(\K[X])\otimes_{\K[X]}\rK_*(\K[X])\right)\bb
\end{CD}\]

\vspace{0.1in}
\noindent is not commutative (e.~g. because only $\Delta_{\rm Kos}$ ``depends'' on $\uu$).
Nevertheless, one has

\begin{proposition}\label{comdiag}
The diagram is commutative up to homotopy: There exist morphisms of $\K[X]$-bimodules \[
h_i: \rB_{*}(\K[X])\to \rK_*(\K[X])\otimes_{\K[X]}\rK_*(\K[X]), \quad |h_i|=-2i-1, \quad i=0,1,\ldots
\]
such that 
\begin{multline}\label{2}
(\Psi\otimes\Psi)\cdot \Delta_{\rm bar}-\Delta_{\rm Kos}\cdot\Psi=\\
=((\ddko+t\ddcurve)\otimes1+1\otimes(\ddko+t\ddcurve))\cdot\sum_{i=0}^\infty h_i t^i+\sum_{i=0}^\infty h_i t^i\cdot (\ddbar+t\ddcurve).
\end{multline}
\end{proposition}
\noindent{\bf Proof} is given in Appendix \ref{B} (page \pageref{comdiagproof}).  

\begin{corollary} The isomorphisms
\begin{eqnarray*}
&&\Psi^*: \cKo^*(\K[{X}],W; M)\to \cHH^*(\K[{X}],W; M),\\ &&\Psi_*:\cHH_*(\K[{X}],W; M)\to \cKo_*(\K[{X}],W; M),
\end{eqnarray*}
induced by the quasi-isomorphisms from Corollary \ref{mainl}, are compatible with the cup and  cap products:
\[
\Psi^*(\xi_1\cup_{\rm Kos} \xi_2)=\Psi^*(\xi_1)\cup_{\rm Hoch} \Psi^*(\xi_2),\quad \forall\, \xi_i\in \cKo^*(\K[{X}],W; M_i),
\]
\[
\Psi_*(\omega)\cap_{\rm Kos} \xi=\Psi_*(\omega\cap_{\rm Hoch}\Psi^*(\xi)) \quad \forall \, \omega\in \cHH_*(\K[{X}],W; M_1), \xi\in \cKo^*(\K[{X}],W; M_2).
\]
\end{corollary}

\subsubsection{A more explicit description of the Koszul calculus}\label{comasp}  Our aim in this section is to rewrite all the structures in the Koszul calculus in a slightly more explicit way. We start by working out a more convenient description of the mixed Koszul complexes.

Let  $\star: {\rm Cl}_N\to {\rm Cl}_N$ be the algebra involution (i.e., a $\K$-linear degree 0 involutive map satisfying $(\xi_1\xi_2)^\star=(-1)^{|\xi_1||\xi_2|}\xi_2^\star \xi_1^\star$) determined by
\[
\theta_i^\star=\theta_i,\quad \partial_{\theta_i}^\star=-\partial_{\theta_i},\quad \forall\, i=1,\ldots, N.
\]
Consider the following perfect pairing between $\K[\theta]$ and 
$\K[\partial_{\theta}]$:
\begin{equation}\label{pair1}
\{\cdot,\cdot\}: \K[\theta]\otimes\K[\partial_{\theta}]\to\K, \quad \{ p(\theta),q(\partial_{\theta})\}:=(-1)^{|p(\theta)||q(\partial_{\theta})|}\ {\rm CT}\left(q(\partial_{\theta})^\star (p(\theta))\right)
\end{equation}
where on the right-hand side we apply the ``differential operator'' $q(\partial_{\theta})^\star$ to the polynomial $p(\theta)$ and take the constant term of the resulting polynomial. One can easily check that
\begin{equation}\label{dual}
\{ \xi(p(\theta)),q(\partial_{\theta})\}=(-1)^{|p(\theta)||\xi|} \{ p(\theta) ,\xi^\star (q(\partial_{\theta}))\}\quad \forall \xi \in{\rm Cl}_N, \, p(\theta)\in\K[\theta],\, q(\partial_{\theta})\in \K[\partial_{\theta}]
\end{equation}
where  $\xi$ acts on $\K[\theta]$ and $\K[\partial_{\theta}]$ as described in Section \ref{CA}.

Let $M$ be a $(\K[{X}],W)$-bimodule. The above pairing gives rise to the pairing
\begin{equation}\label{pair2}
\rK_*(\K[X])\otimes M[\partial_{\theta}]\to M,\quad \{ f_0(x)f_1(y)p(\theta), mq(\partial_{\theta}) \}:=(f_0mf_1)\{ p(\theta),q(\partial_{\theta}) \}
\end{equation}
which induces  an isomorphism  
\begin{equation}\label{TM}
M[\partial_{\theta}]\to\rK^*(\K[{X}], M).
\end{equation}
Under this isomorphism the differentials $\ddk=\delta^\vee_{\rm Kos}$ and $\ddw=\delta^\vee_{\rm curv}$ on the right-hand side correspond to the differentials (denoted by the same symbols)
\begin{eqnarray*}
\ddk:=\sum_{i=1}^N(x_i-y_i)\cdot\partial_{\theta_i},\quad \ddw:=-\sum_{i=1}^N\nabla^{x\to(x,y)}_i(W)\cdot\theta_i
\end{eqnarray*}
on $M[\partial_{\theta}]$ where the $x$-variables act on $M$ from the left and the $y$-variables act from the right.  Thus,
\begin{equation}\label{tm}
\rKo^*(\K[{X}],W; M)\simeq \left(M[\partial_{\theta}], \ddk, \ddw\right).
\end{equation}

There is a similar description of  $\rKo_*(\K[{X}],W; M)$. Namely, observe that 
\begin{equation}\label{OM}
\rK_*(\K[{X}], M)\to M[\theta],\quad m\otimes f_0(x)f_1(y)p(\theta)\mapsto (f_1m f_0)p(\theta)
\end{equation}
is an isomorphism. It transforms the differentials $\bbk=1\otimes\delta_{\rm Kos}$ and $\bbw=1\otimes \delta_{\rm curv}$ on the left-hand side into the differentials 
\begin{eqnarray*}
\bbk:=\sum_{i=1}^N(x_i-y_i)\cdot\partial_{\theta_i},\quad \bbw:=\sum_{i=1}^N\nabla^{x\to(x,y)}_i(W)\cdot\theta_i
\end{eqnarray*}
on $M[\theta]$ where now the $x$-variables act on $M$ from the {\it right} and the $y$-variables act from the {\it left}.  Thus,
\begin{equation}\label{om}
\rKo_*(\K[{X}],W; M)\simeq\left(M[\theta], \bbk, \bbw\right).
\end{equation}

\vspace{0.1in}
 
Our next goal is to describe $\cup_{\rm Kos}$ and $\cap_{\rm Kos}$ in  terms of the mixed complexes on the right-hand side of (\ref{tm}) and (\ref{om}). In what follows, $M_1$ and $M_2$ are two $(\K[X],W)$-bimodules. We start with the cup product.

Take arbitrary two elements $\xi_i\in M_i[\partial_{\theta}]$ ($i=1,2$) and let $\widehat{\xi}_i\in \rK_*(\K[X]), M_i)$ denote their images under (\ref{TM}). Take also any $p(\theta)\in\K[\theta]\subset \K[X]^{\otimes 2}[\theta]=\rK_*(\K[X])$. Then, by our definitions (see (\ref{caup1})) 
\begin{equation*}
\{ p(\theta), {\xi}_1\cup_{\rm Kos} {\xi}_2\}=\langle p(\theta), \widehat{\xi}_1\cup_{\rm Kos} \widehat{\xi}_2\rangle= \langle \Delta_{\rm Kos}(p(\theta)), \widehat{\xi}_1\boxtimes \widehat{\xi}_2\rangle.
\end{equation*}
The pairing  
$\langle\cdot,\cdot\rangle$ in the  last term  is the $\K\bb$-linear extension of the pairing    
\[
\text{\scalebox{0.95}{$
\rK_*(\K[X])\otimes_{\K[X]}\rK_*(\K[X])\otimes \Hom_{\K[X]^{\otimes 2}}(\rK_*(\K[X])\otimes_{\K[X]}\rK_*(\K[X]), M_1\otimes_{\K[X]}M_2)\to M_1\otimes_{\K[X]}M_2.
$}}
\] 
Under the isomorphism 
\[
\Hom_{\K[X]^{\otimes 2}}(\rK_*(\K[X])\otimes_{\K[X]}\rK_*(\K[X]), M_1\otimes_{\K[X]}M_2) \simeq M_1[\partial_{\theta}]\otimes_{\K[X]} M_2[\partial_{\theta}]
\] 
$\widehat{\xi}_1\boxtimes \widehat{\xi}_2$ is just $\xi_1\otimes \xi_2$ and the above pairing  is 
nothing but the pairing (\ref{pair1}) for ${\rm Cl}_{N}\otimes {\rm Cl}_{N}={\rm Cl}_{2N}$ (or, rather, its extension  analogous to (\ref{pair2})). Let us denote it by $\{\!\{\cdot,\cdot\}\!\}$.
Putting everything together, we obtain
\[
\{p(\theta), {\xi}_1\cup_{\rm Kos} {\xi}_2\}=\{\!\{e^{t{\rm H}_W}\Delta_0(p(\theta)), \xi_1\otimes \xi_2\}\!\}.
\]
Furthermore, 
\[
\{\!\{e^{t{\rm H}_W}\Delta_0(p(\theta)), \xi_1\otimes \xi_2\}\!\}=\{\!\{\Delta_0(p(\theta)), e^{t{\rm H}^\star_W}(\xi_1\otimes \xi_2)\}\!\}=\{\!\{\Delta_0(p(\theta)), e^{t{\rm H}_W}(\xi_1\otimes \xi_2)\}\!\}
\]
where we use that ${\rm H}_W$ is even (so no signs pop up upon  using (\ref{dual})) and also  that ${\rm H}^\star_W={\rm H}_W$.
To conclude the calculation of the cup product, we observe that
\begin{equation*}
\{\!\{\Delta_0(p(\theta)), q_1(\partial_{\theta})\otimes q_2(\partial_{\theta})\}\!\}=
\{ p(\theta), q_1(\partial_{\theta})q_2(\partial_{\theta})\}\quad \forall\,p,q_1,q_2
\end{equation*}
where on the right-hand side we multiply  $q_1(\partial_{\theta})$ and $q_2(\partial_{\theta})$ just as elements of $\K[\partial_{\theta}]$. The formula is easily verified  by substituting monomials for $p$, $q_1$ and $q_2$.

Summarizing the above calculations, we obtain the following result. Let us write symbolically 
\[
e^{t{\rm H}_W}=\sum_{h_i,p_i,l} (h_1\otimes h_2\otimes h_3)(p_1\otimes p_2) t^l,\quad h_i\in\K[X],\, p_i\in\K[\theta].
\]
Then for $m_1q_1(\partial_{\theta})\in M_1[\partial_{\theta}]$ and  $m_2q_2(\partial_{\theta})\in M_2[\partial_{\theta}]$
\begin{equation}\label{formulacup}
m_1q_1(\partial_{\theta})\cup m_2q_2(\partial_{\theta})
=\LL(e^{t{\rm H}_W}\otimes m_1q_1(\partial_{\theta})\otimes m_2q_2(\partial_{\theta}))
\end{equation}
where the right-hand side is an element in $(M_1\otimes_{\K[X]}M_2)[\partial_{\theta}]\bb$ given by the formula
\begin{equation*}
\LL(e^{t{\rm H}_W}\otimes m_1q_1\otimes m_2q_2)=
\sum(h_1m_1h_2\otimes m_2h_3)\,\LL(p_1\otimes p_2\otimes q_1\otimes q_2)\,t^l
\end{equation*}
with $\LL$ on the right-hand side being the map  (\ref{mu}).

There is a similar formula for $\cap_{\rm Kos}$, namely: 
\begin{equation}\label{formulacap}
m_1p({\theta})\cap m_2q(\partial_{\theta})
=\LLL(e^{t{\rm H}_W}\otimes m_1p({\theta})\otimes m_2q(\partial_{\theta}))
\end{equation}
where the right-hand side is an element in $(M_1\otimes_{\K[X]}M_2)[{\theta}]\bb$ given by the formula
\begin{equation*}
\LLL(e^{t{\rm H}_W}\otimes m_1p\otimes m_2 q)=
\sum(h_3m_1h_1\otimes m_2h_2)\,\LLL(p_1\otimes p_2\otimes p\otimes q)\,t^l
\end{equation*}
with $\LLL: \K[{\theta}]\otimes \K[{\theta}]\otimes \K[{\theta}]\otimes \K[\partial_{\theta}]\to\K[{\theta}]$ being the unique map such that  
\begin{equation}\label{nu}
\{ \LLL(p_1\otimes p_2\otimes p\otimes q), q'\}=(-1)^{|p|(|p_1|+|p_2|)}\{ p, \LL(p_1\otimes p_2\otimes q\otimes q')\}\quad \forall q'\in \K[\partial_{\theta}].
\end{equation}
The derivation of (\ref{formulacap}) is analogous to that of (\ref{formulacup}) and is left to the reader.

\vspace{0.1in}

\noindent \emph{In the remainder of this section, we will {\it identify} $\rKo^*(\K[{X}],W; M)$ with $(M[\partial_{\theta}], \ddk, \ddw)$ and $\rKo_*(\K[{X}],W; M)$ with $(M[{\theta}], \bbk, \bbw)$. }

\subsection{A multiplicative HKR isomorphism for the Hochschild (co)homology of LG orbifolds}\label{sect5}

\subsubsection{Outline} In this section, $W\in\K[X]$ is again arbitrary (e.g. 0) and we also fix an abelian group $G$ of symmetries of $(X,W)$ of the form specified in Section \ref{SN}.

The results of the previous section yield an explicit isomorphism between the tuples
\[\left(\cHH^*(\K[X],W;\K[X]\rtimes  G), \cup_{\rm Hoch}, \cHH_*(\K[X],W;\K[X]\rtimes  G), \cap_{\rm Hoch}\,\right)\]
and
\[\left(\cKo^*(\K[X],W;\K[X]\rtimes  G), \cup_{\rm Kos}, \cKo_*(\K[X],W;\K[X]\rtimes  G), \cap_{\rm Kos}\,\right).\]
Our aim in this section is to derive a more detailed description of the latter tuple  and also a description of the $G$-action on the Koszul (co)homology that corresponds to the $G$-action on the Hochschild (co)homology (Section \ref{sect3.1}) under the isomorphism.

\subsubsection{$G$-twisted Koszul (co)homology}\label{sect5.1} Let $g$ be an arbitrary element of $G$. Our first goal  is to ``calculate'' the Koszul (co)homology $\cKo^*(\K[{X}],W; \K[X]\otimes g)$ and $\cKo_*(\K[{X}],W; \K[X]\otimes g)$. We start with the  cohomology.

The isomorphism $\K[X]\otimes g\to\K[x]$, $f\otimes g\mapsto f$ induces an isomorphism of mixed complexes
\begin{equation}\label{tmg}
\left((\K[X]\otimes g) [\partial_{\theta}], \ddk, \ddw\right)\simeq \left(\K[X][\partial_{\theta}], \ddkg, \ddwg\right)
\end{equation}
where
\begin{eqnarray*}
\ddkg:=\sum_{i=1}^N(x_i-g(x_i))\,\partial_{\theta_i},\quad \ddwg:=-\sum_{i=1}^N\nabla^{x\to(x,y)}_i(W)|_{y=g(x)}\,\theta_i.
\end{eqnarray*}
Observe that  
$
\ddkg=\sum_{i\in I_{g}}(1-g_i)x_i\,\partial_{\theta_i}
$
and 
$
 \nabla^{x\to(x,y)}_i(W)|_{y=g(x)}=l_i^{x\to(x,g(x))}(\partial_{x_i}W(x))$ for $i\in I^{g}$.
The latter observation suggests one to split $\ddwg$ into two components, namely
\[
\ddkg+t\ddwg=\ddkg+t\ddwg'+t\ddwg''
\]
with
\[
\ddwg':=-\sum_{i\in I_{g}}\nabla^{x\to(x,g(x))}_i(W)\,\theta_i,\qquad \ddwg'':=-\sum_{i\in I^{g}}l_i^{x\to(x,g(x))}(\partial_{x_i}W(x))\, \theta_i.
\]

\begin{lemma}\label{ethg} Let ${\rm H}_{{W,g}}$ be as in (\ref{hg}). 
Then
\begin{eqnarray}\label{conj}
\ddkg+t\ddwg=e^{t  {\rm H}_{W,g}}\cdot\left(\ddkg+t\ddwg''\right)\cdot e^{-t{\rm H}_{{W,g}}}
\end{eqnarray}
 viewed as elements of $\K[X]\bb\otimes {\rm Cl}_N$.
\end{lemma}
\noindent{\bf Proof} is given in Appendix \ref{B} (page \pageref{ethgproof}).

\medskip

By this lemma  the map $e^{t  {\rm H}_{{W,g}}}:\K[X][\partial_{\theta}]\bb\to \K[X][\partial_{\theta}]\bb$ induces a quasi-isomorphism
\begin{equation}\label{iso}
(\K[X][\partial_{\theta}]\bb,\ddkg+t\ddwg'')\stackrel{\sim}{\to} (\K[X][\partial_{\theta}]\bb,\ddkg+t\ddwg).
\end{equation}
Let us calculate the cohomology of the  complex on the left-hand side.

Let
 $\K[\partial_\theta^g]$ denote the subalgebra in ${\rm Cl}_N$ generated by  $\{\partial_{\theta_i}\}_{i\in I^{g}}$.  Consider the projection $\K[\partial_{\theta}]\to \K[\partial_\theta^g]\cdot \partial_{\theta_{I_g}}$ which annihilates monomials containing less than $d_g$ elements $\partial_{\theta_i}$, $i\in I_g$. This projection together with the restriction  ${\sf res}^g: \K[X]\to\K[X^g]$ (Section \ref{GA})  give rise to a morphism of mixed complexes
\begin{eqnarray*}
(\K[X][\partial_{\theta}],\,\ddkg,\,\ddwg'')\to (\K[X^g][\partial_\theta^g]\cdot \partial_{\theta_{I_g}},\,0,\, \sum_{i\in I^{g}}(\partial_{x_i}W^g)\theta_i)
\end{eqnarray*}
which clearly is a quasi-isomorphism. In particular, 
\begin{multline*}
{\rm H}^*(\K[X][\partial_{\theta}]\bb,\ddkg+t\ddwg'')\simeq\\ \simeq{\rm H}^{*-d_g}(\K[X^g][\partial_\theta^g]\bb\cdot \partial_{\theta_{I_g}},\,t\sum_{i\in I^{g}}(\partial_{x_i}W^g)\theta_i)\simeq \cKo^{*-d_g}(\K[X^g],W^g).
\end{multline*}
Combining this observation with (\ref{iso}), we obtain:
\begin{proposition}\label{invsec1} As a $\Z$-graded $\K\bb$-module, 
\[\cKo^*(\K[X],W; \K[X]\otimes g)\simeq \cKo^{*-d_g}(\K[X^g],W^g)(\simeq{\rm H}^{*-d_g}(\wedge^* T_{X^g}, [W^g,\cdot])\bb).\]
\end{proposition}

\begin{remark}\label{explcoh}{\rm 
For the purpose of calculation of the cup products it is important to have explicit representatives of classes in $\cKo^*(\K[X],W; \K[X]\otimes g)$. Let us therefore formulate the above observations more carefully. 
The natural embedding $\K[X][\partial_\theta^g]\cdot \partial_{\theta_{I_g}}\hookrightarrow \K[X][\partial_{\theta}]$ induces a map
\begin{equation}\label{sur}
{\rm Ker} (t\ddwg''|_{\K[X][\partial_\theta^g]\bb\cdot \partial_{\theta_{I_g}}})\to {\rm H}^*(\K[X][\partial_{\theta}]\bb,\,\ddkg+t\ddwg'')
\end{equation}
which by the preceeding discussion is surjective. Consequently, representatives of  classes in $\cKo^*(\K[X],W; \K[X]\otimes g)$ can be obtained by applying the map (\ref{iso}) to elements of the space 
${\rm Ker} (t\ddwg''|_{\K[X][\partial_\theta^g]\bb\cdot \partial_{\theta_{I_g}}})$.  To complete this description, we need to understand what the kernel of (\ref{sur}) looks like. It is easy to describe: it contains ${\rm Im} (t\ddwg''|_{\K[X][\partial_\theta^g]\bb\cdot \partial_{\theta_{I_g}}})$, as well as  the closed elements annihilated by ${\sf res}^g$. Thus, the kernel equals
\begin{equation*}
{\rm Im} (t\ddwg''|_{\K[X][\partial_\theta^g]\bb\cdot \partial_{\theta_{I_g}}})+ {\rm Ker} (t\ddwg''|_{\K[X][\partial_\theta^g]\bb\cdot \partial_{\theta_{I_g}}})\cap \sum_{i\in I_g}x_i\cdot\K[X][\partial_{\theta}^g]\bb\cdot \partial_{\theta_{I_g}}.
\end{equation*}
}
\end{remark}
Let us derive an analogous description for $\cKo_*(\K[X],W; \K[X]\otimes g)$. Since the argument is very similar, we will only sketch it and omit details.

There is  a natural isomorphism of mixed complexes 
\begin{equation}\label{omg}
\left((\K[X]\otimes g) [{\theta}], \bbk, \bbw\right)\simeq \left(\K[X][{\theta}], \bbkg, \bbwg=\bbwg'+\bbwg''\right)
\end{equation}
where $\bbkg:=-\sum_{i\in I_{g}}(1-g_i)x_i\,\partial_{\theta_i}$,
\[
\bbwg':=\sum_{i\in I_{g}}\nabla^{x\to(g(x),x)}_i(W)\,\theta_i,\qquad \bbwg'':=\sum_{i\in I^{g}}l_i^{x\to(g(x),x)}(\partial_{x_i}W(x))\, \theta_i.
\]
By analogy with the previous case, we have a quasi-isomorphism 
\begin{equation}\label{expiso}
e^{t  {\rm H}_{{W,g}}^\dagger}: (\K[X][{\theta}]\bb,\bbkg+t\bbwg'')\stackrel{\sim}{\to} (\K[X][{\theta}]\bb,\bbkg+t\bbwg)
\end{equation}
where this time
\begin{eqnarray*}
{\rm H}_{{W,g}}^\dagger={\rm H}_{{W,g}}^\dagger(x):=\sum_{{i,j\in I_{g},\,\,  j<i}}\frac{1}{1-g_j}\nabla^{x\to(x,x^g)}_j\nabla^{x\to(g(x),x)}_i(W)\,\theta_j\,\theta_i.
\end{eqnarray*}

Let
$\K[\theta^g]$ stand for the subalgebra in ${\rm Cl}_N$ generated by  $\{{\theta_i}\}_{i\in I^{g}}$.  Consider the projection $\K[{\theta}]\to \K[\theta^g]$ that annihilates monomials containing  $\theta_i$, $i\in I_g$. This projection and the homomorphism ${\sf res}^g$ give rise to a quasi-isomorphism of mixed complexes
\begin{eqnarray*}
(\K[X][{\theta}],\,\bbkg,\,\bbwg'')\to (\K[X^g][\theta^g],\,0,\, \sum_{i\in I^{g}}(\partial_{x_i}W^g)\theta_i),
\end{eqnarray*}
which we combine with (\ref{expiso}) to obtain
\begin{proposition}\label{invsec2} As a $\Z$-graded $\K\bb$-module,
\[\cKo_*(\K[X],W; \K[X]\otimes g)\simeq \cKo_{*}(\K[X^g],W^g)(\simeq {\rm H}^{*}(\Omega^*_{X^g}, dW^g\wedge\cdot)\bb).\] 
\end{proposition}

\begin{remark}\label{explhom}{\rm
Note that the embedding $\K[X][\theta^g]\hookrightarrow \K[X][{\theta}]$ induces a surjective map 
\begin{equation}\label{sur1}
{\rm Ker} (t\bbwg''|_{\K[X][\theta^g]\bb})\to {\rm H}^*(\K[X][{\theta}]\bb,\,\bbkg+t\bbwg'')
\end{equation}
and so explicit representatives of classes in $\cKo_*(\K[X],W; \K[X]\otimes g)$ can be obtained by applying (\ref{expiso}) to elements of  ${\rm Ker} (t\bbwg''|_{\K[X][\theta^g]\bb})$.
The kernel of (\ref{sur1}) equals
\[
{\rm Im} (t\bbwg''|_{\K[X][\theta^g]\bb})+  {\rm Ker} (t\bbwg''|_{\K[X][\theta^g]\bb})\cap \sum_{i\in I_g}x_i\K[X][{\theta}^g]\bb.
\]
}
\end{remark}

\begin{remark}\label{nonab1}{\rm Part of the above picture can be generalized to not necessarily abelian subgroups $G\subset GL_N(\K)$. Namely, even if $g$ is not diagonal but {\it can} be diagonalized (say when $\K=\C$), Propositions \ref{invsec1} and \ref{invsec2} hold true for the simple reason that neither the Hochschild (co)homology $\cHH^*(\K[X],W; \K[X]\otimes g)$, $\cHH_*(\K[X],W; \K[X]\otimes g)$ nor the geometrically  defined cohomology ${\rm H}^{*}(\wedge^* T_{X^g}, [W^g,\cdot])$ and ${\rm H}^{*}(\Omega^*_{X^g}, dW^g\wedge\cdot)$ depend on any coordinate systems. Furthermore, even though the explicit description of the cohomology classes that we have derived in Remarks \ref{explcoh}, \ref{explhom} does not apply to non-diagonal elements $g$, the underlying idea -- namely, splitting the differentials into two parts and using exponential twists to simplify complexes -- seems quite universal. In order to calculate the $g$-twisted Koszul (co)homology in a concrete example, one can write  the  differentials $\ddkg$, $\ddwg$, $\bbkg$, $\bbwg$ in linear coordinates  in which $g$ is diagonal, apply the above idea to do calculations, and write the result in terms of the original coordinates.  
}
\end{remark}

\subsubsection{Products}\label{prodkos1} Our next goal is to describe the products 
\begin{eqnarray*}
\cKo^*(\K[X],W;\K[X]\otimes  g)\otimes_{\K\bb} \cKo^*(\K[X],W;\K[X]\otimes  h)\stackrel{\cup}\to \cKo^*(\K[X],W;\K[X]\otimes  gh),\\
\cKo_*(\K[X],W;\K[X]\otimes  g)\otimes_{\K\bb} \cKo^*(\K[X],W;\K[X]\otimes  h)\stackrel{\cap}\to \cKo_*(\K[X],W;\K[X]\otimes  gh).
\end{eqnarray*}
The following proposition is an immediate consequence of (\ref{formulacup}) and (\ref{formulacap}). 
\begin{proposition}\label{prodkos}
Under the isomorphism (\ref{tmg}), (\ref{omg}) the above products transform into 
\begin{equation*}
f_1q_1(\partial_{\theta})\cup f_2q_2(\partial_{\theta})
=f_1\cdot g(f_2)\cdot \LL(e^{t{\rm H}_W(x,g(x),gh(x))}\otimes q_1(\partial_{\theta})\otimes q_2(\partial_{\theta})),
\end{equation*}
\begin{equation*}
f_1p({\theta})\cap f_2q(\partial_{\theta})
=f_1\cdot g(f_2)\cdot \LLL(e^{t{\rm H}_W(g(x),gh(x),x)}\otimes p({\theta}) \otimes q(\partial_{\theta})),
\end{equation*}
respectively, where $f_1,f_2\in\K[X]$
and  $\LL$, $\LLL$ are the $\K[X]$-linear extensions of  (\ref{mu}) and (\ref{nu}).
\end{proposition}

\begin{remark}\label{nonab2}{\rm Note that these formulas are  valid for non-abelian groups as well.

}
\end{remark}

\subsubsection{$G$-actions}\label{sect5.4} Our final goal in this section is to transfer the $G$-actions on the Hochschild (co)homology $\cHH^*(\K[X],W;\K[X]\rtimes  G)$, $\cHH_*(\K[X],W;\K[X]\rtimes  G)$ (see  (\ref{Gactcoh})) to the corresponding Koszul (co)homology groups.  According to Proposition \ref{Ginv}, the (co)invariants of the resulting $G$-actions are isomorphic to  $\cHH^*(\K[X]\rtimes  G,W)$ and $\cHH_*(\K[X]\rtimes  G,W)$.

Since we are dealing with an abelian group, the $G$-actions on $\rHH^*(\K[X],W;\K[X]\rtimes G)$ and  $\rHH_*(\K[X],W;\K[X]\rtimes G)$ preserve the mixed subcomplexes 
$\rHH^*(\K[X],W;\K[X]\otimes g)$ and  $\rHH_*(\K[X],W;\K[X]\otimes g)$, for all $g$. Let us introduce the following $G$-action on the algebra ${\rm Cl}_N$:
\begin{equation}\label{Gtd} 
G\ni g=(g_1, \ldots, g_{N}): \quad \theta_i\mapsto g_i\theta_i,\quad \partial_{\theta_i}\mapsto g_i^{-1}\partial_{\theta_i}.
\end{equation}
Note it is compatible with the ${\rm Cl}_N$-module structures on $\K[\partial_{\theta}]$ and $\K[{\theta}]$ from Section \ref{CA}. The combination of the $G$-actions on $\K[X]$ and  ${\rm Cl}_N$ yields $G$-actions on various spaces we have been studying, e.g.  on $\K[X][\partial_{\theta}]$ and $\K[X][{\theta}]$, and also on the underlying space $\rK_*(\K[{X}])$ of the Koszul resolution of $\K[{X}]$. The latter action commutes with the differentials $\ddko$ and $\ddcurve$. This follows from the observation that  $\ddko$ and $\ddcurve$ are $G$-invariant and from the above-mentioned compatibility of the ${\rm Cl}_N$- and $G$-actions on $\K[{\theta}]$. 
As a result, we obtain $G$-actions on the mixed complexes 
\begin{eqnarray}\label{kco}
&&\rKo^*(\K[{X}],w; \K[X]\otimes g)=\left(\rK^*(\K[{X}], \K[X]\otimes g), \,\,\ddk, \,\,\ddw\right),\nonumber\\
&&\rKo_*(\K[{X}],w; \K[X]\otimes g)=\left(\rK_*(\K[{X}], \K[X]\otimes g), \,\,\bbk,\,\,\bbw\right).
\end{eqnarray}
defined on the underlying spaces of the complexes by analogy with (\ref{Gactcoh}). 

\begin{proposition} (a) The quasi-isomorphisms (see Corollary \ref{mainl})
\begin{eqnarray*}
&&\Psi^*: \rKo^*(\K[{X}],w; \K[X]\otimes g)\to \rHH^*(\K[{X}],w; \K[X]\otimes g),\\
&&\Psi_*: \rHH_*(\K[{X}],w; \K[X]\otimes g)\to \rKo_*(\K[{X}],w; \K[X]\otimes g)
\end{eqnarray*}  are $G$-equivariant. 

\noindent (b) The isomorphisms (\ref{tm}) and (\ref{om}) between the mixed complexes on the right-hand sides of (\ref{kco}) and the mixed complexes
\begin{eqnarray}\label{koskom1}
&&((\K[X]\otimes g)[\partial_{\theta}],\ddk,\ddw)\ (\simeq(\K[X][\partial_{\theta}],\ddkg,\ddwg)),\nonumber\\
&&((\K[X]\otimes g)[{\theta}],\bbk,\bbw)\ (\simeq(\K[X][{\theta}],\bbkg,\bbwg))
\end{eqnarray}
are also $G$-equivariant. 
\end{proposition}
\noindent  Part (a) follows from the observation  that the morphism (\ref{psi}) intertwines the $G$-actions on the bar and  Koszul resolutions. Part (b) is obvious for the second complex; for the first one the claim follows from the fact that the pairing (\ref{pair1}) is $G$-invariant.

Let us also point out that the $G$-actions on (\ref{koskom1}) are compatible with the explicit description of the cohomology classes  that we worked out in Remarks \ref{explcoh}, \ref{explhom}. More precisely, the splitting of the differentials $\ddwg=\ddwg'+\ddwg''$ and $\bbwg=\bbwg'+\bbwg''$ is $G$-invariant and the isomorphisms (\ref{iso}) and (\ref{expiso}) intertwine the $G$-actions. 

\begin{remark}\label{nonab3}{\rm Unlike the results of the previous two sections (cf. Remarks \ref{nonab1}, \ref{nonab2}),  the above picture relies very heavily on the fact that $G$ acts by rescaling the variables. The major problem with other groups is that the morphism (\ref{psi}) will not be equivariant anymore, even if $G\subset GL_N(\K)$. 
}
\end{remark}

%%%%%%%%%%%%%%%%%%%%%%%%%%%%%%%%%

\subsection{Proof of Theorem \ref{mainthm}}\label{sect6.1} In this section, we assume that $W$ has only isolated critical points. Then  each  $W^g\in\K[X^g]$ also has isolated critical points. Indeed (cf.  \cite[Lem.2.5.3]{PV}),  the equality $W(g(x))=W(x)$ implies (by differentiating both sides) 
\[
{\sf res}^g(\partial_{x_i}W)=0\quad \forall\, i\in I_g,\qquad {\sf res}^g(\partial_{x_i}W)=\partial_{x_i}W^g \quad \forall\, i\in I^g
\] 
which in turn implies that ${\sf res}^g$ induces a surjective map ${\rm H}^*(\Omega^*_X,  dW\wedge\cdot)\to {\rm H}^*(\Omega^*_{X^g}, dW^g\wedge\cdot)$.

Let us start by writing down explicit isomorphisms 
\begin{eqnarray}\label{spiso1}
&&\Jaco^*(X,W,G)\simeq \bigoplus_{g\in G}{\rm H}^*(\K[X][\partial_{\theta}]\bb,\ddkg+t\ddwg),\nonumber\\  && \Omeg^*(X,W,G)\simeq \bigoplus_{g\in G} {\rm H}^*(\K[X][{\theta}]\bb,\bbkg+t\bbwg)
\end{eqnarray}
as $(\Z\times G)$-graded $G$-equivariant $\K\bb$-modules. 
Since all $W^g$ have isolated critical points, Propositions \ref{invsec1}, \ref{invsec2} and Remarks \ref{explcoh}, \ref{explhom} imply that
 $e^{t  {\rm H}_{{W,g}}}:\K[X][\partial_{\theta}]\bb\to \K[X][\partial_{\theta}]\bb$ induces an isomorphism
\begin{equation}\label{jako1}
e^{t  {\rm H}_{{W,g}}}: M(W^g)[{\uu}^{\pm1}]\,\partial_{\theta_{I_g}}\simeq {\rm H}^*(\K[X][\partial_{\theta}]\bb,\ddkg+t\ddwg)
\end{equation}
and  $e^{t  {\rm H}_{{W,g}}^\dagger}:\K[X][{\theta}]\bb\to \K[X][{\theta}]\bb$ induces an isomorphism
\begin{equation}\label{jako2}
e^{t  {\rm H}_{{W,g}}^\dagger}:M(W^g)[{\uu}^{\pm1}]\,{\theta_{I^g}}\simeq {\rm H}^*(\K[X][{\theta}]\bb,\bbkg+t\bbwg).
\end{equation}
 The specific isomorphisms (\ref{spiso1}) we are interested in are the $\K\bb$-linear extensions of the maps
\begin{equation}\label{spiso}
 \lfloor f\rfloor_g\cdot \xi_g\mapsto  \lfloor f\rfloor_g\cdot [e^{t  {\rm H}_{{W,g}}}(\partial_{\theta_{I_g}})],\qquad 
  \lfloor f\rfloor_g\cdot \omega_g\mapsto  (-1)^{Nd_g+||I_g||}\, \lfloor f\rfloor_g\cdot [e^{t  {\rm H}_{{W,g}}^\dagger}({\theta_{I^{g}}})]
\end{equation}
($f\in\K[X]$) where $[\,\,]$ denotes the cohomology class  and $||(i_1,i_2,\ldots,i_k)||:=i_1+i_2+\ldots+i_k$. Obviously, these isomorphisms preserve the gradings and respect the $G$-actions (cf. (\ref{Gvw}), (\ref{Gtd})). 

Let us now calculate the cup products between elements of the form $\lfloor f\rfloor_g\cdot [e^{t  {\rm H}_{{W,g}}}(\partial_{\theta_{I_g}})]$.
To begin with, the existence of $\sigma_{g,h}\in M(W^{gh})$ such that
\begin{equation*}
[e^{t  {\rm H}_{{W,g}}}(\partial_{\theta_{I_g}})]\cup [e^{t  {\rm H}_{{W,h}}}(\partial_{\theta_{I_h}})]=t^{d_{g,h}} \cdot \sigma_{g,h} \cdot [e^{t  {\rm H}_{{W,gh}}}(\partial_{\theta_{I_{gh}}})]
\end{equation*} 
is an immediate consequence of (\ref{jako1}). 
That $\sigma_{g,h}$ is given by the formula (\ref{cgh}) is a consequence of the formula for $\cup$  from Proposition \ref{prodkos} since (\ref{cgh}) is nothing but the class in $M(W^{gh})$ of the coefficient at $t^{d_{g,h}} \partial_{\theta_{I_{gh}}}$ in 
$
e^{-t  {\rm H}_{{W,gh}}} \left(e^{t  {\rm H}_{{W,g}}}(\partial_{\theta_{I_g}})\cup e^{t  {\rm H}_{{W,h}}}(\partial_{\theta_{I_h}})\right).
$
(Note that according to the formulas in Proposition \ref{prodkos} we should be using ${\rm H}_W(x,g(x),gh(x))$ instead of ${\rm H}_W(x,g(x),x)$ in (\ref{cgh}). But, obviously, after applying $\lfloor\,\rfloor_{gh}$ the result is the same. Note also that ${\rm H}_{{W,gh}}$ does not appear in (\ref{cgh}) because for any $k\in G$ the element $\partial_{\theta_{I_{k}}}$ is not contained in the image of ${\rm H}_{{W,k}}$; see (\ref{hg}).)

Furthermore, by Proposition \ref{prodkos}
the elements  $\sigma_{g,h}$  determine the products on the entire cohomology groups, namely
\[
\lfloor f_1\rfloor_g\cdot [e^{t  {\rm H}_{{W,g}}}(\partial_{\theta_{I_g}})]\cup \lfloor f_2\rfloor_h\cdot [e^{t  {\rm H}_{{W,h}}}(\partial_{\theta_{I_h}})]=t^{d_{g,h}} \cdot \lfloor f_1g(f_2)\rfloor_{gh} \cdot \sigma_{g,h} \cdot [e^{t  {\rm H}_{{W,gh}}}(\partial_{\theta_{I_{gh}}})].
\]
for all $f_1,f_2\in\K[X]$. (We implicitly use the obvious fact that the operator $\lfloor\,\rfloor_g$ commutes with the $G$-actions on $\K[X]$ and $M(W^g)$.) Note the difference between the right-hand side of the latter formula and that of (\ref{cupjac}): In  (\ref{cupjac}) $f_2$ is not twisted by $g$. We claim that the twist by $g$ is not needed, i.e. 
$\lfloor f_1g(f_2)\rfloor_{gh} \cdot \sigma_{g,h} =\lfloor f_1f_2\rfloor_{gh} \cdot \sigma_{g,h}$ as elements of $M(W^{gh})$.
 This is a consequence of the braided super-commutativity (\ref{bracom}) of the cup product. Indeed, by the braided super-commutativity
\begin{multline}\label{ooh}
t^{d_{g,h}} \cdot \sigma_{g,h} \cdot [e^{t  {\rm H}_{{W,gh}}}(\partial_{\theta_{I_{gh}}})]=[e^{t  {\rm H}_{{W,g}}}(\partial_{\theta_{I_g}})]\cup [e^{t  {\rm H}_{{W,h}}}(\partial_{\theta_{I_h}})]=\\
=(-1)^{d_gd_h}\cdot [e^{t  {\rm H}_{{W,h}}}(\partial_{\theta_{I_h}})]\cup h^{-1}\left([e^{t  {\rm H}_{{W,g}}}(\partial_{\theta_{I_g}})]\right)=\alpha\cdot (-1)^{d_gd_h}\cdot t^{d_{g,h}} \cdot \sigma_{h,g} \cdot [e^{t  {\rm H}_{{W,gh}}}(\partial_{\theta_{I_{gh}}})]
\end{multline}
where the constant $\alpha$ is defined by $h^{-1}(\partial_{\theta_{I_g}})= \alpha\cdot \partial_{\theta_{I_g}}$.
Hence  
\begin{multline*}
t^{d_{g,h}} \cdot \lfloor f_1g(f_2)\rfloor_{gh} \cdot \sigma_{g,h} \cdot [e^{t  {\rm H}_{{W,gh}}}(\partial_{\theta_{I_{gh}}})]
=\lfloor f_1\rfloor_g \cdot[e^{t  {\rm H}_{{W,g}}}(\partial_{\theta_{I_g}})]\cup \lfloor f_2\rfloor_h \cdot[e^{t  {\rm H}_{{W,h}}}(\partial_{\theta_{I_h}})]=\\
=(-1)^{d_gd_h}\cdot \lfloor f_2\rfloor_h \cdot[e^{t  {\rm H}_{{W,h}}}(\partial_{\theta_{I_h}})]\cup h^{-1}\left(\lfloor f_1\rfloor_g\right)\cdot h^{-1}\left([e^{t  {\rm H}_{{W,g}}}(\partial_{\theta_{I_g}})]\right)=\\
=(-1)^{d_gd_h} \cdot\lfloor f_2\rfloor_h\cdot [e^{t  {\rm H}_{{W,h}}}(\partial_{\theta_{I_h}})]\cup \lfloor h^{-1}\left(f_1\right)\rfloor_g \cdot h^{-1}\left([e^{t  {\rm H}_{{W,g}}}(\partial_{\theta_{I_g}})]\right)=\\
=\alpha\cdot (-1)^{d_gd_h} \cdot\lfloor f_2\rfloor_h\cdot [e^{t  {\rm H}_{{W,h}}}(\partial_{\theta_{I_h}})]\cup \lfloor h^{-1}\left(f_1\right)\rfloor_g \cdot [e^{t  {\rm H}_{{W,g}}}(\partial_{\theta_{I_g}})]=\\
=\alpha\cdot (-1)^{d_gd_h} \cdot t^{d_{g,h}} \cdot \lfloor f_2f_1\rfloor_{gh} \cdot \sigma_{h,g} \cdot [e^{t  {\rm H}_{{W,gh}}}(\partial_{\theta_{I_{gh}}})]
\stackrel{(\ref{ooh})}= t^{d_{g,h}} \cdot \lfloor f_1f_2\rfloor_{gh} \cdot \sigma_{g,h} \cdot [e^{t  {\rm H}_{{W,gh}}}(\partial_{\theta_{I_{gh}}})].
\end{multline*}

This  completes the proof of the first half of Theorem \ref{mainthm}. The second half is proved similarly, so we only sketch the argument. 

It follows from (\ref{jako2})  that there exist elements $\widetilde{\sigma}_{g,h}\in M(W^{gh})$ such that
\[
[e^{t  {\rm H}_{{W,g}}^\dagger}({\theta_{I^g}})]\cap [e^{t  {\rm H}_{{W,h}}}(\partial_{\theta_{I_h}})]=t^{d_{g,h}} \cdot \widetilde{\sigma}_{g,h}\cdot [e^{t  {\rm H}^\dagger_{{W,gh}}}({\theta_{I^{gh}}})]
\]
We claim that $\widetilde{\sigma}_{g,h}=(-1)^{N(d_g+d_{gh})+||I_g||+||I_{gh}||}\, \sigma_{g,h}$  (cf. (\ref{spiso})). Indeed, one can show, using  (\ref{formulacap}), that $\widetilde{\sigma}_{e,g}=(-1)^{Nd_g+||I_g||}$. Hence
\begin{multline*}
[e^{t  {\rm H}_{{W,g}}^\dagger}({\theta_{I^g}})]\cap [e^{t  {\rm H}_{{W,h}}}(\partial_{\theta_{I_h}})]=
(-1)^{Nd_g+||I_g||}\left([{\theta_{I^e}}]\cap [e^{t  {\rm H}_{{W,g}}}(\partial_{\theta_{I_g}})]\right)\cap [e^{t  {\rm H}_{{W,h}}}(\partial_{\theta_{I_h}})]=\\
=(-1)^{Nd_g+||I_g||}\cdot[{\theta_{I^e}}]\cap \left([e^{t  {\rm H}_{{W,g}}}(\partial_{\theta_{I_g}})]\cup [e^{t  {\rm H}_{{W,h}}}(\partial_{\theta_{I_h}})]\right)=\\
=(-1)^{Nd_g+||I_g||}\cdot t^{d_{g,h}} \cdot \sigma_{g,h} \cdot [{\theta_{I^e}}]\cap  [e^{t  {\rm H}_{{W,gh}}}(\partial_{\theta_{I_{gh}}})]=\\
=(-1)^{Nd_g+||I_g||}\cdot (-1)^{Nd_{gh}+||I_{gh}||}\cdot t^{d_{g,h}} \cdot \sigma_{g,h} \cdot [e^{t  {\rm H}^\dagger_{{W,gh}}}({\theta_{I^{gh}}})].
\end{multline*}
Thus,
\[
(-1)^{Nd_g+||I_g||}\cdot[e^{t  {\rm H}_{{W,g}}^\dagger}({\theta_{I^g}})]\cap [e^{t  {\rm H}_{{W,h}}}(\partial_{\theta_{I_h}})]=
 t^{d_{g,h}} \cdot \sigma_{g,h} \cdot (-1)^{Nd_{gh}+||I_{gh}||}\cdot[e^{t  {\rm H}^\dagger_{{W,gh}}}({\theta_{I^{gh}}})]
\]
which proves that the second map in (\ref{spiso}) respects the cap products on the generators. The extension to the entire (co)homology groups is completely parallel to the case of the cup product and is left to the reader. 
 
Finally, the claims that $\xi_e$ is the unit of $\Jaco^*(X,W,G)$ and $\omega_e$ is a free generator of $\Omeg^*(X,W,G)$ both follow from the obvious fact that $\sigma_{e,g}=1$ for any $g\in G$.
%%%%%%%%%%%%%%%%%%%%%%%%%%%%%%%%%
\subsection{Proof of Theorem \ref{mf}}\label{MFC} The proof will be outlined very schematically since it involves  notions and results of the theory of curved dg categories and their Hochschild invariants \cite{PP} which are far beyond the scope of the present work. (Unfortunately, we have to assume familiarity with the subject.)  
Also, we only sketch the proof of the isomorphism
\begin{equation}\label{ht}
\left(\ccHH^{*}({\rm MF}_G(X,W)), \cup\,\right)\simeq \left(\ccHH^{\rm II,*}(\K[X]\rtimes  G,W), \cup\,\right).
\end{equation}
The proof of the other half is just a straightforward extension of the argument given below  (for the homology all the arrows below should be reversed). 

To start with, we want to replace the ordinary Hochschild cohomology of the matrix factorizations with their Hochschild cohomology of the second kind. More precisely, 
there is a natural algebra homomorphism \cite[Sect.2.4, (24)]{PP}
\begin{equation}\label{IIto}
\ccHH^{\rm II,*}({\rm MF}_G(X,W))\to\ccHH^*({\rm MF}_G(X,W))
\end{equation}
which is an isomorphism in our case. It is this point where the absence of critical points outside of  $W^{-1}(0)$ is needed; see \cite[Sect.4.8,4.10]{PP}. (Formally speaking, this is spelled out in {\it loc.cit.} only in the non-equivariant setting but the same argument, when combined with the discussion in \cite[Sect.2.5]{PV}, applies in the equivariant case.) 

Furthermore, there is a diagram of morphisms 
\begin{equation}\label{diagr}
\ccHH^{\rm II,*}({\rm MF}_G(X,W))\stackrel{I_1}{\longleftarrow} \ccHH^{\rm II,*}(\mathfrak{A})
\stackrel{I_2}{\longrightarrow}\ccHH^{\rm II,*}(\K[X]\rtimes G,W)
\end{equation}
where 

\noindent (1) $\mathfrak{A}$ is the curved dg category whose objects are the $\ZZ$-graded finitely generated projective $(\K[X]\rtimes G)$-modules $E=E^{\rm even}\oplus E^{\rm odd}$, the curvature of every object $E$ is $W\cdot{\rm id}_E$, and  $\Hom_{\mathfrak{A}}(E_1,E_2)$ is the $\ZZ$-graded space of all $(\K[X]\rtimes G)$-linear maps from $E_1$ to $E_2$; we endow this space with the trivial differential.

\vspace{0.05in}

\noindent (2) $I_2$  is induced by the embedding of curved dg categories $(\K[X]\rtimes G,W)\to \mathfrak{A}$ (we view the former as a curved dg category with a single object) sending  the unique object of $(\K[X]\rtimes G,W)$ to the object with $E^{\rm even}=A$ and $E^{\rm odd}=0$. 

\vspace{0.05in}

\noindent (3) $I_1$ is the dual to the map $I_1^\vee$ from the bar resolution of the second kind of ${\rm MF}_G(X,W)$ to the bar resolution of the second kind of $\mathfrak{A}$ given by the following explicit formula (see \cite[Sect.2.4, (18)]{PP} and \cite[Sect.2.3]{Se}): 
\[
I_1^\vee={\rm exp}\ (\text{``insert $\delta$''})
\] where
\begin{equation*}
\text{``insert $\delta$''}: \phi_{E_{0}E_{1}}\ [\ \phi_{E_{1}E_{2}}\ |\  \ldots|\ \phi_{E_{n-1}E_n}\ ]
\ \phi_{E_{n}E_{n+1}}
\mapsto
 \sum_{i=1}^n\phi_{E_{0}E_{1}}\ [\ \ldots|\ \phi_{E_{i-1}E_i}|\ \delta_{E_i}\ |\ \phi_{E_{i}E_{i+1}}\ | \ldots\ ] \ \phi_{E_{n}E_{n+1}}.
\end{equation*}
On the left-hand side, $(E_i,\delta_{E_i})$ are matrix factorizations and  $\phi_{E_{i}E_{i+1}}$ are morphisms in the category
${\rm MF}_G(X,W)$; on the right-hand side, $E_i$, $\phi_{E_{i}E_{i+1}}$, and $\delta_{E_i}$ are viewed as objects/morphisms  in $\mathfrak{A}$. 

\vspace{0.05in}

\noindent The proof of (\ref{ht}) is completed by noticing that, firstly, both $I_1$ and $I_2$ are isomorphisms (this is obvious for $I_1$ and follows from \cite[Sect.2.6,(45)]{PP} for  $I_2$) and, secondly, both $I_1$ and $I_2$ are morphisms of algebras (this is obvious for $I_2$; for $I_1$ this follows the fact that ``insert $\delta$'' is a coderivation with respect to the coproduct $\Delta_{\rm bar}$ (\ref{coalg0})).

\newpage

\appendix

\section{ Examples and applications ({\it by A. Basalaev and D. Shklyarov})}\label{A}

In Appendix A,  the ground field is $\C$. 

\subsection{Example: Invertible polynomials}\label{appA.1}
The aim of this section is to discuss in more detail the (isomorphism class of the) algebra $\Jacoii(X,W,G)$ from Corollary \ref{maincorii} in the case when $X=\C^N$ and $W\in\C[X]$ is an {\it  invertible} polynomial \cite{Kre,KS}.

Instead of giving a formal definition of invertible polynomials, we recall the only fact that matters to us, namely, the classification of such polynomials obtained in \cite{KS} which says that up to a natural equivalence, any invertible polynomial is the Thom-Sebastiani sum of polynomials of the following  {\it atomic types}: 
\begin{description}
 \item[\it Ferma type] $x_1^{a_1}$,
 \item[\it Chain type] $x_1^{a_1}x_2 + x_2^{a_2}x_3 + \dots + x_{N-1}^{a_{N-1}}x_N + x_N^{a_N}$,
 \item[\it Loop type] $x_1^{a_1}x_2 + x_2^{a_2}x_3 + \dots + x_{N-1}^{a_{N-1}}x_N + x_N^{a_N}x_1$
\end{description}
where $a_k \in {\mathbb{N}}_{\ge 2}$ and $N \ge 2$. (The Thom-Sebastiani sum of two functions $(X_1,W_1)$ and $(X_2,W_2)$ is the function 
$ (X_1\times X_2, W_1\boxplus W_2)$ with $(W_1\boxplus W_2)(x_1,x_2):=W_1(x_1)+W_2(x_2)$.)

Note that every invertible polynomial has an isolated critical point at the origin and, due to its quasi-homogeneity, no other critical points. Thus, the main results of the present work are applicable and yield a description of the Hochschild cohomology, i.e. of the ``closed string algebra'' of any invertible LG model. 
But, as was already mentioned in the Introduction, there already exists \cite{BTW1} a good candidate for this role which, like  $\Jacoii(\C^N,W,G)^G$,  is constructed as the subalgebra of $G$-invariants of a  $(\ZZ\times G)$-graded, $G$-equivariant, and braided super-commutative algebra. This latter algebra is called the $G$-twisted Jacobian algebra of $W$ and denoted by $\Jac'(W,G)$.  Let us recall its description.

It suffices to describe $\Jac'(W,G)$ for each atomic polynomial since the algebra has the following ``K\"unneth property'':  If $W=W_1\boxplus\ldots\boxplus W_l$, where each $W_i$ is a polynomial of one of the three atomic types, and $G=G_1\times\ldots\times G_l$, where $G_i$ is an abelian group of symmetries of $W_i$, then
\begin{equation}\label{decom}
\Jacob(W,G)=\Jacob(W_1,G_1)\otimes\ldots \otimes \Jacob(W_l,G_l).
\end{equation}

Let  $W$ be an atomic polynomial and $G$ be its abelian symmetry group. The algebra $\Jac'(W,G)$ looks exactly the same as $\Jacoii(\C^N,W,G)$ -- that is, it is isomorphic to $\Jacoii(\C^N,W,G)$ as a $(\ZZ\times G)$-graded $G$-equivariant $M(W)$-module -- but it has different ``structure constants'' $\sigma'_{g,h}$ in the products  (\ref{cupjac}), namely 
\begin{eqnarray*}
\sigma'_{e,g}=\sigma'_{g,e}=1,\qquad 
\sigma'_{g,g^{-1}}=e^{- \pi \sqrt{-1} \,\age(g)} \lfloor\det(\partial_{x_i}\partial_{x_j} W)_{i,j\in I_g}\rfloor_e,\qquad \forall g\in G
\end{eqnarray*}
and $\sigma'_{g,h}=0$ otherwise. In the above formula $\lfloor\,\rfloor_e$, we recall, 
 denotes the class of an element in $M(W=W^e)$ and $\age(g) := \sum_{i=1}^N q_i$ with $q_i$ being the  rational numbers satisfying $ 0 \le q_i < 1$ and $g= {\rm diag}(e^{2 \pi \sqrt{-1}q_i})$.   (The actual structure constants of $\Jac'(W,G)$, as defined in \cite{BTW1}, differ from $\sigma'_{g,h}$ by complex factors which we ignore here since they do not affect the isomorphism class.) As explained in \cite[Sect.4]{BTW1}, the map $ \xi_g\otimes \xi_h\mapsto \sigma'_{g,h} \xi_{gh}$ does extend by $M(W)$-linearity to a well-defined associative braided super-commutative product on the whole of $\Jac'(W,G)$.

We propose the following conjecture:

\vspace{0.1in}

\noindent {\bf Conjecture.}
{\it For any invertible $W$ and any abelian symmetry group $G$  there is an isomorphism of $(\ZZ\times G)$-graded $G$-equivariant algebras 
$\Jacoii(\C^N,W,G) \simeq \Jacob(W,G)$.}

\vspace{0.1in}

\begin{remark}{ \rm \noindent (1) Note that Proposition \ref{Kunn} implies that $\Jacoii(\C^N,W,G)$ has the property  (\ref{decom}) as well. Thus, it would suffice to prove the conjecture for polynomials of the three atomic types. The Ferma case is an easy exercise,  the conjecture is interesting only for the other two types. 

\noindent (2) There is also an analog of $\Omegii(\C^N,W,G)$ in \cite{BTW1} which the authors denote by $\Omega'_{W,G}$. Its structure is completely analogous to that of  $\Omegii(\C^N,W,G)$:  it is a free rank 1 $\Jac'(W,G)$-module spanned by a generator of degree $-N$. It is not included in the conjecture for the simple reason that it is defined in {\it loc. cit.} as a {\it left} module.
 }
\end{remark}

In the remainder of this section we check the conjecture for a generic chain type polynomial in two variables and its maximal (hence any) symmetry group:
\[
W=x_1^{a_1}x_2 +  x_2^{a_2}\quad (a_2\geq3),\qquad G:=\{(\zeta_1,\zeta_2)\in(\C^{*})^2\,|\,\zeta_1^{a_1}\zeta_2=1,\,\zeta_2^{a_2}=1\}
\] 

Our strategy is as follows: we are going to show that

\begin{description}
\item[A] $\sigma_{g,h}=0$ provided neither of $g,h,gh$ is the unit of $G$.
\item[B] $\sigma_{g,g^{-1}}=\alpha_g\cdot \sigma'_{g,g^{-1}}$ for some $\alpha_g\in\C^*$.
\end{description}
(It is worthwhile noting that $\sigma'_{g,g^{-1}}\in M(W)$ is always non-zero; see  \cite[Prop.30]{BTW1}.) Let us explain why checking these conditions suffices to verify the claim of the conjecture. 
It follows from the braided super-commutativity of both $\Jacoii(\C^2,W,G)$ and $\Jacob(W,G)$ that 
\[
\sigma_{g^{-1},g}=\det(g)\sigma_{g,g^{-1}}, \quad \sigma'_{g^{-1},g}=\det(g)\sigma'_{g,g^{-1}}.
\]
Together with {\bf B} this implies that $\alpha_g=\alpha_{g^{-1}}$ which, in turn, implies
that the assignment 
\[
\xi^{\Jacoii}_g\mapsto \sqrt{\alpha}_g\cdot \xi^{\Jac'}_g
\]
-- for any choice of the square root of the function $G\to\C^*, g\mapsto \alpha_g$ satisfying $\sqrt{\alpha}_g=\sqrt{\alpha}_{g^{-1}}$ -- extends by $M(W)$-linearity to an algebra isomorphism $\Jacoii(\C^2,W,G)\to \Jacob(W,G)$. 

Let us calculate all the ingredients for the formula (\ref{cgh}) in our case. One has 
\begin{multline*}
    {\rm H}_W(x,y,z) = \left( \frac{x_1^{a_1}-y_1^{a_1}}{x_1-y_1} - \frac{x_1^{a_1}-z_1^{a_1}}{x_1-z_1} \right)\frac{x_2}{y_1-z_1} \,\theta_1\otimes\theta_1+
\\
+\left( \frac{x_2^{a_2}-y_2^{a_2}}{x_2-y_2} - \frac{x_2^{a_2}-z_2^{a_2}}{x_2-z_2} \right)\frac{1}{y_2-z_2} \,\theta_2\otimes\theta_2
+\frac{y_1^{a_1}-z_1^{a_1}}{y_1-z_1} \,\theta_{2}\otimes\theta_{1}
\end{multline*}
which implies that for $g=(\zeta_1,\zeta_2)\neq e$ ($\Leftrightarrow\,\, \zeta_1\neq1$)
\begin{equation}\label{Hwxgxx1}
    {\rm H}_W(x,g(x),x) = \left( \frac{a_1}{1-\zeta_1} - \frac{1-\zeta_1^{a_1}}{(1-\zeta_1)^2} \right) x_1^{a_1-2}x_2\,\theta_1\otimes\theta_1
+\frac{a_2}{1-\zeta_2} x_2^{a_2-2} \,\theta_2\otimes\theta_2
+\frac{1-\zeta_1^{a_1}}{1-\zeta_1} x_1^{a_1-1}\,\theta_{2}\otimes\theta_{1},
\end{equation}
if $\zeta_2\neq1$ ($\Leftrightarrow\,\, d_g=2$), and
\begin{equation}\label{Hwxgxx2}
    {\rm H}_W(x,g(x),x) = \frac{a_1}{1-\zeta_1} x_1^{a_1-2}x_2\,\theta_1\otimes\theta_1
+\frac{a_2(a_2-1)}{2} x_2^{a_2-2}\,\theta_2\otimes\theta_2
+\frac{1-\zeta_1^{a_1}}{1-\zeta_1} x_1^{a_1-1}\,\theta_{2}\otimes\theta_{1},
\end{equation}
if $\zeta_2=1$ ($\Leftrightarrow\,\, d_g=1$). 
Also,
\begin{equation}\label{Hgx}
{\rm H}_{{W,g}}(x) = \frac{\zeta_1^{a_1}}{1 - \zeta_1} x_1^{a_1-1}\, \theta_1\theta_{2} \quad\text{if}\,\, \zeta_2\neq1\qquad \text{and}\qquad {\rm H}_{{W,g}}(x)=0\quad\text{if}\,\,   \zeta_2=1.
\end{equation}

\vspace{0.1in}

We are in a position now to check the above conditions {\bf A} and {\bf B}.

\vspace{0.1in}

\noindent {\bf A:}  Let us fix $g,h\in G$ such that neither of $g,h,gh$ is the unit. We have to consider the following possibilities:
\begin{enumerate}
\item $d_g=d_h=d_{gh}=1$;
\item $d_g=1$, $d_h=d_{gh}=2$ or $d_h=1$, $d_g=d_{gh}=2$;
\item $d_g=d_h=d_{gh}=2$.
\end{enumerate}
However, because of (\ref{non-int}), only the last possibility is non-trivial, so we assume $d_g=d_h=d_{gh}=2$ (consequently, $d_{g,h}=1$). 

By Theorem \ref{mainthm}   $\sigma_{g,h}$ is the coefficient at $\partial_{\theta_1}\partial_{\theta_2}$ in the expression
\[
\LL\left(\lfloor{\rm H}_W(x,g(x),x)\rfloor_{gh}+\lfloor{\rm H}_{{W,g}}(x)\rfloor_{gh}\otimes1+ 1\otimes\lfloor {\rm H}_{{W,h}}(g(x))\rfloor_{gh}\right)\otimes\partial_{\theta_1}\partial_{\theta_2}\otimes \partial_{\theta_1}\partial_{\theta_2}).
\] 
Since $d_{gh}=2$, the operator $\lfloor\,\rfloor_{gh}$ sets both $x_1$ and $x_2$ equal to 0. Thus, by (\ref{Hwxgxx1}) and (\ref{Hgx})  the above expression is 0.

\vspace{0.1in}

\noindent {\bf B:} We fix now $g=(\zeta_1,\zeta_2)\neq e$ and consider two cases:

\vspace{0.1in}

\noindent \underline{$\zeta_2\neq1$ ($\Leftrightarrow\,\, d_g=2$):} In this case $\sigma_{g,g^{-1}}$ is the constant coefficient  in the expression
\[
\frac12\LL\left(\left(\lfloor{\rm H}_W(x,g(x),x)\rfloor_{e}+\lfloor{\rm H}_{{W,g}}(x)\rfloor_{e}\otimes 1+1\otimes \lfloor {\rm H}_{{W,g^{-1}}}(g(x))\rfloor_{e}\right)^2\otimes \partial_{\theta_1}\partial_{\theta_2} \otimes
\partial_{\theta_1}\partial_{\theta_2}\right)
\] 
for ${\rm H}_W(x,g(x),x)$ given by (\ref{Hwxgxx1}). 
Denoting 
\begin{eqnarray}\label{ABC}
{\rm H}_W(x,g(x),x) = A_1\,\theta_1\otimes\theta_1
+A_2 \,\theta_2\otimes\theta_2
+A_3\,\theta_{2}\otimes\theta_{1},\\ {\rm H}_{{W,g}}(x) = B\, \theta_1\theta_{2}, \quad {\rm H}_{W,g^{-1}}(g(x)) = C\, \theta_1\theta_{2},
\end{eqnarray}
one easily checks that $\sigma_{g,g^{-1}}=\lfloor B\cdot C- A_1\cdot A_2\rfloor_e$, that is
\[
\sigma_{g,g^{-1}}=\left\lfloor -\frac{\zeta_1^{a_1}}{(1 - \zeta_1)^2} x_1^{2a_1-2}-\left( \frac{a_1}{1-\zeta_1} - \frac{1-\zeta_1^{a_1}}{(1-\zeta_1)^2} \right)\frac{a_2}{1-\zeta_2}  x_1^{a_1-2}x_2^{a_2-1}\right\rfloor_e
\] 
or, taking account of the facts that $\zeta_1^{a_1}=\zeta_2^{-1}$ and $\lfloor x_1^{2a_1-2}\rfloor_e=\lfloor -a_2x_1^{a_1-2}x_2^{a_2-1}\rfloor_e$,
\[
\sigma_{g,g^{-1}}=-\frac{a_1a_2}{(1 - \zeta_1)(1 - \zeta_2)}  \left\lfloor x_1^{a_1-2}x_2^{a_2-1}\right\rfloor_e.
\] 
The reader is invited to check that the class of the Hessian of $W$ in $M(W)$ is proportional to the above element which completes the verification of {\bf B} in the present case. 

\vspace{0.1in}

\noindent \underline{$\zeta_2=1$ ($\Leftrightarrow\,\, d_g=1$):}  This time $\sigma_{g,g^{-1}}$ is the constant coefficient  in the expression
\[
\LL\left(\lfloor{\rm H}_W(x,g(x),x)\rfloor_{e}\otimes \partial_{\theta_1}\otimes\partial_{\theta_1}\right)
\] 
for ${\rm H}_W(x,g(x),x)$ given by (\ref{Hwxgxx2}). Using the notation (\ref{ABC}), 
one sees that 
\[
\sigma_{g,g^{-1}}=\lfloor - A_1\rfloor_e =-\frac{a_1}{1-\zeta_1} \lfloor x_1^{a_1-2}x_2 \rfloor_e.
\] 
On the other hand, $\partial^2_{x_1} W=a_1(a_1-1)x_1^{a_1-2}x_2$, so  {\bf B} holds  in this case as well. 

\subsection{Application: Hochschild cohomology of Fukaya categories of surfaces}\label{A.2}
Let us fix an integer $\gen\geq2$ and  a symplectic compact connected oriented surface $S$ of genus $\gen$. Let $\mathscr{F}(S)$ denote the ($\ZZ$-graded $\C$-linear) Fukaya  $A_\infty$ category of $S$ as defined in \cite{Sei}. Our aim in this section is to combine Corollary \ref{mf} with the homological mirror symmetry theorem for surfaces established in  \cite{Sei} ($\gen=2$) and \cite{Ef} ($\gen\geq3$) in order to  prove the following claim: 

\begin{theorem}\label{A2} There is an isomorphism of $\ZZ$-graded algebras $\ccHH^*(\mathscr{F}(S))\simeq {\rm H}^*(S,\C)$.
\end{theorem}

\noindent We should emphasize that the theorem itself is not new; we only present a new proof. The claim can also  be deduced -- again in combination with the theorems of \cite{Ef,Sei} which imply that $\mathscr{F}(S)$ is homologically smooth -- from results of \cite{Gan,GPS1,GPS2} (see, in particular, Corollary 7 in \cite[Sect.1.2]{Gan}). In fact, the latter approach yields more than the mere existence of an isomorphism of algebras. It shows that a {\it specific} map ${\rm H}^*({S},\C)\to \ccHH^*(\mathscr{F}({S}))$, the so-called closed-open map, is an isomorphism. 

\vspace{0.1in}

Let us proceed to the proof. 
Of course, the first step is to apply the mirror symmetry theorem of \cite{Ef,Sei} and thereby convert  Theorem \ref{A2} into a claim about Hochschild invariants of specific LG models. Namely,  the category $\mathscr{F}({S})$ is shown in \cite{Ef,Sei} to be derived Morita equivalent to the dg category ${\rm MF}_G(\C^3,W)$ of equivariant matrix factorizations associated with the pair
\begin{equation}\label{LGg}
W:= x_1^{2\gen + 1} + x_2^{2\gen + 1} + x_3^{2\gen+1} - x_1x_2x_3, \qquad G:=\{(\zeta,\zeta,\zeta^{-2})\in(\C^{*})^3\,|\,\zeta^{2\gen+1}=1\}.
\end{equation}
It follows, according to \cite{Ke},  that $
\ccHH^*(\mathscr{F}({S}))\simeq\ccHH^*({\rm MF}_G(\C^3,W))
$
and so proving the theorem reduces to showing  the existence of an algebra isomorphism
$
\ccHH^*({\rm MF}_G(\C^3,W))\simeq {\rm H}^*({S},\C).
$

The problem now is that the above $W$ does not satisfy the condition of Corollary \ref{mf}. In addition to the origin $\orig\in\C^3$,  $W$ has $(2\gen+1)^2\ (2\gen-2)$ other isolated critical points where the critical values are different from 0. To apply Corollary \ref{mf}, we need to replace $\C^3$ by an open affine subset $X\subset\C^3$ containing $\orig$ but none of those extra critical points. (Note that upon restricting the domain in this way we do not alter the  Hochschild homology of the matrix factorizations since the natural dg functor ${\rm MF}_G(\C^3,W)\to{\rm MF}_G(X,W)$ is known to be a dg Morita equivalence.) 

The subtlety is that our new domain $X$ has to be as specified in Section \ref{SN}, i.e. within the range of applicability of the main results of this paper. That is, $X$ should be of the form
$
\C^3\setminus\bigcup_{i,j}\{x_i=\lambda_i^j\}
$ for some $\lambda_i^j\neq0$.
The union of hyperplanes that we remove should be $G$-stable and contain  
all the critical points we want to get rid of.  Such a configuration is easy to construct: we take as $\lambda_i^j$ all the non-zero coordinates of all those redundant critical points. The $G$-invariance of $W$ implies that this configuration is $G$-stable.

Corollary \ref{mf} is applicable now and reduces the theorem to the following claim:

\begin{proposition}
There is an isomorphism of $\ZZ$-graded  algebras
$
\Jacoii(X,W,G)^G\simeq{\rm H}^*({S},\C).
$
\end{proposition}
 
\noindent The remainder of this appendix is devoted to the proof of the proposition.

\vspace{0.1in}

\noindent\underline{\it $\Jacoii(X,W,G)^G$ as a $\ZZ$-graded space:} Let $\zv=(\zeta,\zeta,\zeta^{-2})$ be any cyclic generator of $G$. Let us rename the generators $\xi_g$ of the ``twisted sectors'' in $\Jacoii(X,W,G)$ as follows:
\[
\xi^+_k:=\xi_{\zv^{\, k}}, \quad \xi^-_k:=\xi_{\zv^{\, -k}}, \qquad  k=1,\ldots,\gen.
\]
Since the $g$-fixed locus of any $e\neq g\in G$ is the origin $\orig$, $d_g=3$ and $M(W^g)=\C$  in this case, so
\[
{\rm M}^{\rm even}(X,W,G)=M(W) \xi_e \quad {\rm M}^{\rm odd}(X,W,G)=\bigoplus_{k=1}^{\gen}\left(\C \xi^+_k\oplus\C \xi^-_k\right).
\]
It is easy to see that all the $\xi^\pm_k$ are $G$-invariant (see (\ref{Gvw}) for the definition of the $G$-action), hence
\[
{\rm M}^{\rm even}(X,W,G)^G=M(W)^G \xi_e \quad {\rm M}^{\rm odd}(X,W,G)^G=\bigoplus_{k=1}^{\gen}\left(\C \xi^+_k\oplus\C \xi^-_k\right).
\]
Let us compute $M(W)^G$. An important technical aspect is that we can now treat $W$ as a {\it local germ} at $\orig$ rather than a global function, since $\orig$ is the only critical point of $W$ in $X$.
One can show that $M(W)$ is spanned by the classes of $1$, $x_1x_2x_3$, and $x_i^l$ for $i=1,2,3$,  $l=1,\ldots, 2\gen$. Hence $M(W)^G=\C\oplus \C{\varphi}$ where  ${\varphi}$ is the class of $x_1x_2x_3$. 
Thus, finally,
\[
{\rm M}^{\rm even}(X,W,G)^G=\C \xi_e\oplus\C{\varphi}\xi_e \quad {\rm M}^{\rm odd}(X,W,G)^G=\bigoplus_{k=1}^{\gen}\left(\C \xi^+_k\oplus\C \xi^-_k\right).
\]

%\vspace{0.1in}

\noindent\underline{\it $\Jacoii(X,W,G)^G$ as an algebra:} Being the Hochschild cohomology of something,  $\Jacoii(X,W,G)^G$ is automatically super-commutative (Corollary \ref{cge}).
Furthermore, ${\varphi}^2=0$  since ${\varphi}^2$ is $G$-invariant and so lies in $\C\oplus\C\cdot {\varphi}$ but it has to be nilpotent (recall that we are dealing with the Milnor algebra of an isolated local singularity). Hence
\begin{equation}\label{nilp1}
{\varphi} \xi_e\cup{\varphi} \xi_e=0.
\end{equation}
Also, the class of $x_1x_2x_3$ vanishes in $M(W^g)=\C$ for $g\neq e$, hence
\begin{equation}\label{nilp2}
{\varphi} \xi_e\cup \xi^\pm_k= 0\quad \forall k.
\end{equation}
Next, since $d_g=3$ for all $g\neq e$, (\ref{non-int}) implies that
\begin{equation}\label{trivial}
\xi^+_k\cup \xi^+_l=\xi^-_k\cup \xi^-_l=0\quad \forall k,l; \qquad \xi^+_k\cup \xi^-_l=0\quad k\neq l.
\end{equation}
Hence the only interesting products are 
\begin{equation}\label{nontrivial}
\xi^+_k\cup \xi^-_k=\sigma_{\zv^{\, k},\zv^{\, -k}}\cdot \xi_e,\quad \sigma_{\zv^{\, k},\zv^{\, -k}}\in M(W)^G. 
\end{equation}
The shape of the formula (\ref{cgh}) -- namely, the fact that the formula contains only the second partial difference derivatives of $W$ -- suggests  that $ \sigma_{\zv^{\, k},\zv^{\, -k}}$ will have no constant term, i.e. 
\[
 \sigma_{\zv^{\, k},\zv^{\, -k}}=c_k {\varphi},\quad c_k\in \C.
\] 
We compute $\sigma_{\zv^{\, k},\zv^{\, -k}}$ below and, in particular, show that $c_k\neq 0$ for all $k$. This suffices to conclude the proof of the proposition. Indeed, as an abstract  super-commutative algebra, ${\rm H}^*({S},\C)$ has the form:
\[
{\rm H}^{\rm even}({S},\C)=\C\oplus \C\gamma,\quad {\rm H}^{\rm odd}({S},\C)=\bigoplus_{k=1}^{\gen}\left(\C\cdot \alpha_k\oplus\C\cdot \beta_k\right)
\]
with
\begin{eqnarray*}
&&\gamma\cdot \gamma=0,\qquad \gamma\cdot \alpha_k=\gamma\cdot \beta_k=0\quad \forall k\\
&&\alpha_k\cdot \alpha_l=\beta_k\cdot \beta_l=0\quad \forall k,l;\qquad  \alpha_k\cdot \beta_l=0\quad \forall k\neq l;\qquad \alpha_k\cdot \beta_k =\gamma\quad \forall k.
\end{eqnarray*}
Thus, by (\ref{nilp1}), (\ref{nilp2}), (\ref{trivial}), and (\ref{nontrivial}) the assignment
\[
\xi_e\mapsto 1,\quad \xi^+_k\mapsto \alpha_k,\quad \xi^-_k\mapsto  c_k\beta_k,\quad {\varphi} \xi_e\mapsto  \gamma
\]
extends  to an isomorphism $\Jacoii(X,W,G)^G\to{\rm H}^*({S},\C)$.

So the last step is to prove

\begin{lemma} $\sigma_{\zv^{\, k},\zv^{\, -k}}=\frac{1}{(1-\zeta^k)^2\cdot(1-\zeta^{-2k})}\ {\varphi}.$
\end{lemma}

We have
\begin{multline*}
  {\rm H}_W (x,y,z)=  \left(\frac{x_1^{2\gen+1}-y_1^{2\gen+1}}{x_1-y_1}-\frac{x_1^{2\gen+1}-z_1^{2\gen+1}}{x_1-z_1}\right)\frac{1 }{y_1-z_1}\theta_1 \otimes\theta_1+
  \\
   +\left(\frac{x_2^{2\gen+1}-y_2^{2\gen+1}}{x_2-y_2}-\frac{x_2^{2\gen+1}-z_2^{2\gen+1}}{x_2-z_2}\right)
\frac{1}{y_2-z_2}\theta_2 \otimes\theta_2+\\
  +  \left(\frac{x_3^{2\gen+1}-y_3^{2\gen+1}}{x_3-y_3}-\frac{x_3^{2\gen+1}-z_3^{2\gen+1}}{x_3-z_3}\right) \frac{1}{y_3-z_3}\theta_3 \otimes\theta_3- \\
-  x_3\ \theta_2 \otimes\theta_1 -  y_2\ \theta_3 \otimes\theta_1 -  z_1 \theta_3 \otimes\theta_2
\end{multline*}
and hence 
\begin{multline*}
{\rm H}_W (x,\zv^{\,k}(x),x)=  \frac{2\gen+1}{1-\zeta^k}  x_1^{2\gen-1}\theta_1 \otimes\theta_1+\frac{2\gen+1}{1-\zeta^k}  x_2^{2\gen-1}\theta_2 \otimes\theta_2+\frac{2\gen+1}{1-\zeta^{-2k}}  x_3^{2\gen-1}\theta_3 \otimes\theta_3
  - \\
-  x_3\ \theta_2 \otimes\theta_1 -  \zeta^k x_2\ \theta_3 \otimes\theta_1 -  x_1 \theta_3 \otimes\theta_2.
\end{multline*}
Also
\[
{\rm H}_{W,\zv^{\,k}}(x)=-\frac{\zeta^{k}}{1-\zeta^k} x_3\ \theta_1\theta_2-\frac{\zeta^{2k}}{1-\zeta^k} x_2\ \theta_1\theta_3,\quad {\rm H}_{W,\zv^{-\,k}}(\zv^{\,k}(x))=\frac{\zeta^{-2k}}{1-\zeta^k} x_3\ \theta_1\theta_2+\frac{1}{1-\zeta^k} x_2\ \theta_1\theta_3.
\]

We are looking for $\sigma_{\zv^{\, k},\zv^{\, -k}}$ which is the constant coefficient in the expression
\[
\frac16\LL\left(\left(\lfloor{\rm H}_W (x,\zv^{\,k}(x),x)\rfloor_{e}+\lfloor{\rm H}_{W,\zv^{\,k}}(x)\rfloor_{e}\otimes1+1\otimes\lfloor{\rm H}_{W,\zv^{-\,k}}(\zv^{\,k}(x))\rfloor_{e}\right)^3\otimes \partial_{\theta_1}\partial_{\theta_2}\partial_{\theta_3}\otimes \partial_{\theta_1}
\partial_{\theta_2}\partial_{\theta_3}\right).
\] 
Denoting 
\[
{\rm H}_W (x,\zv^{\,k}(x),x) = A_{11}\,\theta_1\otimes\theta_1
+A_{22} \,\theta_2\otimes\theta_2+A_{33} \,\theta_3\otimes\theta_3
+A_{21}\,\theta_{2}\otimes\theta_{1}+A_{31}\,\theta_{3}\otimes\theta_{1}+A_{32}\,\theta_{3}\otimes\theta_{2},
\]
and  
\[
{\rm H}_{W,\zv^{\,k}}(x)= B_{12}\, \theta_1\theta_{2}+B_{13}\, \theta_1\theta_{3},\quad  {\rm H}_{W,\zv^{-\,k}}(\zv^{\,k}(x)) = C_{12}\, \theta_1\theta_{2} + C_{13}\, \theta_1\theta_{3},
\]
one checks that 
\begin{multline*}
\sigma_{\zv^{\, k},\zv^{\, -k}}=\lfloor A_{11}A_{22}A_{33}-A_{22}B_{13}C_{13}-A_{33}B_{12}C_{12}+A_{32}B_{12}C_{13}\rfloor_{e}=\\
=\frac{(2\gen+1)^3}{(1-\zeta^k)^2(1-\zeta^{-2k})}\lfloor x_1^{2\gen-1}x_2^{2\gen-1}x_3^{2\gen-1}\rfloor_{e}
+\frac{(2\gen+1)\zeta^{2k}}{(1-\zeta^k)^3}\lfloor x_2^{2\gen+1}\rfloor_{e}
+
\\
+
\frac{(2\gen+1)\zeta^{-k}}{(1-\zeta^k)^2(1-\zeta^{-2k})}\lfloor x_3^{2\gen+1}\rfloor_{e}
+
\frac{\zeta^k}{(1-\zeta^k)^2}\lfloor x_1x_2x_3\rfloor_{e}.
\end{multline*}
The first term equals ${\varphi}^{2\gen-1}$ and hence vanishes. Also,
\[
(2\gen+1)\lfloor x_1^{2\gen}\rfloor_{e} = \lfloor x_2x_3 \rfloor_{e},\quad 
(2\gen+1)\lfloor x_2^{2\gen}\rfloor_{e} = \lfloor x_1x_3 \rfloor_{e},\quad
(2\gen+1)\lfloor x_3^{2\gen}\rfloor_{e} = \lfloor x_1x_2 \rfloor_{e}
\]
imply $\lfloor(2\gen+1)x_i^{2\gen+1}\rfloor_{e} = {\varphi}$. Thus,
\begin{equation*}
\sigma_{\zv^{\, k},\zv^{\, -k}}=
\frac{\zeta^{2k}}{(1-\zeta^k)^3}{\varphi}
+
\frac{\zeta^{-k}}{(1-\zeta^k)^2(1-\zeta^{-2k})}{\varphi}
+
\frac{\zeta^k}{(1-\zeta^k)^2}{\varphi}
=\frac{1}{(1-\zeta^k)^2(1-\zeta^{-2k})}\ {\varphi}.
\end{equation*}

\newpage

\section{ Proofs of intermediate results}\label{B}

\vspace{0.1in}

\subsubsection*{\bf Proof of Proposition \ref{Kunn}}\label{Kunnproof} 
Consider the $\Z$-graded $(A\otimes A')$-bimodule 
\[
\rD_*(A\otimes A'):=\rB_*(A)\otimes_{} \rB_*(A')
\]
and denote by $\delta_{\rm bar}^\rD$ and $\delta_{\rm curv}^\rD$ the differentials $\ddbar\otimes 1+1\otimes \ddbar$ and $\ddcurve\otimes 1+1\otimes \ddcurve$ thereon. Repeating the definitions from Section \ref{sect2.3} for the triple $(\rD_*(A\otimes A'),\delta_{\rm bar}^\rD, \delta_{\rm curv}^\rD)$ instead of $(\rB_*(A\otimes A'),\delta_{\rm bar}, \delta_{\rm curv})$, one obtains functors 
\[
\cHD^*(A\otimes A', W\otimes 1+1\otimes W';-),\quad \cHD_*(A\otimes A', W\otimes 1+1\otimes W';-)
\] on the category of $(A\otimes A', W\otimes 1+1\otimes W')$-bimodules which share all the properties of the Hochschild (co)homology. In particular, the new functors come with their cup and cap products defined just as $\cup_{\rm Hoch}$ and $\cap_{\rm Hoch}$ using  the coproduct
\[
\Delta_\rD: \rD_*(A\otimes A')\to \rD_*(A\otimes A')\otimes\rD_*(A\otimes A'),
\]
\[
\Delta_\rD(\underline{a}'\otimes \underline{a}'')=(-1)^{|\underline{a}'_{(2)}||\underline{a}''_{(1)}|}(\underline{a}'_{(1)}\otimes \underline{a}''_{(1)})\otimes (\underline{a}'_{(2)}\otimes \underline{a}''_{(2)})
\]
where $\underline{a}'\in \rB_*(A)$, $\underline{a}''\in \rB_*(A')$, and
$
\underline{a}'_{(1)}\otimes \underline{a}'_{(2)}:=\Delta_{\rm bar}(\underline{a}')
$, $\underline{a}''_{(1)}\otimes \underline{a}''_{(2)}:=\Delta_{\rm bar}(\underline{a}'')$.

One can formulate a version of  Proposition \ref{Kunn} with $\cHD^*(A\otimes A', W\otimes 1+1\otimes W';-)$ and $\cHD_*(A\otimes A', W\otimes 1+1\otimes W';-)$ on the right-hand sides of  (\ref{kunn}) instead of the Hochschild (co)homology. We claim that this version of the assertion does hold: The isomorphisms in this case are induced by the obvious maps
\begin{eqnarray}\label{2maps}
&&\text{\scalebox{0.95}{$
\Hom_{A\otimes {A}^{\rm op}}(\rB_*(A), M)\otimes \Hom_{A'\otimes {A'}^{\rm op}}(\rB_*(A'), M')\to \Hom_{(A\otimes A')\otimes {(A\otimes A')}^{\rm op}}(\rD_*(A\otimes A'), M\otimes M'),
$}}\nonumber
\\
&&\text{\scalebox{0.95}{$
(M\otimes_{A\otimes {A}^{\rm op}}\rB_*(A))\otimes (M'\otimes_{A'\otimes {A'}^{\rm op}}\rB_*(A'))\to (M\otimes M')\otimes_{(A\otimes A')\otimes {(A\otimes A')}^{\rm op}}\rD_*(A\otimes A').
$}}
\end{eqnarray}
Despite the seeming simplicity of the claim, there are two subtleties that require clarification. (In fact, it is this part of the proof that relies on the {\it smoothness} assumption.)
  
Firstly, we need to explain why (\ref{2maps}) induce quasi-isomorphisms of  the corresponding mixed complexes. For the second map it is straightforward: the map is an isomorphism even on the cochain level. The first map, on the other hand, is in general only an inclusion. However, as explained in the proof of Theorem 3.1 on page 210 in \cite{CE}, the assumption that $A$ and $A'$ admit resolutions by finitely generated projective bimodules suffices to claim that the induced map is still a quasi-isomorphism. 

Secondly, we should explain why the quasi-isomorphisms of the mixed complexes induced by (\ref{2maps}) yield isomorphisms of  $\K\bb$-modules
\begin{eqnarray}\label{2moremaps}
\cHH^*(A,W;M)\otimes_{\K\bb} \cHH^*(A',W';M') \simeq \cHD^*(A\otimes A', W\otimes 1+1\otimes W'; M\otimes M'),\nonumber\\
\cHH_*(A,W;M)\otimes_{\K\bb} \cHH_*(A',W';M') \simeq \cHD_*(A\otimes A', W\otimes 1+1\otimes W'; M\otimes M').
\end{eqnarray}
In general, the functor $(\rC,b,B)\mapsto
{\rm H}^*(\rC\bbb, b+ {\uu} B)$ is only {\it lax} monoidal, i.~e.  for generic  mixed complexes $(\rC, b, B)$ and $(\rC', b', B')$ the natural map
\begin{equation}\label{monoid}
{\rm H}^*(\rC\bbb, b+ {\uu} B)\otimes_{\K\bb} {\rm H}^*(\rC'\bbb, b'+ {\uu} B')\to {\rm H}^*((\rC\otimes\rC')\bbb, (b+ {\uu} B)\otimes1+1\otimes (b'+ {\uu} B'))
\end{equation}
need not to be an isomorphism. 

\vspace{0.1in}
\noindent {\bf Lemma.} {\it If ${\rm H}^n(\rC,b)=0$ and ${\rm H}^n(\rC',b')=0$ for $n\ll0$ 
then (\ref{monoid}) is an isomorphism.}

\vspace{0.1in}

\noindent To prove the lemma, we pick $n_0$ so that  ${\rm H}^n(\rC,b)=0$ and ${\rm H}^n(\rC',b')=0$ vanish in degrees $< n_0$ and consider the corresponding truncated  mixed complexes
$
(\tau\rC,  b, B) 
$
and 
$
(\tau\rC',  b', B') 
$
where
\[
\tau\rC_n:=\begin{cases} 0 & n<n_0\\ {\rm Coker} (b: \rC_{n_0-1}\to \rC_{n_0}) & n=n_0\\ \rC_n & n>n_0\end{cases}
\]
and the same for $\tau\rC'$.  Since the canonical projections  
\[
p: (\rC,  b, B)\to(\tau\rC,  b, B),\quad p': (\rC',  b', B')\to(\tau\rC',  b', B')
\]
are quasi-isomorphisms, their tensor product
\[
p\otimes p': (\rC,  b, B)\otimes (\rC',  b', B')\to (\tau\rC,  b, B)\otimes(\tau\rC',  b', B')
\]
is a quasi-isomorphism as well. Consider the commutative diagram
\[\begin{CD}
\text{\scalebox{0.9}{$
{\rm H}^*(\rC\bbb, b+ {\uu} B)\otimes_{\K\bb} {\rm H}^*(\rC'\bbb, b'+ {\uu} B')
$}} 
@>>>
\text{\scalebox{0.9}{$
{\rm H}^*((\rC\otimes\rC')\bbb, (b+ {\uu} B)\otimes1+1\otimes (b'+ {\uu} B'))
$}}\\
@VVV @VVV\\
\text{\scalebox{0.9}{$
{\rm H}^*(\tau\rC\bbb, b+ {\uu} B)\otimes_{\K\bb} {\rm H}^*(\tau\rC'\bbb, b'+ {\uu} B')
$}} 
@>>>
\text{\scalebox{0.9}{$
{\rm H}^*((\tau\rC\otimes\tau\rC')\bbb, (b+ {\uu} B)\otimes1+1\otimes (b'+ {\uu} B'))
$}}
\end{CD}\]
where the vertical maps are the isomorphisms induced by the above canonical projections and the horizontal maps come from the monoidal structure of the periodic cohomology functor. The claim of the lemma follows immediately from the observation that the lower horizontal map in the diagram is also an isomorphism -- this is a consequence of the isomorphisms 
\[
\tau\rC\bbb\simeq \tau\rC\otimes \K\bb,\quad \tau\rC'\bbb\simeq \tau\rC'\otimes \K\bb,\quad (\tau\rC\otimes\tau\rC')\bbb\simeq (\tau\rC\otimes\tau\rC')\otimes \K\bb
\]
which themselves follow from the fact that all the complexes involved are bounded below; see (\ref{Ct}). 

Note that the mixed complex $\rHH^*(A,W;M)$ always satisfies the condition of the lemma and the mixed complex $\rHH_*(A,W;M)$ satisfies this condition provided $A$ is smooth. Thus, (\ref{2moremaps}) are isomorphisms.

To conclude the proof of Proposition \ref{Kunn} it remains to establish a product-preserving equivalence between the (co)homology functors 
$\cHD^*(A\otimes A', W\otimes 1+1\otimes W';-)$ and $\cHD_*(A\otimes A', W\otimes 1+1\otimes W';-)$ and the ordinary Hochschild (co)homology functors $\cHH^*(A\otimes A', W\otimes 1+1\otimes W';-)$ and $\cHH_*(A\otimes A', W\otimes 1+1\otimes W';-)$

Consider the map of $A\otimes A'$-bimodules 
$
{\sf Sh}:\rD_*(A\otimes A')=\rB_*(A)\otimes_{} \rB_*(A')\to \rB_*(A\otimes A')
$
determined by 
\[
{\sf Sh}: (1[a_1|\ldots |a_n]1)\otimes (1[a'_1|\ldots |a'_m]1)\mapsto (1\otimes1){\sf sh}[a_1\otimes1|\ldots|a_n\otimes1|1\otimes a'_1|\ldots|1\otimes a'_m](1\otimes1)
\]
where ${\sf sh}$ stands for the sum over all the permutations that shuffle the $a$'s with the $a'$'s while preserving the order within the two groups, and each summand is multiplied by the sign of the corresponding permutation. It is a classical fact (cf. Section 6, Chapter XI in \cite{CE}) that ${\sf Sh}$ is compatible with the differential $\ddbar$ on the bar resolutions and induces a quasi-isomorphism of complexes. It also turns out to be compatible with $\ddcurve$; this is a more straightforward observation and we leave the proof to the reader. Thus, ${\sf Sh}$ gives rise to isomorphisms of functors
\begin{eqnarray*}
{\sf Sh}^*: \cHH^*(A\otimes A', W\otimes 1+1\otimes W';-)\to \cHD^*(A\otimes A', W\otimes 1+1\otimes W';-),\\
 {\sf Sh}_*: \cHD_*(A\otimes A', W\otimes 1+1\otimes W';-)\to \cHH_*(A\otimes A', W\otimes 1+1\otimes W';-).
\end{eqnarray*}
To show the compatibility of the isomorphisms with the products, it suffices to show that
the diagram
\[\begin{CD}
\rD_*(A\otimes A') @>{\Delta_\rD}>>\rD_*(A\otimes A')\otimes \rD_*(A\otimes A')\\
@VV{{\sf Sh}}V @VV{{\sf Sh}\otimes {\sf Sh}}V\\
\rB_*(A\otimes A') @>{\Delta_{\rm bar}}>>\rB_*(A\otimes A')\otimes \rB_*(A\otimes A')
\end{CD}\]
is commutative or,  denoting  $x_i:=a_i\otimes1$, $y_j:=1\otimes a'_j$, that
\begin{multline*}
\Delta_{\rm bar}\left({\sf sh}[x_1|\ldots|x_n|y_1|\ldots|y_m]\right)=
\\
=\sum_{i=0}^n\sum_{j=0}^m(-1)^{j(n-i)}
{\sf sh}[x_1|\ldots|x_i|y_1|\ldots|y_j]\otimes {\sf sh}[x_{i+1}|\ldots|x_n|y_{j+1}|\ldots|y_m].
\end{multline*}
The latter is an easy combinatorial exercise.

\subsubsection*{\bf Proof of Proposition \ref{Ginv}.}\label{123} The claim is a straightforward generalization of well-known results in the non-curved setting, so we will only sketch the proof. 

Consider the following two maps of graded spaces: 

\noindent (1)  
$\Xi^*: \rB^*(A,A\rtimes  G)^G\to\rB^*(A\rtimes  G,A\rtimes  G)$ defined as the restriction to the $G$-invariants of the map  
$
\rB^*(A,A\rtimes  G)\to\rB^*(A\rtimes  G,A\rtimes  G)
$ that sends  the element
$
1[a_1|a_2|\ldots |a_n]1\mapsto
\widetilde{D}(a_1|a_2|\ldots |a_n)\otimes g
$ 
 of $\rB^*(A,A\rtimes  G)$ to 
\begin{eqnarray*}
1[a_1\otimes g_1|a_2\otimes g_2|\ldots |a_n\otimes g_n]1\mapsto
\widetilde{D}(a_1|g_1(a_2)|\ldots |g_1g_2\ldots g_{n-1}(a_n))\otimes gg_1g_2\ldots g_n
\end{eqnarray*}

\noindent (2) $\Xi_*: \rB_*(A\rtimes  G,A\rtimes  G)\to\rB_*(A,A\rtimes  G)_G$ induced by the map 
$
\rB_*(A\rtimes  G,A\rtimes  G)\to\rB_*(A,A\rtimes  G)
$
given by
\begin{multline*}
(a\otimes g)\otimes \left(1[a_1\otimes g_1|a_2\otimes g_2|\ldots |a_n\otimes g_n]1\right)\mapsto\\
(a\otimes gg_1g_2\ldots g_n)\otimes \left(1[g_n^{-1}\ldots g_1^{-1}(a_1)|g_n^{-1}\ldots g_2^{-1}(a_2)|\ldots |g_n^{-1}(a_n)]1\right)
\end{multline*}
The maps $\Xi^*$ and $\Xi_*$ can be shown to define morphisms of mixed complexes 
\[
\Xi^*: \rHH^*(A,W; A\rtimes  G)^G\to\rHH^*(A\rtimes  G,W),\quad \Xi_*: \rHH_*(A\rtimes  G,W)\to\rHH_*(A,W;A\rtimes  G)_G,
\]
and the key observation is that they are quasi-isomorphisms; see \cite[Thm.5.4]{CGW} and \cite[Prop.8]{Ba}.  Finally, $\Xi^*$ and $\Xi_*$  are easily seen to be compatible with the cup and cap products in the following sense:
\[
\Xi^*(D_1\cup D_2)=\Xi^*(D_1)\cup \Xi^*(D_2),\quad \forall\, D_i\in \rB^*(A, A\rtimes  G)^G,
\]
\[
\Xi_*(\omega)\cap D=\Xi_*(\omega\cap\Xi^*(D)) \quad \forall \, \omega\in \rB_*(A\rtimes  G, A\rtimes  G), D\in \rB^*(A, A\rtimes  G)^G.
\]

\subsubsection*{\bf Proof of Proposition \ref{brcomm}.}\label{Gtwcomproof}
The idea is to show that the diagram of $\K\bb$-linear complexes
\[\xymatrixcolsep{5pc}
\xymatrix{
\rB^*(A,A\otimes g)\bbb\otimes_{\K\bb}\rB^*(A,A\otimes h)\bbb \ar[rd]^\cup \ar[d]^{({\rm id}\otimes h^{-1})\cdot\sigma} &\\
\rB^*(A,A\otimes h)\bbb\otimes_{\K\bb} \rB^*(A,A\otimes h^{-1}gh)\bbb \ar[r]^-\cup  & \rB^*(A,A\otimes gh)\bbb,
}
\]
where $\sigma$ denote the (graded) transposition of the terms, commutes {\it up to homotopy} or, in other words, that there exists a series $\sum_{i=i_0}^\infty \chi_it^i$, where $\chi_i$ is a degree $-2i-1$ linear operator from $\rB^*(A,A\otimes g)\otimes\rB^*(A,A\otimes h)$ to 
$\rB^*(A,A\otimes gh)$, such that
\[
\cup -\cup\cdot ({\rm id}\otimes h^{-1})\cdot \sigma 
=(\ddh+t\ddw)\sum_{i} \chi_i t^i  +\sum_{i} \chi_i t^i ((\ddh+t\ddw)\otimes1+1\otimes(\ddh+t\ddw)). 
\]
In fact, we will show that there exists such a series $\sum_i \chi_it^i$ with $\chi_i=0$ for $i\neq0$, i.e. we will construct a degree $-1$ operator $\chi=\chi_0: \rB^*(A,A\otimes g)\otimes\rB^*(A,A\otimes h)\to 
\rB^*(A,A\otimes gh)$ such that
\begin{eqnarray}\label{homotop}
\cup -\cup\cdot ({\rm id}\otimes h^{-1})\cdot \sigma 
&=&\ddh\cdot \chi +\chi\cdot (\ddh\otimes1+1\otimes\ddh),\nonumber \\
0 &=& \ddw\cdot \chi+\chi\cdot(\ddw\otimes1+1\otimes\ddw).
\end{eqnarray}

The homotopy $\chi$ is given by a slight modification of M. Gerstenhaber's well-known formula in the non-equivariant setting \cite[Thm.3]{Ge}. Namely, let us fix $D_1\in \rB^{l}(A,A\otimes g)$ and $D_2\in \rB^{m}(A,A\otimes h)$. So,
\begin{eqnarray}\label{Dtil}
D_1(a_0[a_1|\ldots |a_n]a_{n+1})&=&\begin{cases} a_0\cdot\widetilde{D}_1(a_1|\ldots |a_l)\cdot g(a_{l+1})\otimes g, & n=l\\ 0, & \text{otherwise}\end{cases},\\
D_2(a_0[a_1|\ldots |a_n]a_{n+1})&=&\begin{cases} a_0\cdot \widetilde{D}_2(a_1|\ldots |a_m)\cdot h(a_{m+1})\otimes h, & n=m\\ 0, & \text{otherwise}\end{cases}\nonumber
\end{eqnarray}
for some linear maps $\widetilde{D}_1: A^{\otimes l}\to A$ and $\widetilde{D}_2: A^{\otimes m}\to A$. (This notation will be used throughout the proof: $\widetilde{D}$ will stand for the map $A^{\otimes *}\to A$ associated with a Hochschild cochain $D\in \rB^{*}(A,A\otimes -)$ by the above rule.) We define 
$\chi(D_1,D_2)\in \rB^{l+m-1}(A,A\otimes gh)$
by
\begin{multline*}
\widetilde{\chi(D_1,D_2)}(a_1|\ldots |a_{l+m-1}):=\\
=\sum_{i=1}^{l}(-1)^{(i-1)(m-1)+l}
\widetilde{D}_1(a_1|\ldots|a_{i-1}|\widetilde{D}_2(a_i|\ldots |a_{i+m-1})|h(a_{i+m})|\ldots |h(a_{l+m-1})).
\end{multline*}
Let us compute $\partial\chi(D_1,D_2)+\chi(\partial D_1,D_2)+(-1)^l\chi(D_1,\partial D_2)$ for $\partial:=\ddh$.
We have
\begin{multline*}
\widetilde{\partial\chi(D_1,D_2)}(a_1|\ldots |a_{l+m})=
\\
=(-1)^{l+m-1}a_1\cdot \widetilde{\chi(D_1,D_2)}(a_2|\ldots |a_{l+m})-\widetilde{\chi(D_1,D_2)}(a_1|\ldots |a_{l+m-1})\cdot gh(a_{l+m})
\\
+\sum_{j=1}^{l+m-1}(-1)^{j+l+m-1}\widetilde{\chi(D_1,D_2)}(a_1|\ldots |a_{j}a_{j+1}|\ldots|a_{l+m}).
\end{multline*}
Unfolding the definitions,
\begin{multline*}
(-1)^{l+m-1}a_1\cdot \widetilde{\chi(D_1,D_2)}(a_2|\ldots |a_{l+m})=
\\
=\sum_{i=1}^{l}(-1)^{i(m-1)}
a_1\cdot \widetilde{D}_1(a_2|\ldots|a_{i}|\widetilde{D}_2(a_{i+1}|\ldots |a_{i+m})|h(a_{i+m+1})|\ldots |h(a_{l+m})),
\end{multline*}
\begin{multline*}
-\widetilde{\chi(D_1,D_2)}(a_1|\ldots |a_{l+m-1})\cdot gh(a_{l+m})=
\\
=-\sum_{i=1}^{l}(-1)^{(i-1)(m-1)+l}
\widetilde{D}_1(a_1|\ldots|a_{i-1}|\widetilde{D}_2(a_i|\ldots |a_{i+m-1})|h(a_{i+m})|\ldots |h(a_{l+m-1}))\cdot gh(a_{l+m})
\end{multline*}
and
\begin{multline*}
\sum_{j=1}^{l+m-1}(-1)^{j+l+m-1}\widetilde{\chi(D_1,D_2)}(a_1|\ldots |a_{j}a_{j+1}|\ldots|a_{l+m})=
\\
=
\sum_{j=m+1}^{l+m-1}\sum_{i=1}^{j-m}(-1)^{i(m-1)+j}\widetilde{D}_1(a_1|\ldots|a_{i-1}|\widetilde{D}_2(a_i|\ldots |a_{i+m-1})|\ldots|h(a_{j}a_{j+1})|\ldots |h(a_{l+m}))
\\
+\sum_{j=m}^{l+m-1}\sum_{i=j-m+1}^{j}(-1)^{i(m-1)+j}
\widetilde{D}_1(a_1|\ldots|a_{i-1}|\widetilde{D}_2(a_i|\ldots|a_{j}a_{j+1}|\ldots |a_{i+m})|h(a_{i+m+1})|\ldots |h(a_{l+m}))
\\
+\sum_{j=1}^{l-2}\sum_{i=j+2}^{l}(-1)^{(i-1)(m-1)+j}
\widetilde{D}_1(a_1|\ldots|a_{j}a_{j+1}|\ldots|a_{i-1}|\widetilde{D}_2(a_i|\ldots |a_{i+m-1})|h(a_{i+m})|\ldots |h(a_{l+m})).
\end{multline*}
Furthermore, 
\begin{multline*}
\widetilde{\chi(\partial D_1,D_2)}(a_1|\ldots |a_{l+m})=
\\
=\sum_{i=1}^{l+1}(-1)^{(i-1)(m-1)+l+1}
\widetilde{\partial D}_1(a_1|\ldots|a_{i-1}|\widetilde{D}_2(a_i|\ldots |a_{i+m-1})|h(a_{i+m})|\ldots |h(a_{l+m}))=
\\
=-\widetilde{D}_2(a_1|\ldots |a_{m})\cdot \widetilde{D}_1(h(a_{1+m})|\ldots |h(a_{l+m}))
+(-1)^{lm}\widetilde{D}_1(a_1|\ldots|a_{l})\cdot g\left(\widetilde{D}_2(a_{l+1}|\ldots |a_{l+m})\right)
\\
-\sum_{i=2}^{l+1}(-1)^{(i-1)(m-1)}
a_1\cdot \widetilde{D}_1(a_2|\ldots|a_{i-1}|\widetilde{D}_2(a_i|\ldots |a_{i+m-1})|h(a_{i+m})|\ldots |h(a_{l+m}))
\\
+\sum_{i=1}^{l}(-1)^{(i-1)(m-1)+l}
\widetilde{D}_1(a_1|\ldots|a_{i-1}|\widetilde{D}_2(a_i|\ldots |a_{i+m-1})|h(a_{i+m})|\ldots |h(a_{l+m-1}))\cdot gh(a_{l+m})
\\
-\sum_{i=3}^{l+1}\sum_{j=1}^{i-2}(-1)^{(i-1)(m-1)+j}
(-1)^j\widetilde{D}_1(a_1|\ldots|a_ja_{j+1}|\ldots|a_{i-1}|\widetilde{D}_2(a_i|\ldots |a_{i+m-1})|h(a_{i+m})|\ldots |h(a_{l+m}))
\\
-\sum_{i=2}^{l+1}(-1)^{(i-1)m}
\widetilde{D}_1(a_1|\ldots|a_{i-1}\cdot\widetilde{D}_2(a_i|\ldots |a_{i+m-1})|h(a_{i+m})|\ldots |h(a_{l+m}))
\\
+\sum_{i=1}^{l}(-1)^{(i-1)m}
\widetilde{D}_1(a_1|\ldots|a_{i-1}|\widetilde{D}_2(a_i|\ldots |a_{i+m-1})\cdot h(a_{i+m})|\ldots |h(a_{l+m}))
\\
-\sum_{i=1}^{l-1}\sum_{j=i+m}^{l+m-1}(-1)^{i(m-1)+j}
\widetilde{D}_1(a_1|\ldots|a_{i-1}|\widetilde{D}_2(a_i|\ldots |a_{i+m-1})|h(a_{i+m})|\ldots|h(a_ja_{j+1})|\ldots |h(a_{l+m})).
\end{multline*}

\smallskip
Finally, 
\begin{multline*}
(-1)^l\widetilde{\chi( D_1,\partial D_2)}(a_1|\ldots |a_{l+m})=
\\
=\sum_{i=1}^{l}(-1)^{(i-1)m}
\widetilde{D}_1(a_1|\ldots|a_{i-1}|\widetilde{\partial D}_2(a_i|\ldots |a_{i+m})|h(a_{i+m+1})|\ldots |h(a_{l+m}))=
\\
=\sum_{i=1}^{l}(-1)^{im}
\widetilde{D}_1(a_1|\ldots|a_{i-1}|a_i\cdot \widetilde{D}_2(a_{i+1}|\ldots |a_{i+m})|h(a_{i+m+1})|\ldots |h(a_{l+m}))
\\
-\sum_{i=1}^{l}(-1)^{(i-1)m}
\widetilde{D}_1(a_1|\ldots|a_{i-1}|\widetilde{D}_2(a_i|\ldots |a_{i+m-1})\cdot h(a_{i+m})|h(a_{i+m+1})|\ldots |h(a_{l+m}))
\\
-\sum_{i=1}^{l}\sum_{j=i}^{i+m-1}(-1)^{i(m-1)+j}
\widetilde{D}_1(a_1|\ldots|a_{i-1}|\widetilde{D}_2(a_i|\ldots|a_ja_{j+1}|\ldots |a_{i+m})|h(a_{i+m+1})|\ldots |h(a_{l+m}))
\end{multline*}

Summing all the above equalities results in 
\begin{multline*}
\widetilde{\partial\chi(D_1,D_2)}(a_1|\ldots |a_{l+m})+\widetilde{\chi(\partial D_1,D_2)}(a_1|\ldots |a_{l+m})+
(-1)^l\widetilde{\chi( D_1,\partial D_2)}(a_1|\ldots |a_{l+m})=
\\
=-\widetilde{D}_2(a_1|\ldots |a_{m})\cdot \widetilde{D}_1(h(a_{1+m})|\ldots |h(a_{l+m}))
+(-1)^{lm}\widetilde{D}_1(a_1|\ldots|a_{l})\cdot g\left(\widetilde{D}_2(a_{l+1}|\ldots |a_{l+m})\right)
\end{multline*}
which is equivalent to the first equality in (\ref{homotop}). The proof of the second equality  is easier and is left to the reader.

\subsubsection*{\bf Proof of Proposition \ref{main}.}\label{mainproof}
Both equalities in (\ref{morph}) are proven by a direct calculation. To avoid {\it too} long formulas, we are including calculations in a special case, hoping that the reader will see the patterns.
 
Let us calculate $\ddko(\Psi(f_0[f_1|f_2|f_3]f_{4}))$. 
To begin with, since $\theta_i^2=0$ for all $i$, we can replace the strict inequalities in the index set in the right-hand side of  (\ref{psi}) by the non-strict ones, that is
\begin{equation*}
\ddko(\Psi(f_0[f_1|f_2|f_3]f_{4}))=
\sum_{1\leq j_1{\leq}j_2{\leq} j_3\leq N}l_1(f_0)\,\ddko\left(\nabla_{j_1}(f_1)\theta_{j_1}\,\nabla_{j_2}(f_2)\theta_{j_2}\, \nabla_{j_3}(f_3)\theta_{j_3}\right)\,l_{N+1}(f_4)
\end{equation*}
The commutation relations in ${\rm Cl}_N$ imply
\begin{multline*}
\ddko(\Psi(f_0[f_1|f_2|f_3]f_{4}))=\\
=\sum_{1\leq j_1{\leq}j_2{\leq} j_3\leq N}
l_1(f_0)\nabla_{j_1}(f_1)(x_{j_1}-y_{j_1})\, \nabla_{j_2}(f_2)\theta_{j_2}\,\nabla_{j_3}(f_3)\theta_{j_3}l_{N+1}(f_{4})-\\
-\sum_{1\leq j_1{\leq}j_2{\leq} j_3\leq N}l_1(f_0)\nabla_{j_1}(f_1)\theta_{j_1}\,\nabla_{j_2}(f_2)(x_{j_2}-y_{j_2})\,\nabla_{j_3}(f_3)\theta_{j_3}l_{N+1}(f_{4})+\\
+\sum_{1\leq j_1{\leq}j_2{\leq} j_3\leq N}l_1(f_0)\nabla_{j_1}(f_1)\theta_{j_1}\,\nabla_{j_2}(f_2)\theta_{j_2}\,\nabla_{j_3}(f_3)(x_{j_3}-y_{j_3})l_{N+1}(f_{4})
\end{multline*}
which by (\ref{deltai})  equals
\begin{multline*}
\sum_{1\leq j_1{\leq}j_2{\leq} j_3\leq N}
l_1(f_0)(l_{j_1}(f_1)-l_{j_1+1}(f_1))\, \nabla_{j_2}(f_2)\theta_{j_2}\,\nabla_{j_3}(f_3)\theta_{j_3}l_{N+1}(f_{4})-\\
-\sum_{1\leq j_1{\leq}j_2{\leq} j_3\leq N}l_1(f_0)\nabla_{j_1}(f_1)\theta_{j_1}\,(l_{j_2}(f_2)-l_{j_2+1}(f_2))\,\nabla_{j_3}(f_3)\theta_{j_3}l_{N+1}(f_{4})+\\
+\sum_{1\leq j_1{\leq}j_2{\leq} j_3\leq N}l_1(f_0)\nabla_{j_1}(f_1)\theta_{j_1}\,\nabla_{j_2}(f_2)\theta_{j_2}\,(l_{j_3}(f_3)-l_{j_3+1}(f_3))l_{N+1}(f_{4}).
\end{multline*}
Observe that 
\begin{eqnarray}\label{suml}
\sum_{j=\alpha}^{\beta}(l_{j}(f)-l_{j+1}(f))=l_{\alpha}(f)-l_{\beta+1}(f),
\end{eqnarray} 
so the previous expression equals
\begin{multline*}
\sum_{1\leq j_2{\leq} j_3\leq N}
l_1(f_0)(l_1(f_1)-l_{j_2+1}(f_1))\, \nabla_{j_2}(f_2)\theta_{j_2}\,\nabla_{j_3}(f_3)\theta_{j_3}l_{N+1}(f_{4})-\\
-\sum_{1\leq j_1{\leq} j_3\leq N}l_1(f_0)\nabla_{j_1}(f_1)\theta_{j_1}\,(l_{j_1}(f_2)-l_{j_3+1}(f_2))\,\nabla_{j_3}(f_3)\theta_{j_3}l_{N+1}(f_{4})+\\
+\sum_{1\leq j_1{\leq}j_2\leq N}l_1(f_0)\nabla_{j_1}(f_1)\theta_{j_1}\,\nabla_{j_2}(f_2)\theta_{j_2}\,(l_{j_2}(f_3)-l_{N+1}(f_3))l_{N+1}(f_{4})
\end{multline*}
or, renaming the indices,
\begin{multline*}
\sum_{1\leq \alpha{\leq} \beta\leq N}
l_1(f_0)(l_1(f_1)-l_{\alpha+1}(f_1))\, \nabla_{\alpha}(f_2)\theta_{\alpha}\,\nabla_{\beta}(f_3)\theta_{\beta}l_{N+1}(f_{4}) -\\
-\sum_{1\leq \alpha{\leq} \beta\leq N}l_1(f_0)\nabla_{\alpha}(f_1)\theta_{\alpha}\,(l_{\alpha}(f_2)-l_{\beta+1}(f_2))\,\nabla_{\beta}(f_3)\theta_{\beta}l_{N+1}(f_{4})+\\
+\sum_{1\leq \alpha{\leq}\beta\leq N}l_1(f_0)\nabla_{\alpha}(f_1)\theta_{\alpha}\,\nabla_{\beta}(f_2)\theta_{\beta}\,(l_{\beta}(f_3)-l_{N+1}(f_3))l_{N+1}(f_{4}).
\end{multline*}
The latter can be simplified by regrouping the summands  
\begin{multline*}
=\sum_{1\leq \alpha{\leq} \beta\leq N}
l_1(f_0)l_1(f_1)\, \nabla_{\alpha}(f_2)\theta_{\alpha}\,\nabla_{\beta}(f_3)\theta_{\beta}l_{N+1}(f_{4})
\\
-\sum_{1\leq \alpha{\leq} \beta\leq N}
l_1(f_0)(l_{\alpha+1}(f_1)\, \nabla_{\alpha}(f_2)+\nabla_{\alpha}(f_1)\,l_{\alpha}(f_2))\theta_{\alpha}\,\nabla_{\beta}(f_3)\theta_{\beta}l_{N+1}(f_{4})
\\
+\sum_{1\leq \alpha{\leq} \beta\leq N}l_1(f_0)\nabla_{\alpha}(f_1)\theta_{\alpha}\,(l_{\beta+1}(f_2)\,\nabla_{\beta}(f_3)+\nabla_{\beta}(f_2)\,l_{\beta}(f_3))\theta_{\beta}l_{N+1}(f_{4})
\\
-\sum_{1\leq \alpha{\leq}\beta\leq N}l_1(f_0)\nabla_{\alpha}(f_1)\theta_{\alpha}\,\nabla_{\beta}(f_2)\theta_{\beta}\,l_{N+1}(f_3)l_{N+1}(f_{4})
\end{multline*}
and employing the obvious equality
\begin{eqnarray}\label{leib}
\nabla_i(fg)=\nabla_i(f)l_{i}(g)+l_{i+1}(f)\nabla_i(g)
\end{eqnarray}
which yields
\begin{multline*}
=\sum_{1\leq \alpha{\leq} \beta\leq N}
l_1(f_0f_1)\, \nabla_{\alpha}(f_2)\theta_{\alpha}\,\nabla_{\beta}(f_3)\theta_{\beta}l_{N+1}(f_{4})
\\
-\sum_{1\leq \alpha{\leq} \beta\leq N}
l_1(f_0)\nabla_{\alpha}(f_1f_2)\theta_{\alpha}\,\nabla_{\beta}(f_3)\theta_{\beta}l_{N+1}(f_{4})
\\
+\sum_{1\leq \alpha{\leq} \beta\leq N}l_1(f_0)\nabla_{\alpha}(f_1)\theta_{\alpha}\,\nabla_{\beta}(f_2f_3)\theta_{\beta}l_{N+1}(f_{4})
\\
-\sum_{1\leq \alpha{\leq}\beta\leq N}l_1(f_0)\nabla_{\alpha}(f_1)\theta_{\alpha}\,\nabla_{\beta}(f_2)\theta_{\beta}\,l_{N+1}(f_3f_{4})
\end{multline*}
which is precisely $\Psi(\ddbar(f_0[f_1|f_2|f_3]f_{4}))$.

The second equality in (\ref{morph}) is easier to prove. Let us again demonstrate it in a special case:
\begin{multline*}
\ddcurve(\Psi(f_0[f_1|f_2|f_3]f_{4}))=
\sum_{1\leq j_1{<}j_2{<} j_3\leq N}l_1(f_0)\ddcurve\left(\nabla_{j_1}(f_1)\theta_{j_1}\,\nabla_{j_2}(f_2)\theta_{j_2}\, \nabla_{j_3}(f_3)\theta_{j_3}\right)l_{N+1}(f_4)=
\\
=\sum_{1\leq j_1{<}j_2{<} j_3\leq N}\sum_{j=1}^Nl_1(f_0)\nabla_{j}(W)\theta_{j}\nabla_{j_1}(f_1)\theta_{j_1}\,\nabla_{j_2}(f_2)\theta_{j_2}\, \nabla_{j_3}(f_3)\theta_{j_3}l_{N+1}(f_4). 
\end{multline*}
The sum over $j$ splits into $\sum_{1\leq j<j_1}+\sum_{j_1<j<j_2}+\sum_{j_2<j<j_3 }+\sum_{j_3<j\leq N}$ and we get
\begin{multline*}
\sum_{1\leq j< j_1{<}j_2{<} j_3\leq N}l_1(f_0)\nabla_{j}(W)\theta_{j}\nabla_{j_1}(f_1)\theta_{j_1}\,\nabla_{j_2}(f_2)\theta_{j_2}\, \nabla_{j_3}(f_3)\theta_{j_3}l_{N+1}(f_4)-
\\
-\sum_{1\leq j_1{<j<}j_2{<} j_3\leq N}l_1(f_0)\nabla_{j_1}(f_1)\theta_{j_1}\,\nabla_{j}(W)\theta_{j}\nabla_{j_2}(f_2)\theta_{j_2}\, \nabla_{j_3}(f_3)\theta_{j_3}l_{N+1}(f_4)+
\\
+\sum_{1\leq j_1{<}j_2{<} j < j_3\leq N}l_1(f_0)\nabla_{j_1}(f_1)\theta_{j_1}\,\nabla_{j_2}(f_2)\theta_{j_2}\, \nabla_{j}(W)\theta_{j}\nabla_{j_3}(f_3)\theta_{j_3}l_{N+1}(f_4)-
\\
-\sum_{1\leq j_1{<}j_2{<} j_3<j\leq N}l_1(f_0)\nabla_{j_1}(f_1)\theta_{j_1}\,\nabla_{j_2}(f_2)\theta_{j_2}\, \nabla_{j_3}(f_3)\theta_{j_3}\nabla_{j}(W)\theta_{j}l_{N+1}(f_4)
\end{multline*}
which is precisely $\Psi(\ddcurve(f_0[f_1|f_2|f_3]f_{4}))$.

Finally, to prove the last claim in the proposition, we recall that there is the following well-known
quasi-isomorphism of complexes of bimodules $i: (\rK_*(\K[{X}]),\ddko)\to(\rB_*(\K[{X}]),\ddbar)$
\begin{equation}\label{mori}
 i: f_0(x)f_1(y)\theta_{k_1}\ldots \theta_{k_n}\mapsto \sum_{\sigma\in S_n} {\rm sgn}(\sigma) \cdot f_0[x_{k_{\sigma(1)}}|\ldots|x_{k_{\sigma(n)}}]f_1.
\end{equation}
It remains to notice that $\Psi\cdot i={\rm id}_{\rK_*(\K[{X}])}$. 

\subsubsection*{\bf Proof of Proposition \ref{coprod}.}\label{coprodproof}
To shorten formulas, we will identify ${\rm Cl}_N^{\otimes 2}$ with ${\rm Cl}_{2N}=\K\langle \theta_i, \eta_i,{\partial_{\theta_i}}, {\partial_{ \eta_i}}\rangle$ via $\theta_i\otimes1\mapsto \theta_i$, $1\otimes \theta_i\mapsto  \eta_i$, etc.  
In the computation below we  use the fact that    
$
{\rm H}_{W}(x,y,z)=\sum_{i,j=1}^N \nabla^{y\to(y,z)}_j\nabla^{x\to(x,y)}_i(W)\,\theta_i\eta_j
$ because $\nabla^{y\to(y,z)}_j\nabla^{x\to(x,y)}_i=0$ for $i\leq j$. 

Abbreviating 
$
(\delta\otimes1+1\otimes \delta)\cdot\Delta-\Delta\cdot \delta
$
to $[\delta,\Delta]$, the claim is that
\begin{equation}\label{coalg}
[\ddko,\Delta_0]=0,\quad [\ddcurve,\Delta_{-2l}]+[\ddko,\Delta_{-2l-2}]=0,\quad \forall\, l
\end{equation}
According to our notation, 
$
\ddko=\sum_{i=1}^N(x_i-z_i)\partial_{\theta_i}
$
and
\[
\ddko\otimes1+1\otimes \ddko=\sum_{i=1}^N\left((x_i-y_i)\partial_{\theta_i}+(y_i-z_i)\partial_{ \eta_i}\right).
\] 
The first equality in (\ref{coalg}) follows from 
these formulas and the following obvious relation: 
\[
((x_i-y_i)\partial_{\theta_i}+(y_i-z_i)\partial_{ \eta_i})\cdot\Delta _0=\Delta_0\cdot (x_i-z_i)\partial_{\theta_i}\quad \forall\, i.
\]

Furthermore, 
$
\ddcurve=\sum_{i=1}^N\nabla^{x\to(x,z)}_i(W)\cdot\theta_i
$
and therefore $
\Delta_0\cdot \ddcurve=(\delta'_{\rm curv}\otimes1+1\otimes \delta'_{\rm curv})\cdot \Delta_0
$
where
\[
\delta'_{\rm curv}\otimes1+1\otimes \delta'_{\rm curv}=\sum_{i=1}^N\left(\nabla^{x\to(x,z)}_i(W)\cdot\theta_i+\nabla^{x\to(x,z)}_i(W) \eta_i\right).
\]
This observation, together with the first equality in (\ref{coalg}) and the fact that ${\rm H}_{W}$ is even (hence central  in $\K[X]^{\otimes 3}[\theta, \eta]$), reduces the proof of the remaining relations in (\ref{coalg}) to proving the following equality:
\begin{equation*}
(\ddcurve\otimes1+1\otimes \ddcurve)-(\delta'_{\rm curv}\otimes1+1\otimes \delta'_{\rm curv})=-(\ddko\otimes1+1\otimes \ddko)({\rm H}_{W})
\end{equation*}
or, in our new notation,
\begin{multline}\label{3}
\sum_{i=1}^N\left((\nabla^{x\to(x,y)}_i(W)-\nabla^{x\to(x,z)}_i(W))\theta_i+(\nabla^{y\to(y,z)}_i(W)-
\nabla^{x\to(x,z)}_i(W)) \eta_i\right)=
\\
=
-\sum_{j=1}^N\left((x_j-y_j)\partial_{\theta_j}+(y_j-z_j)\partial_{ \eta_j}\right)({\rm H}_{W}).
\end{multline}
Obviously, for any collection of polynomials $\{W_{\alpha\beta}\}$
\begin{equation*}
\sum_{j=1}^N\left((x_j-y_j)\partial_{\theta_j}+(y_j-z_j)\partial_{ \eta_j}\right)\sum_{\alpha,\beta=1}^N W_{\alpha\beta}\theta_\alpha \eta_\beta=
\sum_{j,\beta=1}^N(x_j-y_j) W_{j\beta} \eta_\beta-\sum_{j,\alpha=1}^N(y_j-z_j)W_{\alpha j}\theta_\alpha.
\end{equation*} 
Thus, setting $W_{\alpha\beta}=\nabla^{y\to(y,z)}_\beta\nabla^{x\to(x,y)}_\alpha(W)$, (\ref{3}) amounts to 
\begin{eqnarray*}
\nabla^{x\to(x,y)}_i(W)-\nabla^{x\to(x,z)}_i(W)=\sum_{j=1}^N(y_j-z_j)\nabla^{y\to(y,z)}_j\nabla^{x\to(x,y)}_i(W),\nonumber\\ 
\nabla^{x\to(x,z)}_i(W)-\nabla^{y\to(y,z)}_i(W)=\sum_{j=1}^N(x_j-y_j)\nabla^{y\to (y,z)}_i\nabla^{x\to(x,y)}_j(W).
\end{eqnarray*}
The first equality follows from (\ref{diffder}). The second one will also be a consequence of (\ref{diffder}) once we use the formula
\begin{equation*}
\nabla^{y\to(y,z)}_{i}\nabla^{x\to(x,y)}_{j}=\nabla^{x\to(x,y)}_{j}\nabla^{x\to(x,z)}_{i} 
\end{equation*}
which is easy to check by applying both hand sides to monomials in $x$.

\subsubsection*{\bf Proof of Proposition \ref{comdiag}.}\label{comdiagproof}
The left-hand side of (\ref{2}) is a series (in fact, a polynomial) in $t$. Let us denote the coefficients of this series by $\psi_i$:
\[
(\Psi\otimes\Psi)\cdot \Delta_{\rm bar}-\Delta_{\rm Kos}\cdot\Psi=:\sum_{i=0}^\infty \psi_it^i,\quad |\psi_i|=-2i.
\]
We will view the $\psi_i$s as elements of the complex of morphisms of complexes of $\K[X]$-bimodules:
\begin{equation}\label{cofm}
\Hom^*\left(\, (\rB_{*}(\K[X]),\ddbar)\,,\, (\rK_*(\K[X])\otimes_{\K[X]}\rK_*(\K[X]), \ddko\otimes1+1\otimes\ddko)\,\right).
\end{equation}
Since ${\rm H}^*(\rK_*(\K[X])\otimes_{\K[X]}\rK_*(\K[X]), \ddko\otimes1+1\otimes\ddko)\simeq \K[X]$, the complex (\ref{cofm}) is quasi-isomorphic to the complex $
\Hom^*\left((\rB_{*}(\K[X]),\ddbar),\K[X]\,\right)$, i.~e. to the  Hochschild cochain complex of $\K[X]$. Thus, (\ref{cofm}) has non-trivial cohomology groups only in non-negative degrees. This observation and the fact that the degress of the $\psi_i$s are non-positive imply that in order to prove the existence of the $h_i$s in (\ref{2}), it would suffice to show that the ``constant term'' $\psi_0$ defines a trivial class in the cohomology of (\ref{cofm}).  Let us explain why this would be enough.

Let us denote the differential in (\ref{2}) by $\widehat{\delta}$ 
\[
\widehat{\delta}(\psi)=(\ddko\otimes1+1\otimes\ddko)\cdot \psi-(-1)^{|\psi|}\psi\cdot \ddbar.
\]
and the second -- associated with $\ddcurve$ -- differential (of degree $-1$) on this complex by $\widehat{\delta}_{\rm curv}$: 
\[
\widehat{\delta}_{\rm curv}(\psi)=(\ddcurve\otimes1+1\otimes\ddcurve)\cdot \psi-(-1)^{|\psi|}\psi\cdot \ddcurve.
\]
Note that the two differential anti-commute, so $\widehat{\delta}+t\widehat{\delta}_{\rm curv}$ squares to 0. 

By Propositions \ref{main} and \ref{coprod} 
\[
(\widehat{\delta}+t\widehat{\delta}_{\rm curv})(\sum_{i=0}^\infty \psi_it^i)=0.
\]
We are looking for a series $\sum_{i=0}^\infty h_it^i$ such that 
\[
\sum_{i=0}^\infty \psi_it^i=(\widehat{\delta}+t\widehat{\delta}_{\rm curv})(\sum_{i=0}^\infty h_it^i).
\]
Assume we can find a degree $-1$ bimodule map $h_0$ such that $\psi_0=\widehat{\delta}(h_0)$ and consider the series 
\[
\sum_{i=1}^\infty \psi'_i\,t^i:= \sum_{i=0}^\infty \psi_it^i-(\widehat{\delta}+t\widehat{\delta}_{\rm curv})(h_0).
\]
It also satisfies 
$
(\widehat{\delta}+t\widehat{\delta}_{\rm curv})(\sum_{i} \psi'_i\,t^i)=0
$
but $|\psi'_1|=-2$ so, by the above observation, $\psi'_1$ defines a trivial class in the cohomology of (\ref{cofm}), i.~e. there exists a degree $-3$ bimodule map $h'_1$ such that $\psi'_1=\widehat{\delta}(h'_1)$, etc.    

So, to finish the proof of the proposition, it remains to prove that the class of $\psi_0$ is trivial, i.~e. that there exists an $h_0$ as above. By Proposition \ref{main} $\Psi: (\rB_*(\K[{X}]),\ddbar)\to (\rK_*(\K[{X}]),\ddko)$ is a quasi-isomorphism of complexes of $\K[X]$-bimodules. Since both complexes are K-projective, $\Psi$ is a homotopy equivalence. Moreover, we already know its right homotopy inverse: it is  the quasi-isomorphism $i: (\rK_*(\K[{X}]),\ddko)\to(\rB_*(\K[{X}]),\ddbar)$  defined in (\ref{mori}). By a standard general argument, $i$ is also a left homotopy inverse of $\Psi$, i.~e. $i\cdot \Psi={\rm id}_{\rB_*(\K[{X}])}-\ddbar\cdot h-h\cdot \ddbar$ for some degree $-1$ $\K[X]$-bimodule map $h:\rB_*(\K[{X}])\to\rB_*(\K[{X}])$.

It is easy to show that $\Delta_{\rm bar}\cdot i=(i\otimes i)\cdot \Delta_0$ where $\Delta_0$ is the morphism  defined in (\ref{del0}). Since $\psi_0=(\Psi\otimes\Psi)\cdot \Delta_{\rm bar}-\Delta_0\cdot\Psi$, we get
\[
\psi_0\cdot i=(\Psi\otimes\Psi)\cdot \Delta_{\rm bar}\cdot i-\Delta_0\cdot\Psi\cdot i=(\Psi\otimes\Psi)\cdot (i\otimes i)\cdot \Delta_0-\Delta_0=0.
\]
Then
\begin{multline*}
\psi_0=\psi_0\cdot {\rm id}_{\rB_*(\K[{X}])}=\psi_0\cdot(i\cdot \Psi+\ddbar\cdot h+h\cdot \ddbar)=
(\psi_0\cdot i)\cdot \Psi+\psi_0\cdot\ddbar\cdot h+\psi_0\cdot h\cdot\ddbar=
\\=
(\ddbar\otimes1+1\otimes\ddbar)\cdot (\psi_0\cdot h)+(\psi_0\cdot h)\cdot \ddbar=
\widehat{\delta}(\psi_0\cdot h).
\end{multline*}

\subsubsection*{\bf Proof of Lemma \ref{ethg}.}\label{ethgproof}
 The condition (\ref{conj}) is equivalent to 
 \[
 \ddkg+t\ddwg'=e^{t  {\rm H}_{{W,g}}}\cdot\ddkg\cdot e^{-t\cdot {\rm H}_{{W,g}}}
 \] which, in turn,  is easily seen to be equivalent to
\begin{eqnarray}\label{dkoshg}
\sum_{i\in I_{g}}(1-g_i)x_i\,\partial_{\theta_i}({\rm H}_{{W,g}})=\sum_{i\in I_{g}}\nabla^{x\to(x,g(x))}_i(W)\,\theta_i
\end{eqnarray}
where both hand sides are viewed as elements of $\K[X][\theta]$. Let us show that the element (\ref{hg}) satisfies the latter condition.

In the following calculation $W_i:=\nabla^{x\to(x,g(x))}_i(W)$, the indices $\alpha,i,j$ belong to $I_{g}$, and $l_i$, $\nabla_i$ stand for $l_i^{x\to(x,x^g)}$, $\nabla^{x\to(x,x^g)}_i$, respectively. 
We have:
\begin{multline*}
\sum_{\alpha}(1-g_\alpha)x_\alpha\,\partial_{\theta_\alpha}({\rm H}_{{W,g}})=
\sum_{\alpha}\sum_{j<i}\frac{1-g_\alpha}{1-g_j}x_\alpha\nabla_j(W_i)\partial_{\theta_\alpha}(\theta_j\,\theta_i)=\\
=\sum_{j<i}x_j\nabla_j(W_i)\theta_i-\sum_{j<i}\frac{1-g_i}{1-g_j}x_i\nabla_j(W_i)\theta_j.
\end{multline*}
Note that 
\[
\sum_{j<i}x_j\nabla_j(W_i)\theta_i=\sum_{j<i}(l_j(W_i)-l_{j+1}(W_i))\theta_i=\sum_{i}(W_i-l_{i}(W_i))\theta_i
\]
where we have used the fact that $l_{j+1}^{x\to(x,x^g)}=l_{j'}^{x\to(x,x^g)}$ for any pair $j<j'$ of {\it consecutive} elements of $I_g$.
Thus,
\begin{eqnarray}\label{temp}
\sum_{\alpha}(1-g_\alpha)x_\alpha\,\partial_{\theta_\alpha}({\rm H}_{{W,g}})=
\sum_{i}W_i\theta_i-\sum_{i}l_{i}(W_i)\theta_i-\sum_{j<i}\frac{1-g_i}{1-g_j}x_i\nabla_j(W_i)\theta_j.
\end{eqnarray}
Furthermore,
\begin{equation}\label{nablax}
l_j(x_i)=\begin{cases} x_i\quad j\leq i \\ 0\quad \text{otherwise}\end{cases}
\end{equation}
and by (\ref{leib})
\begin{multline*}
\sum_{j<i}\frac{1-g_i}{1-g_j}x_i\nabla_j(W_i)\theta_j=\sum_{j<i}\frac{1-g_i}{1-g_j}l_{j+1}(x_i)\nabla_j(W_i)\theta_j=\\
=\sum_{j<i}\frac{1-g_i}{1-g_j}\nabla_j(x_iW_i)\theta_j-\sum_{j<i}\frac{1-g_i}{1-g_j}\nabla_j(x_i)l_{j}(W_i)\theta_j=\\
=\sum_{j<i}\frac{1-g_i}{1-g_j}\nabla_j(x_iW_i)\theta_j =-\sum_{j\geq i}\frac{1-g_i}{1-g_j}\nabla_j(x_iW_i)\theta_j.
\end{multline*}
where the last equality follows from $\sum_{i}(1-g_i)x_iW_i=0$. Because of (\ref{nablax}) $\nabla_j(x_iW_i)=0$ for $j>i$. Using this fact, together with (\ref{leib}) and (\ref{nablax}), we obtain:
\begin{multline*}
\sum_{j<i}\frac{1-g_i}{1-g_j}x_i\nabla_j(W_i)\theta_j=-\sum_{i}\nabla_i(x_iW_i)\theta_i=\\=-\sum_{i}\left(\nabla_i(x_i)l_i(W_i)+l_{i+1}(x_i)\nabla_i(W_i)\right)\theta_i=-\sum_{i}l_i(W_i)\theta_i.
\end{multline*}
Substituting this result into (\ref{temp}) we get
$
\sum_{\alpha}(1-g_\alpha)x_\alpha\,\partial_{\theta_\alpha}({\rm H}_{{W,g}})=\sum_{i}W_i\theta_i.
$

\end{document}